\documentclass[11pt]{article}
\usepackage{amssymb,latexsym,amsmath, amsthm}
\usepackage{amsfonts,amssymb, latexsym, amsmath, amsthm,mathrsfs, verbatim,geometry,multicol}

\usepackage{amsfonts, graphics}
\usepackage{mathtools}
\usepackage{hyperref} 
\usepackage{csquotes} 
\usepackage{enumitem,moreenum} 
\usepackage{float}
\usepackage{dsfont,bbm} 
\usepackage{tikz}
\numberwithin{equation}{section}

\def\N{\mathbb N}

\hypersetup{
    colorlinks,
    citecolor=red,
    filecolor=green,
    linkcolor=blue,
    urlcolor=black
}

\def\N{\mathbb N}
\def\Z{\mathbb Z}

\def\D{\mathcal{D}}
\def\be{\begin{equation}}
\def\ee{\end{equation}}
\def\bea{\begin{eqnarray}}
\def\eea{\end{eqnarray}}
\def\beas{\begin{eqnarray*}}
\def\eeas{\end{eqnarray*}}

\def\pa{\partial }
\def\S{\mathcal S}

\def\D{{\mathcal D}}

\def\sign{{\mathrm{sign}}}

\def\pa{\partial}
\def\l{\lambda}

\def\al{\alpha}
\def\de{\delta}
\def\ga{\gamma}

\def\H{\mathcal H}
\def\lv{\left\vert}
\def\rv{\right\vert}
\def\Emin{E_{\text{min}}}
\def\k{\eps}

\def\l{\lambda}

\def\al{\alpha}
\def\de{\delta}
\def\ga{\gamma}

\let\th\relax
\newcommand{\th}{\theta}
\newcommand{\eps}{\varepsilon}
\newcommand\om{\omega}

\newcommand{\beq}{\begin{equation}}
\newcommand{\eeq}{\end{equation}}
\newcommand{\beqs}{\begin{equation*}}
\newcommand{\eeqs}{\end{equation*}}
\newcommand{\beqa}{\begin{equation}\begin{aligned}}
\newcommand{\eeqa}{\end{aligned}\end{equation}}
\newcommand{\beqas}{\begin{equation*}\begin{aligned}}
\newcommand{\eeqas}{\end{aligned}\end{equation*}}

\let\originalleft\left
\let\originalright\right
\renewcommand{\left}{\mathopen{}\mathclose\bgroup\originalleft}
\renewcommand{\right}{\aftergroup\egroup\originalright}

\newcommand*\diff{\mathop{}\!\mathrm{d}} 

\newtheorem{theorem}{Theorem}[section]

\newtheorem{corollary}[theorem]{Corollary}

\newtheorem{lemma}[theorem]{Lemma}

\newtheorem{remark}[theorem]{Remark}

\usepackage{tikz}
\usepackage{amsfonts, graphics}

\usetikzlibrary{calc}
\usetikzlibrary{arrows.meta}

\usepackage{tikz}
\usetikzlibrary{shapes,backgrounds,shadows,arrows} 
\usetikzlibrary{arrows.meta,calc,decorations.markings,math}
\definecolor{ubtgreen}{RGB}{3,138,94}

\def\bcr{\begin{color}{red}}
\def\bcb{\begin{color}{blue}}
\def\bcc{\begin{color}{violet}}
\def\ec{\end{color}}
\def\be{\begin{equation}}
\def\ee{\end{equation}}

\def\Tmin{T_{\text{min}}}

\def\Om{\Omega}

\def\S{\mathbb S}

\def\E{\mathcal E}

\definecolor{ubtgreen}{RGB}{3,138,94}

\def\che{\chi_{\text{ell}}}
\def\chh{\chi_{\text{hyp}}}
\def\chext{\chi_{\text{ext}}}
\def\Iell{I^{\text{ell}}}
\def\Ihyp{I_1^{\text{hyp}}}
\def\Ihypall{I^{\text{hyp}}}
\def\Iext{I^{\text{ext}}}
\def\pml{p_{m,\ell}}

\def\Eml{E_{\ell}^1}
\def\Emld{E_{\ell}^{1,\delta}}

\def\tE{\mathcal E}
\def\fin{f_{\text{in}}}
\def\rD{\overline{\text{ran}\,\D}}

\begin{document}

\title{On absence of embedded eigenvalues and stability of BGK waves}
\author{M.~Had\v zi\'c\thanks{University College London, UK. Email: m.hadzic@ucl.ac.uk}, M. Moreno\thanks{University College London, UK. Email: matias.bustamante.23@ucl.ac.uk}}

\maketitle

\begin{abstract}
We consider space-periodic and inhomogeneous steady states of the one-dimensional electrostatic Vlasov-Poisson system, known as the Bernstein-Greene-Kruskal (BGK) waves.
We prove that there exists a large class of fixed background ion densities and spatial periods, so that the corresponding linearised operator around the associated BGK-equilibria has no embedded eigenvalues inside the essential spectrum. As a consequence we conclude a nonquantitative version of Landau damping around a subclass of such equilibria with monotone dependence on particle energy. 

The BGK equilibria under investigation feature trapped electrons  which lead to presence of both elliptic and hyperbolic critical points in the characteristic phase-space diagram. They also feature a small parameter, which roughly speaking governs the size of the trapped zone - also referred to as electron hole.  Our argument 
uses action-angle variables and a careful analysis of the associated period function. 
To exclude embedded eigenvalues we develop an energy-based approach which deals with resonant interactions between the energy (action)-space and the angle frequencies; their singular structure and summability properties are the key technical challenge. Our approach is robust and applicable to other spectral problems featuring elliptic and hyperbolic critical points.
\end{abstract}

\tableofcontents

\section{Introduction}

A classical kinetic model of an electrostatic plasma is the Vlasov-Poisson system describing the dynamics of electrons against a fixed ion background. 
This model allows for a large class of steady states whose stability is one of the central topics in plasma physics. In this paper we work with 
space-inhomogeneous spatially periodic equilibria also known as Bernstein-Greene-Kruskal (BGK) solutions~\cite{BGK1957} (also BGK waves or BGK modes). 

Very little is known rigorously about asymptotic stability of the BGK waves. A nonquantitative version of Landau damping has been shown for the so-called Boltzmannian steady state in a pioneering work of Despr\'es~\cite{DESP2019}, wherein the author considers space-inhomogeneous ion backrgounds.  If the ion background is constant-in-space, Guo and Lin~\cite{GuoLin} identified a class of BGK equilibria that are stable and possess no embedded eigenvalues inside the essential spectrum, which also suggests the validity of Landau damping around such steady states. According to numerical evidence, BGK equilibria can act as attractors for the dynamics near  homogeneous equilibria~\cite{GuoLin,DeHoll1991,DaAnDr2004,BrCaPe2000,Ma1997,HSc2000,MeDiRoSh1998}, which makes them particularly important. In particular Lin and Zeng~\cite{LinZeng2011} showed that in sufficiently ``low" regularity the BGK equilibria prevent the Landau damping around homogeneous equilibria. 
We emphasise that many periodic BGK equilibria are also rigorously known to be unstable to perturbations with periods that are multiples (larger than $1$) of the period of the wave~\cite{GuoStrauss1995,Lin2001,Lin2005,PaAl2014}.

From the point of view of spectral theory, BGK equilibria lead to presence of (severe) particle trapping and both elliptic and hyperbolic critical points are present in the characteristic system. This is a challenge for any kind of decay result, and indeed decay can only be formulated in the orthogonal complement of the nontrivial kernel of the linearised operator. 
Inspired by our earlier work on the problem of gravitational Landau damping~\cite{HRSS2023}, which sees the effects of elliptic (but not hyperbolic critical points), our aim is to show that the essential spectra around a large class of BGK waves with {\em nonconstant} background ion densities have no embedded eigenvalues. This will imply a nonquantiative version of Landau damping around a subset of such equilibria that feature monotone dependence on the particle energy. The use of action-angle formalism provides a natural frequency variable for the problem (the angle) and is the first crucial tool in both~\cite{HRSS2023} and this work, while the second essential ingredient in both~\cite{HRSS2023, DESP2019} as well as the current work is the use of a suitable small parameter $\eps$ which we explain next.

We consider the Vlasov-Poisson system problem in the $1+1$-dimensional setting 
where the unknowns are the electron density $F(t,x,v)$ and the associated electrostatic potential $\varphi(t,x)$. Here $x$ belongs to a periodic interval $[0,\eps]$, where the period $\eps>0$ is some positive number, and $v\in\mathbb R$. We further make the common assumption that 
the ion spatial density distribution $\tilde\rho_+$ is time-independent and acts as a static background. The resulting system of equations reads
\begin{align}
\pa_t F +v \pa_x F + \pa_x \varphi \pa_v F & = 0,  \ \ (x,v)\in \eps\mathbb T^1\times\mathbb R, \label{E:VLASOV}\\
-\pa_{xx}\varphi & = \tilde \rho_+(x) - \int F\diff v, \ \ x\in \eps\mathbb T^1.\label{E:POISSON}
\end{align}
where $\eps\mathbb T_1$ stands for $\mathbb R/(\eps\,\mathbb Z)$ - the 1-dimensional torus with period $\eps$.

We now rescale the VP-system by letting
\be\label{E:SCALINGFF}
F(t,x,v) = f(t, \frac x\eps, \frac v\eps), \ \ \varphi(t,x) = \eps^2 \varphi(t,\frac x\eps).
\ee
Letting $\rho_+(\frac x\eps)=\tilde\rho_+(x)$ and renaming the scaled phase-space coordinates $(\frac x\eps,\frac v\eps)$ into $(x,v)$ again, the VP-system~\eqref{E:VLASOV}--\eqref{E:POISSON} takes the form
\begin{align}
\pa_t f + v \pa_xf + \pa_x\varphi \pa_vf & = 0, \ \ (x,v)\in\mathbb T^1\times \mathbb R \label{E:VLASOV2}\\
-\pa_{xx}\varphi& = \rho_{+}(x) - \eps \int_{\mathbb R} f \diff  v, \ \ x\in\mathbb T^1. \label{E:POISSON2}
\end{align}

From now on we shall work with the system~\eqref{E:VLASOV2}--\eqref{E:POISSON2}. 
By undoing the scaling~\eqref{E:SCALINGFF} we can map solutions of~\eqref{E:VLASOV2}--\eqref{E:POISSON2} to the original VP-system. 
In our analysis $0<\eps\ll1$ will be small and thus the solutions of the original system~\eqref{E:VLASOV}--\eqref{E:POISSON} will be spatially supported on a torus of very small period. Scaling~\eqref{E:SCALINGFF} is however not the only one that brings the VP-system into the form~\eqref{E:VLASOV2}--\eqref{E:POISSON2}. There is in fact a two-parameter family of rescalings - scaling~\eqref{E:SCALINGFF} is the only one that preserves the $L^\infty$ norms of $\rho_+$ and $f$ simultaneously. For more details see Appendix~\ref{A:SCALING}.


\subsection{Steady states.}


We consider a class of steady states $(f_0,\varphi_0)$ of the VP-system~\eqref{E:VLASOV2}--\eqref{E:POISSON2}
of the form
\begin{align}\label{E:EDEF}
f_0(x,v) = 
\mu(\k^2 E(x,v)), \ \ E(x,v) = \frac12 v^2 - \varphi_0(x),
\end{align}
where the ansatz function $\mu$ is a microscopic equation of state and $0<\eps\ll1$ a small parameter.
Under mild regularity assumptions on $\mu$ the Vlasov equation~\eqref{E:VLASOV2} is clearly satisfied by $f_0$, and thus it remains to check that
$\varphi_0$ solves the Poisson equation
\begin{align}
-\varphi_0''(x) + \eps \int \mu(\eps^2\big(\frac12 v^2-\varphi_0(x)\big)) \diff v = \rho_+(x).
\end{align}

Instead of prescribing the reference ion background $\rho_+$, we {\em prescribe} the potential $\varphi_0:\mathbb T^1\to\mathbb R_{\ge0}$ (as done for example by Despr\'es~\cite{DESP2019}) so that  
\begin{align}\label{E:RHOPLUS}
\rho_+(x) : = -\varphi_0''(x) + \eps \int \mu(\eps^2\big(\frac12 v^2-\varphi_0(x)\big)) \diff v >0, \ \ x\in\mathbb T^1.
\end{align}

\begin{figure}
\begin{center}
\begin{tikzpicture}
\begin{scope}[scale=0.5, transform shape]

\coordinate [label=below:$0$] (A) at (-3,7){};
\coordinate [label=below:$h$] (B) at (-7,7){};
\coordinate [] (C) at (-7,11){};
\coordinate [label=below:$e$] (D) at (1,7){};
\coordinate [label=above:$\mu(e)$] (F) at (-5.9,11){};
\coordinate [label=below:$\Emin$] (N) at (-5,7){};

\draw[very thick] (-3,7)--(-8,7);
\draw[-{>},very thick] (A)--(1,7);
\draw[very thick] (A)--(-3,11);
\draw[dashed, thick] (B) -- (-7,13);
\draw (-6.8,14) .. controls +(0.8,-4) and +(-1,0.5) .. (-3,9);
\draw (-3,9) .. controls +(1,-0.5) and +(-1,0) .. (1,7.8);
\draw[dashed] (N) -- (-5,10.2);

\draw[fill=black] (-3,7) circle (2pt);
\draw[fill=black] (-7,7) circle (2pt);
\draw[fill=black] (-3,9) circle (2pt);
\draw[fill=black] (N) circle (2pt);
\draw[fill=black] (-5,10.2) circle (2pt);

\coordinate [label=below:$0$] (E) at (10,7){};
\coordinate [label=below:$h$] (J) at (6,7){};
\coordinate [] (G) at (8,11){};
\coordinate [label=below:$e$] (H) at (14,7){};
\coordinate [label=above:$\mu(e)$] (I) at (9.1,10.5){};
\coordinate [label=below:$\Emin$] (M) at (7,7){};

\draw[very thick] (E)--(5,7);
\draw[-{>},very thick] (E)--(H);
\draw[very thick] (E)--(10,11);
\draw[dashed, thick] (J) -- (6,13);
\draw (7,7.7) .. controls +(0.8,2.1) and +(-1,0) .. (10,10);
\draw (10,10) .. controls +(1,0) and +(-1,0) .. (14,7.8);
\draw[dashed] (7,7) -- (7,7.7);

\draw[fill=black] (10,7) circle (2pt);
\draw[fill=black] (6,7) circle (2pt);
\draw[fill=black] (10,10) circle (2pt);
\draw[fill=black] (7,7.7) circle (2pt);
\draw[fill=black] (7,7) circle (2pt);

\end{scope}
\end{tikzpicture}
\caption{Schematic depiction of the regular  microscopic equations of state $e\mapsto \mu(e)$. On the left we have a monotone profile, and on the right a profile with a single local maximum.}
\label{F:MU}
\end{center}
\end{figure}
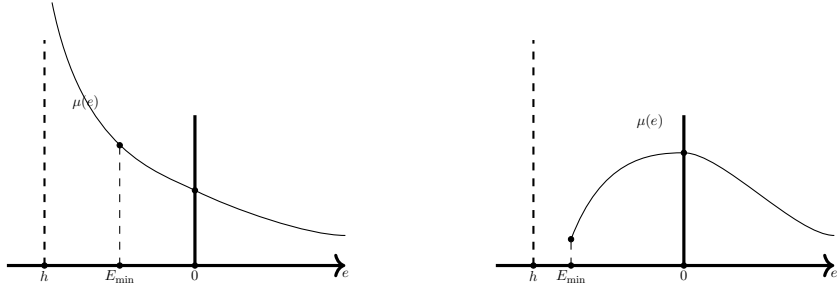
 
It will not be hard to check that the positivity of $\rho_+$ can be ensured for a large class of choices of $\mu$ and $\varphi_0$. 
We next provide precise assumptions on the microscopic equation of state function $\mu$. 
We assume that for some $h\in[-\infty,0)$, the profile $\mu:(h,\infty)\to\mathbb R_+$ has the following properties:
\begin{enumerate}
\item[($\mu1$)] ({\em Regularity})  We assume that 
\begin{align}
\mu\in C^{0,\frac12}((h,\infty),\mathbb R_+)  \ \ \text{(global H\"older continuity).} \label{E:MUREG} 
\end{align}
If we let 
\begin{align}
 \mu_-(e) = \mu(e)\chi_{(h,0]}(e), \ \ \mu_+(e) = \mu(e)\chi_{(0,\infty)}(e), 
\end{align}
then we assume
\begin{align}\label{E:MUMINUSPLUSREG}
\mu_-\in C^2((h,0],\mathbb R_+), \ \ \mu_+\in C^0([0,\infty),\mathbb R_+)\cap C^2((0,\infty),\mathbb R_+).
\end{align}
\item[($\mu2$)] ({\em Asymptotics at small positive energies})
There exists a $C>0$ such that
\begin{align}\label{E:MUBLOWUP}
|\mu'(e)| \le  Ce^{-\frac12}, \ \ |\mu''(e)|\le C e^{-\frac32}, \ \ 0<e<1.  
\end{align}
\item[($\mu3$)] ({\em Tail behaviour at high energies}) There exists a $C>0$ such that
\begin{align}
 |\mu'(e)|  \le \frac C{1+e}, \ \ |\mu''(e)|& \le \frac C{1+e^2}, \ \ e\ge1, \label{E:MUTAIL} \\
 \int_0^\infty |\mu'(e)|\sqrt e \diff e & <\infty. \label{E:COMP}
\end{align}
\item[($\mu4$)] 
\begin{align}
\text{$\mu_-$ is either strictly decreasing or strictly increasing, $\mu_+$ is strictly decreasing.} \label{E:MU4}
\end{align}
\end{enumerate}

\begin{remark}
By our assumptions, if $\mu$ is $C^2$ in a neighbourhood of $0$, then it is either monotone decreasing or it increases for $e<0$, has a positive maximum at $e=0$, and then decreases, see~Figure~\ref{F:MU}.
\end{remark}

\begin{remark}
Condition~\eqref{E:COMP} would be implied if we replaced the first bound in~\eqref{E:MUTAIL} by a slightly stronger decay bound
 $|\mu'(e)|  \le \frac C{1+e^{\frac{3}{2}+\de}}$ for some $\de>0$ and $e>1$.
\end{remark}

In the following we provide a few canonical examples of steady states covered in our work.

{\em Regular equations of state.}
If we choose $\mu\in C^2((h,\infty);\mathbb R)$ then conditions~\eqref{E:MUREG}--\eqref{E:MUBLOWUP} are trivially satisfied as the function is $C^2$ on all of its domain. We refer to such BGK waves as {\em regular}.
\begin{enumerate}
\item (Monotone regular BGK waves)
A particularly convenient class of regular steady states are the ones where $\mu$ is monotonically decreasing.
One family of examples of such solutions satisfying~\eqref{E:MUREG}--\eqref{E:MUTAIL}  studied first by Despr\'es~\cite{DESP2019} are the Boltzmannian steady states where 
\be\label{E:BOLTZMANNIAN}
\mu(e) = \exp(-\beta e), \ \ \beta>0, \ \ e\in\mathbb R.
\ee 
Exponential ansatz has been used in the physics literature (see e.g.~\cite{Hutch,HSc1979}), although it is not clear whether Maxwell-Boltzmann-like
distribution should play a distinguished role in the long-time asymptotics (see e.g. the discussion in~\cite{Hutch} and in the gravitational context~\cite{HRSS2023}). 

Another class of examples are algebraically decaying 
functions of the form 
\be\label{E:POLYMON}
\mu(e) = (e^{2k+1}-h^{2k+1})^{-1}, \ \ e>h,  \ k\in\mathbb N_0,
\ee
which we refer to as polytropes.

\item
(Non-monotone regular BGK-waves) A typical example of a non-monotone profile to which our results apply are even functions of the form
\begin{align}
\mu(e) = \frac{1}{1+e^{2k}}, \ \ e\in\mathbb R, \ \  k\in\mathbb N.
\end{align}
$\mu(\cdot)$ has a local maximum at the separatrix $e=0$ and is not monotone in $e$.
\end{enumerate}

{\em $C^{0,\frac12}$-equations of state (irregular examples).}
Here we are motivated by the
{\em Schamel-type} steady state, introduced by Schamel~\cite{HSc1979,HSc1986, Hutch} in the study of the electron trapping (electron ``holes") in stationary
plasmas. We use the same particle energy dependency, but work in a periodic interval and not the whole line (as done in the original works); nevertheless this model captures the essential features of trapping and the behaviour near the critical points. Note also that in~\cite{HSc1986} background ion density $\rho_+$ is assumed to be $1$, unlike this work. 
\begin{enumerate}
\item[3.]
We consider here a family of examples parametrised by the choice of the  lab frame velocity $\al\ge0$ and the negative of the inverse temperature $\beta>0$. We then define the Schamel-type distribution~\cite{HSc1986,Hutch} via
\begin{align}\label{E:MUALBE}
\mu_{\al,\beta}(e) =  \begin{cases} \exp(\beta e-\al^2) & e\in(-\infty,0], \\ \exp(-(\sqrt{e}+\al)^2) & e>0. \end{cases}
\end{align}
We may therefore write $\mu_-(e)= \exp{(\beta e-\al^2)}$ and $\mu_+(e)=\exp{(-(\sqrt e+\al)^2)}$ 
and note that $\mu$ is continuous. Moreover $\mu$ is $C^\infty$ on $\mathbb R\setminus\{0\}$ and $C^{0,\frac12}$ at $e=0$ when $\al\neq0$ (and if $\al=0$, then it is $C^{0,1}$ at $e=0$). In particular ansatz~\eqref{E:MUALBE} leads to a solution of the Vlasov equation pointwise everywhere except along the separatrix - which is a closed 1-dimensional curve in the phase-space given implicitly by the formula $\frac12 v^2 = \varphi_0(x)$.  It is readily checked that $\rho_{\al,\beta}(x)=\int_{\mathbb R}\mu(e(x,v))\diff v$ is $C^3$ on $\mathbb T^1$ and the conditions~\eqref{E:MUBLOWUP} and~\eqref{E:MUTAIL} are met. Properties $(\mu1)$--$(\mu4)$ are easy to verify.
\end{enumerate}

\begin{remark}
The steady states considered in~\cite{HSc1986} are not even in $v$ and therefore contain two branches in the region $e>0$ depending on the sign of $v$. For a discussion of the multi-branched case we refer to Appendix~\ref{A:SCHAMEL}.
\end{remark}

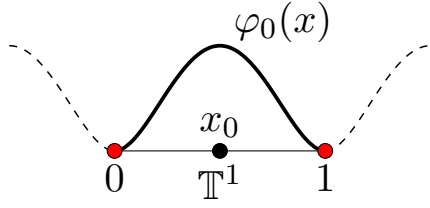
\begin{figure}[h!]
\begin{center}
		\begin{tikzpicture}
		\begin{scope}[scale=1.4, transform shape]
\draw[black, line width = 0.50mm]   plot[smooth,domain=0:2.0] (\x, {\x*\x*(2.0-\x)*(2.0-\x)});
\draw[black, dashed, line width = 0.2mm]   plot[smooth,domain=-1:0] (\x, {\x*\x*(\x+2)*(\x+2)});
\draw[black, dashed, line width = 0.2mm]   plot[smooth,domain=2:3] (\x, {(\x-2)*(\x-2)*(4-\x)*(4-\x)});

\coordinate[label=below:$0$] (A) at (0,0){};
\coordinate[label=below:$1$] (B) at (2,0){};
\coordinate[label=above:$\varphi_0(x)$] (C) at (1.6,.9){};
\coordinate[label=below:$\mathbb T^1$] (D) at (1,0){};
\coordinate[label=above:$x_0$] (E) at (1,0){};

\draw[black] (A) -- (B);
\draw[fill=red] (A) circle (2pt);
\draw[fill=red] (B) circle (2pt);
\draw[fill=black] (E) circle (2pt);

\end{scope}
		\end{tikzpicture}
	\end{center}
	\caption{Schematic depiction of the steady state potential $x\mapsto \varphi_0(x)$. See similar depiction in~\cite[Fig. 1]{Hutch} and~\cite{DESP2019}.}
	\label{F:POTENTIAL}
\end{figure}

\subsubsection{Assumptions on $\varphi_0$.} For a given $\mu$ satisfying~\eqref{E:MUREG}--\eqref{E:MUTAIL} above, we assume
\begin{enumerate}
\item[($\varphi 1$)] $\varphi_0$ is $C^3(\mathbb T^1,\mathbb R_{\ge0})$. Viewed as a function on the periodic interval $[0,1]$, $\varphi_0$ attains a unique maximum at some $x_0\in(0,1)$ and global minima at $x_0=0$ and $x_0=1$. Moreover 
\begin{align}
\varphi_0(0)&=\varphi_0(1)=0, \quad \varphi_0(x_0)>0, \quad  -\varphi_0(x_0) >h, \label{E:PHI_ASMP} 
\end{align}
where we remind the reader that $(h,\infty)$ is the domain of definition of $\mu$.
\item[($\varphi 2$)] $\varphi_0$ is strictly increasing in $(0,x_0)$ and strictly decreasing in $(x_0,0)$. We also assume that $\varphi_0$ is strcitly convex in a neighbourhood of $0$ and $1$ and strictly concave in a neighbourhood of $x_0$.
\item[($\varphi 3$)] (Existence of steady states) $\varphi_0''(x)<\int \mu(v^2/2)\diff v$ for all $x\in\mathbb T^1$. Together with~\eqref{E:MUREG}--\eqref{E:MUTAIL}, this assumption guarantees that for $\eps>0$ sufficiently small the ion reference density $\rho_+$ defined through~\eqref{E:RHOPLUS} is strictly positive on $\mathbb T^1$. 
\end{enumerate}

\begin{figure}
\begin{center}
\begin{tikzpicture}[>=stealth]
\def\wx{3}      
\def\wy{1.3}    
\def\wz{0.8}   

\draw[thick,red] (0,0) ellipse (\wx cm and \wy cm);

\draw[thick] (0,0) ellipse (2 cm and 0.8 cm);

\draw (-3,-3) -> (-3,3); 

\draw (3,-3) -> (3,3); 

\draw (-3,0) -> (3,0);

\coordinate[label=left:$0$] (C) at (-3,0.){};

\draw[fill=red] (C) circle (2pt);

\coordinate[label=right:$1$] (D) at (3,0.){};

\draw[fill=red] (D) circle (2pt);

\coordinate[label=below:$x_-(E)$] (A) at (-2,0.){};

\draw[fill=black] (A) circle (2pt);

\coordinate[label=below:$x_+(E)$] (B) at (2,0){};

\draw[fill=black] (B) circle (2pt);

\coordinate[label=above:$E>0\,(v>0)$] (E) at (0,1.9){};

\coordinate[label=below:$E>0\,(v<0)$] (F) at (0,-1.9){};

\coordinate[label=below:$E<0$] (G) at (0,1.3){};

\coordinate[label=below:$x_0$] (H) at (0.2,0){};

\draw[fill=black] (H) circle (2pt);

\draw[thick] (-3,2.3) .. controls +(1,0.7) and +(-1,0.7) .. (3,2.3);

\draw[thick] (-3,-2.3) .. controls +(1,-0.7) and +(-1,-0.7) .. (3,-2.3);

\end{tikzpicture}
\end{center}
\caption{Red trajectory is the separatrix corresponding to $E=0$. This is a prototypical phase-space diagram associated with the presence of
both elliptic and hyperbolic critical points considered by Despr\'es~\cite[Fig. 1]{DESP2019}, Guo-Lin~\cite{GuoLin}, Hutchinson~\cite[Fig. 2]{Hutch}, Schamel~\cite{HSc1979} and many others.\label{F:PS}}
\end{figure}
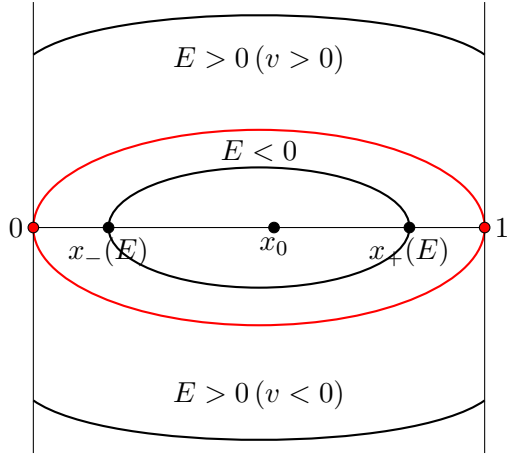

The key requirement in our analysis is that the period of trapped particle trajectories associated with the characteristic system
\begin{align}\label{E:CS}
\dot x = v, \ \ \dot v = \varphi_0'(x).
\end{align}
is monotone as a function of particle energy $E$. Concretely, since $E(x,v)$ is conserved along the flow of~\eqref{E:CS}, to each $E\in(\Emin,0)$ we associate the period function
\begin{align}
T(E) = 2 \int_{x_-(E)}^{x_+(E)} \frac{1}{\sqrt{2(E+\varphi_0(y))}}\diff y, \label{E:PERIODDEF}
\end{align} 
where $0<x_-(E)<x_+(E)<1$ are the turning points of the flow given by 
\begin{align}
\label{E:TURPOINTS}
\varphi_0(x_{\pm}(E))=-E, \ \ E\in (\Emin,0),
\end{align}
see Figure~\ref{F:PS}. Here 
\begin{align}
\Emin:=-\varphi_0(x_0) \label{E:EMINDEF}
\end{align} 
is the minimal energy level attained by the solutions of~\eqref{E:CS}. 
It can be shown that $T\in C^3([\Emin,0))$. Our final requirement is that
\begin{enumerate}
\item[($\varphi 4$)] (Monotonicity of the period function) 
\begin{align}\label{E:PERIODMON}
T'(E)\neq 0, \ \ E\in [\Emin,0).
\end{align}
We note that $T'(\Emin)>0$~\cite[Lemma A.12]{HRSS2023}, and therefore~\eqref{E:PERIODMON} implies $T'>0$ on $[\Emin,0)$.
\end{enumerate}

A simple way to guarantee condition~\eqref{E:PERIODMON} is to make a somewhat stronger assumption:
\begin{enumerate}
\item[($\varphi 4^\ast$)]
$(0,1)\ni x\mapsto \varphi_0'(x)^2+2(\varphi_0(x_0)-\varphi_0(x))\varphi_0''(x)$ is positive on $(0,1)\setminus\{x_0\}$. This is well-known to imply the monotonicity of $T$~\cite{HRS2021,ChWa1986,Ch1987}.
This assumption is also made in~\cite{DESP2019} to guarantee the monotonicity of the period.
\end{enumerate}

\begin{remark}
We note that similar assumptions are made in~\cite{DESP2019} (and in particular assumptions $(\varphi 3)$ and $(\varphi 4^\ast)$ are the same). 
\end{remark}

\begin{remark}[Monotonicity of the period]
It is well-known that monotonicity of the period plays an important role in description of phase-mixing properties of self-gravitating kinetic systems near trapping equilibria~\cite{LB1962, RiSa2020,HRS2021, HRSS2023,ChLu2022, MoRiVa2022, HRSS2024, ChLu2024}. 
In this work, by {\em choice} of the steady electrostatic potential $\varphi_0$ we can guarantee condition $(\varphi4)$ by for example demanding $(\varphi 4^\ast)$.
However the problem of proving the monotonicity of the period in the gravitational problem is challenging - this is an open question for a large class of 
important steady states of the (3d) gravitational VP-system~\cite{HRS2021,Straub2024}.
\end{remark}

\begin{remark}\label{R:INHOM}
As pointed out in~\cite{GuoLin}, if we assume that $\rho_+(x)\equiv c$ for some constant $c>0$, upon differentiating~\eqref{E:RHOPLUS} with respect to $x$ and then integrating against $\pa_x\varphi$ over $\mathbb T^1$ we obtain  the identity $\int_{\mathbb T^1} |\pa_x\varphi_0|^2\diff x -\eps \int_{\mathbb T^1\times\mathbb R}\mu'(e)|\pa_x\varphi_0|^2\diff (x,v)=0$. For monotone equlibria with $\mu'<0$ this is clearly impossible and implies that the background ion density is necessarily space-inhomogeneous.
\end{remark}

\begin{remark}[Ion background $\rho_+$ can be arbitrarily close to a constant in $L^\infty$]
We can choose $\varphi_0$ satisfying ($\varphi1$)--($\varphi4$) so that 
$\left\|\rho_+- 1\right\|_{L^\infty}\ll 1$. For instance choose the steady profile $\mu$ so that $\int\mu(v^2/2)\diff v=1$ and let $\varphi_0:=\de\tilde{\varphi_0}$, $\de>0$, where $\tilde{\varphi_0}$ satisfies the assumptions  ($\varphi1$)--($\varphi4$) with $\|\tilde{\varphi_0}\|_{C^2(\mathbb T^1)}<1$. Then clearly $\varphi_0$ also satisfies the same assumptions and therefore 
\begin{align}
	|\rho_+(x)-1|&\leq\de\|\tilde{\varphi_0}''\|_{\infty}+\eps \left|\int_{\mathbb R}\mu\left(\eps^2\big(\frac{v^2}{2}-\de\tilde{\varphi}_0(x)\big)\right)\diff v -1\right|\\
&\leq \de\|\tilde{\varphi_0}''\|_{\infty}+\left|\int_{\mathbb R}\mu\left(\frac{\tilde v^2}{2}´-\k^2\varphi_0(x)\right)\diff \tilde v -1\right|
\end{align} 
It follows that for $0<\eps\ll1$ sufficiently small, $\|\rho_+-1\|_{L^\infty}\le 2\de$ for any $\de>0$.
\end{remark}


\subsection{Linearisation and action-angle coordinates}


The formal linearisation of the Vlasov-Poisson system takes the form
\begin{align}
\pa_tf+\D f + \k^2 \mu'(\k^2 E)v\pa_x \varphi_f & = 0 \label{E:VLASOVLIN}\\
\pa_{xx}\varphi_f & =\eps  \int f\diff v, \label{E:POISSONLIN}
\end{align}
where the antisymmetric operator $\D$ is given by 
\begin{align}\label{E:DDEF}
\D = v\pa_x +\varphi_0'(x)\pa_v.
\end{align}

Letting $\Om:=\mathbb T^1\times \mathbb R$, the natural Hilbert space for our analysis is given by
\begin{align}\label{E:HILBERTSPACE}
\mathcal H : = \big\{f:\Omega\to\mathbb R \text{ meas.} \,\big| \iint_{\Om}\frac{f(x,v)^2}{|\mu'(\eps^2 E(x,v))|}\diff (x,v)<\infty, \ \ \iint f(x,v)\diff (x,v)=0\big\},
\end{align}
and it is standard to check that $\D: D(\D)\subset \H\to\H$ is indeed skew-adjoint with the domain of definition $D(\D)$ given by all $f\in\H$ so that $\D f\in \H$ weakly.
The linearised operator $\mathscr L: D(\mathscr L)=D(\D)\subset \H \to \H$ is given by
\begin{align}
\mathscr Lf & = \D f + \eps^2 \mu'(\eps^2E) v\pa_x \varphi_f = \D (f+ \eps^2 \mu'(\eps^2 E)\varphi_f).
 \label{E:DL}
\end{align}
Lemma~\ref{L:FORCEBOUND} shows that $\mathscr L$ is well-defined, in fact the operator $\H \ni f\mapsto \eps^2 \mu'(\eps^2E) v\pa_x \varphi_f \in \H$ is bounded.

The core tool in our work are the action-angle variables, 
that among other things  allow us to easily describe the essential spectrum of $\mathscr L$. 
To define them we first observe that the stationary points of the characteristic system~\eqref{E:CS} are given by 
\begin{enumerate}
\item  $(x_0,0)$; this corresponds to the minimal energy level $\Emin=-\varphi_0(x_0)$ and we refer to it as the {\em elliptic} critical point,
\item  $(0,0)=(1,0)$; this point has energy zero $E=0$ and we refer to it as the {\em hyperbolic} critical point.
\end{enumerate}

%

We introduce the action support $I$ given by 
\begin{align}
I : = I_1 \cup I_2, \ \ \ 
I_1 : = (\Emin,0),\ \ \
I_2  : = (0,\infty).
\end{align}
Particles with energy $E\in I_1$ are trapped and lead to periodic orbits. As such they present a major obstacle to any sort of decay. Particles with positive energies, i.e. $E>0$ are outside the separatrix and wind around the torus with either the positive velocity ($v>0$) or negative velocity ($v<0$).

{\em The trapped region $E\in I_1$.}
The period of this motion is given by the standard formula~\eqref{E:PERIODDEF}.
For $E=\frac12 v^2-\varphi_0(x)\in I_1$ we define
the angle by
\begin{align}\label{E:THETADEF}
\theta(x,v) = \frac1{T(E)}\int_{x_-(E)}^x \frac{1}{\sqrt{2(E+\varphi_0(y))}}\diff y \in[0,\frac12], \ \ v\ge0,
\end{align}
where we recall~\eqref{E:TURPOINTS}.
For $v<0$ and $E\in I_1$ we extend the above definition by letting $\th(x,v) := 1-\th(x,-v)$, so that the trapped region
$\{x,v\} = \{\Emin<E(x,v)<0\}$ is parametrised by the cylinder $\mathbb S^1\times I_1$.

\begin{remark}[Trapping]
The particle trapping phenomenon is typical of BGK waves and is tied to the presence of so-called electron/ion holes. It plays an important role in the description of global dynamics of plasmas~\cite{HSc1979,HSc1986,DeHoll1991,LuSch2005,Hutch,CGNA2023}. In the gravitational case the trapping appears very naturally as space-homogeneous steady states are inadequate for the description of isolated galaxies - see the pioneering discussion in~\cite{LB1962,LB1967} as well as for example~\cite{Kalnajs1971,Kalnajs1977,BaOlYa,HRS2021,HRSS2023}.
\end{remark}

{\em The outside region $E\in I_2$.}
In the range $E\in I_2$ the period function is given by the formula
\begin{align}
T(E) = \int_{0}^{1} \frac{1}{\sqrt{2(E+\varphi_0(y))}}\diff y.
\end{align}
For any $x\in [0,1]$ the angle $\th(x,v)$ is given by 
\begin{align}
\th(x,v) = \frac1{T(E)}\int_{0}^x \frac{1}{\sqrt{2(E+\varphi_0(y))}}\diff y \in [0,1].
\end{align}
Thus the particles with strictly positive energies $\{E(x,v)>0\}$ are parametrised by the cylinder $\mathbb S^1\times I_2$.
We shall also use the notation
\begin{align}\label{E:TURN_POINTS_EXT}
x_-(E):=0,\ x_+(E):=1\text{ for any }E\in I_2.
\end{align}
{\em The map $(\th,E)\mapsto (x,v)$.}
The change of variables $(x,v)\mapsto(\th,E)\in \mathbb S^1\times I$ is 1-1 on the set $\Om\setminus\{(x_0,0), (0,0),(1,0)\}$.
To explicitly invert the map we  let 
\begin{align}\label{E:XVTHETA}
x(\theta,E):=X(\theta T(E),E), \ \ v(\theta,E):=V(\theta T(E),E), \ \ \th\in[0,1],
\end{align}
where $s\mapsto (X(s,E),V(s,E))$ is the solution to~\eqref{E:CS} with data
\begin{align}
X(0,E)=x_{-}(E),\ \ V(0,E)=0, \ \ E\in I.
\end{align}
Our fundamental new variables are thus given by $(\theta,E)$. The change of variables $(\th,E)\mapsto (x,v)$ is an invertible differentiable map  for all $(\th,E)\in\mathbb S^1\times I$
and the volume element transforms according to $\diff x \diff v=T(E)\diff \th \diff E$. This transformation does not preserve the phase-space volume, but nevertheless we continue to refer to $(\th,E)$ as {\em action-angle} variables.  It will often be convenient to use the frequency function
\begin{align}\label{E:FREQUENCYDEF}
\om(E): = \frac1{T(E)}, \ \ E\in I.
\end{align}

{\em The transport operator.}
It is well known~\cite{GuoLin,HRSS2023} that the transport operator 
\begin{align}
\D : = v\pa_x +\varphi_0'(x)\pa_v
\end{align}
transforms as
\begin{align}\label{E:DPOS}
\D g= \om(E) \pa_\theta g, \ \ \Emin<E(x,v)<0,
\end{align}
and 
\begin{align}\label{E:DNEG}
\D g= \sign (v) \om(E) \pa_\theta g, \ \ E(x,v)>0.
\end{align}
The sign distinction in~\eqref{E:DNEG} captures the different winding orientation
that particles outside the separatrix have depending whether they are located ``above" the separatrix ($v>0$)
or ``below" ($v<0$). The following lemma, the proof of which is found in Appendix~\ref{A:SPECTRAL}, states some well-known
spectral properties of the operator $\mathscr L$.


\begin{lemma}\label{L:ESSENTIAL}
Let $(f_0,\varphi_0)$ be a steady state of the VP-system~\eqref{E:VLASOV2}--\eqref{E:POISSON2}, so that 
$(f_0,\varphi_0)$ is of the form~\eqref{E:EDEF}, 
the microscopic equation of state $\mu(\cdot)$ satisfies properties~\eqref{E:MUREG}--\eqref{E:MU4}
and the associated 
electrostatic potential $\varphi_0$ satisfies properties $(\varphi1)$--$(\varphi4)$.
Consider the linearised operator $\mathscr L:D(\D)\to \H$ and the pure transport operator $\D:D(\D)\to\H$. Then
\[
\sigma_{\text{ess}}(\mathscr L) = \sigma_{\text{ess}}(\D) = i \mathbb R,
\]
where $\sigma_{\text{ess}}(A)$ is the essential spectrum of $A$.  The kernel of $\D$ is infinite-dimensional and is given by 
\be\label{E:KERD0}
\text{ker}\,\D = \big\{f\in \H\,\big|\, f = \phi(E) \ \text{a.e.}\big\}. 
\ee
Moreover the orthogonal complement of the kernel is given by
\begin{align}\label{E:KERD}
(\text{ker}\,\D)^\perp = \rD= \big\{f\in\H \, \big| \, \int_{\S^1}f(\th,E)\diff \th = 0 \ \ \text{a.e.} \ E\in I\big\}. 
\end{align}
\end{lemma}


\subsection{Main results}


Our main theorem states that there are no embedded eigenvalues inside the essential spectrum of $\mathscr L$.


\begin{theorem}[Absence of embedded eigenvalues]\label{T:MAIN}
Let $(f_0,\varphi_0)$ be a steady state of the VP-system~\eqref{E:VLASOV2}--\eqref{E:POISSON2}, so that 
$(f_0,\varphi_0)$ is of the form~\eqref{E:EDEF}, 
the microscopic equation of state $\mu(\cdot)$ satisfies properties~\eqref{E:MUREG}--\eqref{E:MU4}
and the associated 
electrostatic potential $\varphi_0$ satisfies properties $(\varphi1)$--$(\varphi4)$ (or $(\varphi1)$--$(\varphi3)$ and $(\varphi4^*)$).
Then there exists an $\eps_0>0$ such that for any $0<\eps<\eps_0$ the  operator~\eqref{E:DL} $\mathscr L:D(\mathscr L)\subset \mathcal H\to\mathcal H$ has no non-zero embedded eigenvalues.
\end{theorem}


\begin{remark}
See Appendix~\ref{A:SCHAMEL} for the analogous result for ``multi-branched" steady state functions $\mu$.
\end{remark}

We next specialise to {\em monotonically} decreasing regular equilibria, i.e. we assume $\mu'<0$ and $\mu\in C^2((h,\infty),\mathbb R_+)$. Examples are~\eqref{E:BOLTZMANNIAN} and~\eqref{E:POLYMON}.
We introduce the operator 
\be\label{E:LDEF}
L f = f + \eps^2 \mu'(\eps^2 E) (\varphi_f-\int\varphi_f\diff x).
\ee
It is a bounded linear operator on $\H$, as shown in Lemma~\ref{L:FORCEBOUND}.
Since $\mu'<0$ 
it
is easy to check that $L: \H\to \H$ is in fact self-adjoint\footnote{Without the monotonicity assumption on $\mu$, $L$ is in general not symmetric.}
  with the associated quadratic form
\be\label{E:FLF}
(f,Lf)_{\H}  = \iint_{\Om}\frac{f^2}{|\mu'(\eps^2 E)|}\diff(x,v) + \eps \int_{\mathbb T}|\pa_x\varphi_f|^2 \diff x.
\ee
By Lemma~\ref{L:FORCEBOUND}
it follows that $\|L-\text{Id}\|_\H \le C\eps$ for some $C>0$, and therefore the operator $L$ is invertible for $0<\eps\ll1$ sufficiently small.
The linearised operator~\eqref{E:DL} can be decomposed as
 \begin{align}
\mathscr L & = \D L. 
\end{align}
It is then easy to check that $\mathscr L $ is  skew-adjoint on the Hilbert space $\H$ equipped with the inner product
\begin{align}\label{E:INNERPRODL}
\left(f,g\right)_{L}: = (f, Lg)_{\H}.
\end{align}
Due to the bijectivity of $L:\H\to\H$, the kernel of $\mathscr L$ is in $1-1$ relation with $\text{ker}\D$ and is given by
\be\label{E:KERNEL1}
\text{ker}(\mathscr L) = L^{-1}\left(\text{ker}\D\right).
\ee
The only meaningful discussion of Landau damping happens if we restrict the dynamics to the $(\cdot,\cdot)_L$-orthogonal complement of $\text{ker}\,\mathscr L $, which is given by $\overline{\text{ran}\,\D}$.\footnote{Note that $\overline{\text{ran}\,\D}$ is an invariant subspace for the linearised dynamics.} 
We therefore study the spectrum of the operator $\mathscr L\big|_{\rD} :\rD\cap D(\D) \to \rD$.

As an application of the RAGE theorem~\cite{CyFrKiSi}, inspired by a related argument for 2D fluids~\cite{LinZeng2022} and the gravitational VP-system~\cite{HRSS2023}, we use Theorem~\ref{T:MAIN} to show that at least 
in a time-averaged sense, the $L^2$-norm of the electrostatic force decays to $0$ for initial data in the orthogonal complement
of the kernel. This should be viewed as a statement about the validity of (nonquantitative) Landau damping around regular monotone equilibria in question.


\begin{corollary}\label{C:LANDAU}
Let $(f_0,\varphi_0)$ be a steady state of the VP-system~\eqref{E:VLASOV2}--\eqref{E:POISSON2} of the form~\eqref{E:EDEF}, so that $\mu$ is $C^2$ and strictly decreasing (in particular $(\mu1)$, $(\mu2)$, and $(\mu4)$ hold automatically). Assume additionally that $(\mu 3)$ holds. 
Then there exists an $\eps_0>0$ such that for any $0<\eps<\eps_0$ and any $\fin\in \rD\cap D(\D)\subset\H$
\begin{align}
\lim_{T\to\infty}\frac1T\int_0^T\|\pa_x\varphi_f(t,\cdot)\|_{L^2}\,dt = 0.
\end{align}
where $f(t,\cdot,\cdot)$ is the solution of the linearised flow~\eqref{E:VLASOVLIN}--\eqref{E:POISSONLIN}.  
\end{corollary}

Proofs of Theorem~\ref{T:MAIN} and Corollary~\ref{C:LANDAU} are presented in Section~\ref{S:PROOFS}.

Our setup is close to the work of Despr\'es~\cite{DESP2019} where the author analysed Landau damping around the particular Boltzmannian background~\eqref{E:BOLTZMANNIAN} with a  small parameter, which is equivalent to our small parameter $0<\eps\ll1$ upon a suitable rescaling.
For such a background, result of~\cite{DESP2019} is strictly stronger than Corollary~\ref{C:LANDAU} as it states pointwise (but also nonquantitative)
decay  $\lim_{t\to\infty}\|\pa_x\varphi_f(t,\cdot)\|_{L^\infty}=0$. The methods of~\cite{DESP2019} are however different from our approach; the author constructs a basis of generalised eigenvectors for the linearised dynamics using ideas from scattering theory, by analysing a Lippmann-Schwinger type equation (for a discussion of the scattering theory approach see also~\cite{Despres2021}). As we shall see below our method is strongly rooted in the use of action-angle formalism and simple energy arguments. 

In both the present work and~\cite{DESP2019} the background ion density $\rho_+$ is space-inhomogeneous. This follows from the way we choose the steady state: we first prescribe the potential and then use it to define the background ion density. In fact by Remark~\ref{R:INHOM} inhomogeneity is necessary if the $\mu$-profile is monotone like in Corollary~\ref{C:LANDAU}.

In contrast to us and~\cite{DESP2019}, Guo and Lin~\cite{GuoLin} exhibit a class of spectrally stable (space-inhomogeneous) BGK waves with constant ion backgrounds ($\rho_+=1$), which contain no embedded eigenvalues in their spectrum. This result is hard as the homogeneity of $\rho_+$ makes the analysis complicated - in this regime the steady state potential $\varphi_0$ is not explicitly prescribed, but instead solves a second order nonlinear ODE. On the other hand, the stable Guo-Lin BGK waves by construction have  constant phase-space density throughout the trapped region and slightly outside it. In particular, their stability analysis ``feels" the trapping only through the pure transport operator $\D$ and not through the nonlocal operator $\eps^2\mu'(\eps^2 E)v\pa_x$ in~\eqref{E:VLASOVLIN}, which then effectively restricts the analysis to the region outside and strictly away from trapping.

{\em Space-homogeneous nontrivial equilibria.}
The stability of space-homogeneous
stationary solutions is a classical topic that goes back to Landau~\cite{Landau1946} who showed that small perturbations
around spectrally stable spatially periodic homogeneous equilibria damp away for the linearised problem. The associated decay of macroscopic
observables such as the spatial density and electrostatic potential is due to the physical phenomenon of phase-mixing and is referred to as Landau damping.
In this case there is a vast literature on the question of stability. For spectrally stable steady states satisfying the so-called Penrose stability condition, Mouhot and Villani~\cite{MoVi2011} proved nonlinear Landau damping in Gevrey spaces, see also~\cite{BeMaMo2016} and~\cite{GrNgRo2020} for different and shorter proofs.  A comparison of different proofs and general overview of Landau damping related literature can be found in~\cite{Be2022}, see also~\cite{IoPaWaWi2024}. The former results apply to spatially periodic settings where phase-mixing is the dominant stabilisation mechanism. On the other hand, the whole domain stability problem is still under investigation. The linearised stability analysis was given in~\cite{HKNgRo2022,BeMaMo2022}, see also~\cite{GlSc94}; for a discussion of the linear decay and the so-called survival threshold we refer to the recent work~\cite{Ng2023}. For a nonlinear stability result around the Poisson equilibrium on $\mathbb R^3$ see~\cite{IoPaWaWi2022} and also~\cite{IoPaWaWi} for the discussion of linear decay around general equilibria. For the so-called screened case we refer to~\cite{BeMaMo2018,NgTT,HKNgRo2021,HuNgXu2,IaRoWi} and references therein. For results and more references in the case of massless electrons we refer the reader to~\cite{GaIa,HuNgXu1}. 


\subsection{Proof strategy}


The proof of Theorem~\ref{T:MAIN} is based on a contradiction argument. We assume that $g\in D(\D)\subset \H$ is an eigenfunction associated to an eigenvalue $\l\in i\mathbb R$, $\l\neq0$. This implies that $g(\cdot,E)\in L^2(\mathbb S^1)$ for almost every $E$. We then pass to Fourier series in the angle-variable so that 
\begin{align*}
g(\th,E) = \sum_{\ell\in\mathbb Z_\ast} \widehat g(\ell,E)e^{2\pi i \ell\th}, \ \ \varphi_g(\th,E) = \sum_{\ell\in\mathbb Z_\ast} \widehat{\varphi_g}(\ell,E)e^{2\pi i \ell\th}, \ \ E\in I,
\end{align*}
where $\Z_\ast=\Z\setminus\{0\}$.
Note that we abuse notation here and continue to use the same letter $g$ to write $g(x,v)=g(\th,E)$.
The key idea in our analysis is to express the electrostatic energy, by virtue of the Plancharel theorem, as a summation over all angular frequencies of integrals over the 
action range $I$. 
For simplicity we assume that $\mu$ is $C^2$ and $\mu'<0$, just to isolate the key ideas of the argument more easily.
The eigenvalue equation reads
\begin{align}\label{E:EVEQN}
\D (g - \k^2 |\mu'(\k^2 E)| \varphi_g) = \l g,
\end{align}
and we use it to express the Fourier coefficients $\widehat{g}(\ell,E)$ in terms of $\widehat{\varphi_g}(\ell,E)$. In the trapped region $E<0$ equation~\eqref{E:EVEQN}
gives the relation (Lemma~\ref{L:FOURIERFORMULA})
\begin{align}\label{E:GFORMULAINTRO}
\widehat g(\ell,E) = \varepsilon^2\frac{|\mu'(\k^2E)|\widehat{\varphi_g}(\ell,E)}{1-\frac{T(E)}{\ell q}}, \ \ (\ell,E)\in\Z_\ast\times I_1, \ \ \ell q\neq T(E),
\end{align}
where 
\[
q : = \frac{2\pi i }{\l} \in \mathbb R.
\] 
Unlike the classical
setting of space-homogeneous steady states~\cite{Landau1946,MoVi2011,BeMaMo2016,GrNgRo2020} the map $g\mapsto\varphi_g$ does not separate frequencies, i.e. $\widehat{\varphi_g}(\ell,E)$
depends on all the Fourier modes $\widehat{g}(\ell',E)$, $\ell'\in\mathbb Z_\ast$, and therefore the dispersion relation~\eqref{E:GFORMULAINTRO} does not reduce to an equation on $\widehat{g}(\ell,E)$ alone. Instead we resort to an energy-type argument which is global in nature and allows us to relate the information on individual Fourier coefficients to the $L^2$-norm of the electrostatic force. In other words we observe that
\begin{align}
\int_{\mathbb T^1}|\pa_x\varphi_g|^2 \diff x & = - \int \varphi_g \pa_{xx}\varphi_g\diff x = - \eps \iint g(x,v) \varphi_g(x) \diff(x,v) \notag\\
& = -\eps \sum_{\ell\in\Z_\ast} \int_{I} \widehat{g}(\ell,E) \overline{\widehat{\varphi_g}(\ell,E)} T(E)\diff E \notag\\
& = \eps^3  \sum_{\ell\in\Z_\ast} \int_{I} \frac{T(E)|\mu'(\k^2E)|}{\frac{T(E)}{\ell q}-1}|\widehat{\varphi}_g(\ell, E)|^2\diff E, \label{E:ENERGYINTRO}
\end{align}
where we have used~\eqref{E:GFORMULAINTRO} in the last line. Our goal is to show that 
\begin{align} 
\text{ right-hand side of~\eqref{E:ENERGYINTRO} }  \le C \eps^a\|\pa_x\varphi_g\|_{L^2(\mathbb T^1)}^2
\end{align}
for some  $a>0$. This then leads to contradiction for $0<\eps\ll1$ sufficiently small.
This is however easier said than done and there are broadly speaking three types of obstacles to proving this.

{\em 1. Small denominators.} The above strategy is predicated on the formula~\eqref{E:GFORMULAINTRO} which is valid only for frequency-energy pairs $(\ell,E)$ where $T(E)\neq \ell q$. Clearly, when $|\frac{T(E)}{\ell q}-1|\ll1$ the denominator in~\eqref{E:ENERGYINTRO} is small and this is a serious difficulty. We go around this by rewriting
\begin{align}
\frac1{\frac{T(E)}{\ell q}-1} =\frac{\ell q}{T'(E)} \frac{d}{dE} \log(\frac{T(E)}{\ell q}-1),
\end{align}
which is possible due to the fundamental period monotonicity assumption~\eqref{E:PERIODMON}; this idea was introduced in~\cite{HRSS2023} in the gravitational setting.
In doing so we integrate-by-parts in~\eqref{E:ENERGYINTRO} and replace the $\frac1{\frac{T(E)}{\ell q}-1} $-singularity by the milder $\log$-singularity, at the expense of 
$E$-derivatives hitting the electrostatic potential. We therefore need to bound terms of the form
\begin{align}
&\ell q \eps^3 \int _I  \log(\frac{T(E)}{\ell q}-1) \frac{d}{dE}\big( \frac{T(E)|\mu'(\k^2E)|}{T'(E)}|\widehat{\varphi}_g(\ell, E)|^2\big) \diff E \notag\\
& = 2\ell q \eps^3 \int _I  \log(\frac{T(E)}{\ell q}-1)  \frac{T(E)|\mu'(\k^2E)|}{T'(E)}\pa_E\widehat{\varphi}_g(\ell, E)\overline{\widehat{\varphi}_g(\ell, E)} \diff E\notag\\
& \ \ \ \ +\ell q \eps^3 \int _I  \log(\frac{T(E)}{\ell q}-1) \frac{d}{dE}\big( \frac{T(E)|\mu'(\k^2E)|}{T'(E)}\big) |\widehat{\varphi}_g(\ell, E)|^2\diff E =: \mathcal A+\mathcal B.
\end{align}

To control  $\mathcal A$ we must bound $E$-derivatives of $\varphi_g$ in terms of $\pa_x\varphi_g$, which naturally leads to the estimate 
\begin{align}
|\pa_E\varphi_g| \le |\pa_x\varphi_g| |\pa_E x(\th,E)|.\label{E:PARTIALPHI}
\end{align}
Due to the degeneracy of action-angle coordinates at the critical points, $|\pa_E x(\th,E)|$ blows up at both $E=\Emin$ and $E=0$. To have a chance of bounding $\mathcal A$ we need a precise rate at which $|\pa_E x(\th,E)|$ degenerates.

{\em 2. Elliptic/hyperbolic degeneracy.}
Near the energy minimum $\Emin$ we have the bound
\begin{align}
|\pa_E x(\th,E)| \lesssim \frac{1}{\sqrt{E-\Emin}}, \ \  (\th,E)\in \mathbb S^1\times I_1. \label{E:ELLDEGEN}
\end{align}
A related issue near the elliptic point has already been encountered in~\cite{HRSS2023} and we use similar strategy. We estimate $\mathcal A$ via Cauchy-Schwarz, but prior to that we desingularise $\pa_E\varphi_g$ by offloading some of the negative powers of $E-\Emin$ coming from~\eqref{E:ELLDEGEN} via~\eqref{E:PARTIALPHI} to the (very) mildly singular $\log$-term.
 
Unlike the elliptic singular point, near $E=0$ the period function and its derivatives blow up 
\begin{align}
& T(E)\sim |\log |E||, \ \ |T'(E)|\sim\frac{1}{|E|}, \ \ |T''(E)|\sim \frac1{|E|^2}, \ \ E \text{ near } \ 0, \label{E:HYPDEGENINTRO1}\\ 
& |\pa_E x(\th,E)| \lesssim \frac{1}{|E|}.  \label{E:HYPDEGENINTRO2}
\end{align}
We note that~\eqref{E:HYPDEGENINTRO2} is a priori {\em more} singular than~\eqref{E:ELLDEGEN}. However, the presence of various $T$-dependent quantities in $\mathcal A$ introduces cancellations near $E\sim0$, which interact with the $\frac1{|E|}$-singularity in a way that allows us to show integrability near $E=0$.

{\em 3. Summability in $\ell\in\Z_\ast$ and uniformity in $q$.}
To illustrate the above point we focus on the summability problem in the hyperbolic zone; one typical term that we need to bound is of the form
\begin{align}\label{E:TILDEA}
\tilde{\mathcal A}:=\sum_{\ell>0}\ell^2 q^2 \eps^3 \int_{I_\ell}  |\log(\frac{T(E)}{\ell q}-1)|^2  \frac{T(E)^2|\mu'(\k^2E)|}{T'(E)}|\widehat{\varphi}_g(\ell, E)|^2 \diff E,
\end{align}
which is one contribution coming from $\mathcal A$ upon using the Cauchy-Schwarz inequality. Here for any fixed frequency $\ell$, $I_\ell$ encodes energy levels that are nearly resonant in the sense that 
\[
|\frac{T(E)}{\ell q}-1|<\delta,
\]
for some control parameter $\delta>0$.
Since we understand the asymptotics of the period function very well (this is the content of Theorem~\ref{T:PERIODASYMP}), we have a very good control over $I_\ell$ in the hyperbolic zone $|E|\ll1$.
The idea is to estimate $|\widehat{\varphi_g}(\ell,E)|^2\le \frac1{4\pi\ell^2}\int_{\mathbb S^1}|\pa_\th \varphi_g(\th,E)|^2\diff \th$, plug it into~\eqref{E:TILDEA}, and change variables $\th\mapsto x$ to obtain
\begin{equation}
\tilde{\mathcal A} \le C \eps^3 q^2 \|\partial_x\varphi_g\|_{L^2}^2   \sum_{\ell>0}\int_{I_\ell}  |\log(\frac{T(E)}{\ell q}-1)|^2  \frac{T(E)^2|\mu'(\k^2E)|}{T'(E)} \diff E. \label{E:AELLONEXINTRO}
\end{equation}
A crucial consequence of Theorem~\ref{T:PERIODASYMP} is that in the hyperbolic zone we have sharp asymptotics from~\eqref{E:HYPDEGENINTRO1}
\begin{align}\label{E:KEYHYPINTRO}
\exp(T(E)) \lesssim |T'(E)| \lesssim  \exp(T(E)), \ \ 0<|E|\ll1.
\end{align}
Upon changing variables $\frac{T(E)}{\ell q} = t$ and carefully examining the interval $I_\ell$, estimate~\eqref{E:AELLONEXINTRO} leads to
\begin{align}
\tilde{\mathcal A} &\le C \eps^3  \|\partial_x\varphi_g\|_{L^2}^2 q^5\sum_{\ell>0}\ell^3  \int_{1}^{1+\delta}\left|\log\left(t-1\right)\right|^2 t^2 \, e^{-(2-\alpha) \ell q t}\diff t \notag\\
& \le C  \eps^3  \|\partial_x\varphi_g\|_{L^2}^2 q^5\sum_{\ell>0}\ell^3  e^{-2\ell q } \le C  \eps^3  \|\partial_x\varphi_g\|_{L^2}^2,
\end{align}
where the last bound can be inferred from the integral criterion and the exponential decay in $q\ell$ to obtain a uniform-in-$q$ bound.
Precise statements, the additional subtle treatment of the elliptic degeneracy $|E-\Emin|\ll1$ and the high energy degeneracy $E\gg 1$ are captured in key energy bound of Lemma~\ref{L:ENERGYBASIC} and the estimates of Lemmas~\ref{L:TRAPLEMMA},~\ref{L:EXTLEMMA}, and~\ref{L:NRLEMMA}.

\bigskip

{\bf Plan of the paper.}
In Section~\ref{S:PERIOD} we formulate the key lemmas that quantify the behaviour of the period function and its derivatives near the critical points. This is then followed up by 
Theorem~\ref{T:DXDE} which quantifies the singular nature of the map $E\mapsto x(\th,E)$ at the critical points. Section~\ref{S:ABSENCE} is dedicated to the proof of the main theorem. We start with a suitable subdivision of the frequency-action phase space that is necessary for tracking the small denominators - see Section~\ref{S:FS}. The contradiction argument based on an energy identity sketched above is carried out in Section~\ref{S:EE}. This is the most important part of the paper. Finally, Theorem~\ref{T:MAIN} and Corollary~\ref{C:LANDAU} are proved in Section~\ref{S:PROOFS}. Proof of the key asymptotic properties of the period function stated in Section~\ref{S:PERIOD} is contained in Appendix~\ref{A:PERIOD}. 
Appendix~\ref{A:SPECTRAL} contains some basic spectral theoretic statements about the linearised operator. Appendix~\ref{A:SCHAMEL} treats the
so-called multi-branched case, where we allow steady states with different dependency on $E$ depending on the sign of $v$, when $E>0$.
Appendix~\ref{A:SCALING} clarifies the scaling freedom in the problem.

\bigskip
{\bf Acknowledgements.}
M. Had\v zi\'c’s research is supported by the EPSRC Early Career Fellowship EP/S02218X/1. M. Moreno's research is supported by the 
National Agency for Research and Development (ANID)/Becas de Doctorado en el extranjero Becas Chile 72220280.


\section{Asymptotic behaviour of the period function at elliptic/hyperbolic points}\label{S:PERIOD}


The singular behaviour of the action-angle variables and the period function depends strongly on 
how close the energy is to a critical point and whether the critical point is hyperbolic (i.e. $E=0$) or elliptic (i.e. $E=\Emin$).
To handle this we introduce a constant
\begin{align}
\Emin <  E^\ast < 0
\end{align}
satisfying the following property: 
\begin{align}\label{E:PROP_E}
\varphi_0''(x)>0\text{ for all }[0,x_{-}(E^\ast)]\cup[x_{+}(E^\ast),1].
\end{align} 
Existence of such an $E^\ast$ follows from assumptions $(\varphi 1)$ and $(\varphi 2)$. 
Importantly,~\eqref{E:PROP_E} remains true if we replace $E^\ast$ by any $E'\in (E^\ast,0)$ so we can choose $E^\ast$ constant as close to zero as techincally required.
 Property~\eqref{E:PROP_E} will be important in describing the asymptotic properties of $T$, $T'$ and $T''$ close to the separatrix stated in Theorem~\ref{T:PERIODASYMP} below
 and all the universal constants will implicitly depend on $E^\ast$. Condition~\eqref{E:PROP_E} is a consequence of the hyperbolicity of the critical points $x=0$ and $x=1$ and Theorem~\ref{T:PERIODASYMP} thus links this geometric property to the asymptotic properties of the period function near the separatrix.

We next organise the energy levels inside $I$ into three overlapping sets $\Iell$, $\Ihypall$, and $\Iext$:
\begin{align}
I &= \Iell \cup \Ihypall \cup \Iext, \\ 
\Iell &:= (\Emin, \frac{E^\ast}{2}), \label{E:IELLDEF}\\ 
\Ihypall &:= (E^\ast, 0)\cup (0, |E^\ast |)  \label{E:IHYPDEF} \\
\Iext & : = (\frac{|E^\ast|}{2},\infty).
\end{align}
We shall often distinguish the bounds in the three regions above by using the suitable characteristic functions. We let 
\begin{align}\label{E:CHISDEF}
\che: = \chi_{\Iell}, \ \ \chh:=\chi_{\Ihypall}, \ \ \chext : = \chi_{\Iext}.
\end{align}

\begin{figure}[h!]
	\begin{center}
		\begin{tikzpicture}
			\tikzmath{\xscale = 2.5; \yscale = 3.;}
			
			\draw[-{>},black] (-2.6*\xscale,0.) -- (1.9*\xscale,0.) node[below] {$E$};
\coordinate[label=below:$\Emin$] (A) at (-2.6*\xscale,0.){};
\coordinate[label=below:$0$] (B) at (-0.8*\xscale,0.){};
\coordinate[label=below:$E^\ast$] (C) at (-1.8*\xscale,0.){};
\coordinate[label=below:$\frac{E^\ast}{2}$] (D) at (-1.3*\xscale,0.){};
\coordinate[label=below:$|E^\ast|$] (E) at (.2*\xscale,0.){};
\coordinate[label=below:$\frac{|E^\ast|}{2}$] (F) at (-0.32*\xscale,0.){};


\draw[fill=red] (A) circle (2pt);
\draw[fill=red] (B) circle (2pt);
\draw[fill=black] (C) circle (2pt);
\draw[fill=black] (D) circle (2pt);
\draw[fill=black] (E) circle (2pt);
\draw[fill=black] (F) circle (2pt);

		\end{tikzpicture}
	\end{center}
	\caption{The subdivision of the energy range into the three overlapping zones $\Iell$, $\Ihypall$, and $\Iext$.}
	\label{F:ERANGE}
\end{figure}
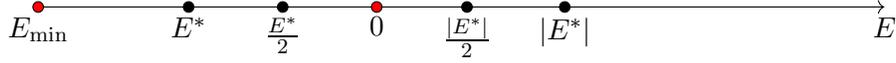




The asymptotic behavior of $T$, $T'$ and $T''$ is given by the following theorem
\begin{theorem}[Regularity, monotonicity and asymptotics of the period function $T$]
    \label{T:PERIODASYMP}
Let $\varphi_0$ satisfy the assumptions $(\varphi1)$--$(\varphi4)$.
\begin{enumerate}
\item For any $\varphi_0\in C^3(\mathbb T^1)$ we have $T\in C^{2}([\Emin,0)\cup(0,\infty))$. 
\item The period function $T:[\Emin,0)\cup (0,\infty)$ is strictly increasing on $[\Emin,0)$ and strictly decreasing
on $(0,\infty)$.
\item
Moreover, there exist $\Emin<E^\ast<0$ satisfying~\eqref{E:PROP_E} and a constant $C>0$ such that for any $E\in I$ the following bounds hold
    \begin{equation}
        \frac{1}{C} \Big(\big|\log|E|\big| \chh +  \che + |E|^{-\frac12} \chext \Big)\leq T(E)\leq C\Big( \big|\log|E|\big|\chh+  \che + |E|^{-\frac12}\chext\Big)\label{E:TASYMPBOUND}
    \end{equation}
    \begin{equation}
       \frac{1}{C} \Big(\frac{1}{|E|}\chh+ \che + |E|^{-\frac32}\chext\Big)\leq |T'(E)|\leq C \Big(\frac{1}{|E|}\chh+ \che+|E|^{-\frac32}\chext\Big)\label{E:TASYMPBOUND_DER}
    \end{equation}
    \begin{equation}
       \frac{1}{C}\Big(\frac{1}{|E|^2}\chh+\che+|E|^{-\frac52}\chext\Big) \leq |T''(E)|\leq C \Big(\frac{1}{|E|^2}\chh+ \che+|E|^{-\frac52}\chext\Big)\label{E:TASYMPBOUND_SECDER}
    \end{equation}
\end{enumerate}
\end{theorem}


\begin{remark}
Proof of Theorem~\ref{T:PERIODASYMP} can be found in Appendix~\ref{A:PERIOD}.
The asymptotic bounds~\eqref{E:TASYMPBOUND} were also shown by Després~\cite[Lemmas A.1, 6.2, 6.4 and 6.5]{DESP2019} for $\varphi_0\in W^{3,\infty}(\mathbb{T}^1)$. 
\end{remark}

The following lemma describes the degeneration rates of the energy-derivative of the flow near the critical points.

\begin{theorem}\label{T:DXDE}
Let $x(\th,E)$ be be given by~\eqref{E:XVTHETA} and introduce the shorthand
\begin{align}
\tE : = E-\Emin, \ \ E\in I.
\end{align}
 Then there exist $\Emin<E^\ast<0$ satisfying~\eqref{E:PROP_E} and a constant $C>0$  such that
\begin{align}
|\pa_Ex(\th,E)| \le C\tE^{-\frac12} \che(E) + \frac C{|E|}\chh(E) +\frac C{|E|}\chext(E), \ \ (\th,E)\in \mathbb S^1\times I. \notag
\end{align} 
\end{theorem}

\begin{proof}
If we differentiate the identity $x(\th(x,E),E)=x$ with respect to $E\in I$ we obtain
\begin{align}
\pa_E x(\theta(x,E),E) &= - \pa_\th x(\theta(x,E),E)\pa_E \th(x,E)\notag\\
& = - \pa_s X(\th(x,E)T(E),E) T(E) \pa_E\th (x,E).\notag
\end{align}
Since $|\partial_sX(\th(x,E)T(E),E)|=\sqrt{2(E+\varphi_0(x))}$, then
\begin{align}
|\pa_Ex(\th(x,E),E)|\leq \sqrt{2(E+\varphi_0(x))} \big(|\pa_E(T(E)\th(x,E))|+|T'(E)|\big). \notag
\end{align}
The claim now follows from Theorem~\ref{T:PERIODASYMP} and Lemma~\ref{L:RSLEMMA}.
\end{proof}

\begin{remark}
Theorem~\ref{T:DXDE} quantifies the loss of regularity of the change of variables $(x,v)\mapsto (\th,E)$ at the elliptic and the hyperbolic critical point. 
\end{remark}


\section{Absence of embedded eigenvalues}\label{S:ABSENCE}


We first express the linearised flow~\eqref{E:VLASOVLIN}--\eqref{E:POISSONLIN} in action angle variables. Due to different orientations of
the particles outside the separatrix characterised by the sign of the velocity (see~\eqref{E:DPOS}--\eqref{E:DNEG}), we introduce
\begin{align}
g_-(t,x,v) = g(t,x,v) \chi_{\{E(x,v)>0,v<0\}}, \ \ g_+ = g- g_-.
\end{align}
Upon moving to the action-angle variables $(\theta,E)$ and using~\eqref{E:DPOS}--\eqref{E:DNEG}, the linearised flow
transforms into
\begin{align}
\pa_t g_+ + \om(E) \pa_\theta \left(g_+ + \k^2 \mu'(\k^2E)\varphi_{g}\right) & = 0,  \label{E:VLASOVLIN1}\\ 
\pa_t g_- - \om(E) \pa_\theta \left(g_- + \k^2 \mu'(\k^2E)\varphi_{g}\right) & = 0, \label{E:VLASOVLIN2}\\
\pa_{xx}\varphi_g & = \eps \int g\diff v. \label{E:POISSONLIN1}
\end{align}

\begin{lemma}\label{L:FOURIERFORMULA}
Under the assumptions of Theorem~\ref{T:MAIN}
assume that $\lambda=\frac{2\pi i}{q}$ is an eigenvalue of the operator $\mathcal L$ with an associated eigenfunction $g\in\mathcal H$.
\begin{enumerate}
\item[(a)]
Then for any $(\ell,E)\in\Z_\ast\times I$ such that $T(E)\neq  \ell q$ we have the identity
\begin{align}\label{E:GFORMULA}
    \widehat{g_+}(\ell,E)&=-\varepsilon^2\frac{\mu'(\k^2E)\widehat{\varphi_g}(\ell,E)}{1-\frac{T(E)}{\ell q}},
    \end{align}
    and for any  $(\ell,E)\in\Z_\ast\times I_2$ such that $T(E)\neq  -\ell q$  we have
       \begin{align}
        \widehat{g_-}(\ell,E)& =-\varepsilon^2\frac{\mu'(\k^2E)\widehat{\varphi_g}(\ell,E)}{1+\frac{T(E)}{\ell q}}.\label{E:GMINUSFORMULA}
\end{align}
\item[(b)]
If on the other hand $T(E)=\pm\ell q$, then $\widehat{\varphi_g}(\ell,E)=0$.
\end{enumerate}
\end{lemma}


\begin{proof}
Let $(\ell,E)\in\mathbb Z_\ast\times I$ be such that $T(E)\neq\ell q$.
Note that the eigenvalue equation for $g_+$ reads
\begin{equation*}
    \frac{1}{T(E)}\partial_{\theta}(g_++ \k^2 \mu'(\k^2E)\varphi_g)=\frac{2\pi i }{q}g_+.
\end{equation*}
Taking the Fourier transform of this relation with respect to $\th\in\mathbb S^1$
we obtain
\begin{equation*}
    \frac{\ell}{T(E)}(\widehat{g_+}(\ell,E)+ \k^2 \mu'(\k^2E)\widehat{\varphi_g}(\ell,E))=\frac{1}{q}\widehat{g_+}(\ell,E),
\end{equation*}
and the claimed formula for $g_+$ is straightforward. Similarly we obtain~\eqref{E:GMINUSFORMULA}.
If on the other hand $T(E)=\pm \ell q$, we infer from the above that $\widehat{\varphi_g}(\ell,E)=0$
\end{proof}


\subsection{Frequency splitting}\label{S:FS}


It is evident from Lemma~\ref{L:FOURIERFORMULA} that the denominator $(\frac{T(E)}{\ell q}-1)^{-1}$ on the right-hand side of~\eqref{E:GFORMULA} requires a more refined understanding of how the (angle-) frequencies and (action-) energies can conspire to create small denominators.
We first introduce  a ``control" parameter $0<\delta<1$ such that
\begin{align}\label{eq:COND_DELTA}
\de<\min\left\{1,\frac{T(\frac{E^\ast}{2})-T(E^\ast)}{T(\frac{E^\ast}{2})},\frac{T(\frac{|E^\ast|}{2})-T(|E^\ast|)}{T(|E^\ast|)}\right\},
\end{align}
where we assume without lost that the constants $E^\ast$ used in Theorem~\ref{T:PERIODASYMP} and Theorem~\ref{T:DXDE} are equal.

{\em Trapped region $I_1$.}
For any given $\ell\in\Z_\ast$, if $\ell q\ge \Tmin$
we let $E_{\ell}^1\in I_1$ be the unique negative energy value such that $ T(E_{\ell}^1)=\ell q$ and $\Emin$ otherwise.
Similarly, if $(1\pm\de)\ell q\ge \Tmin$ we let $E_{\ell}^{1,\pm\de}\in I_1$ be the unique negative energy value such that $ T(E_{\ell}^{1,\pm\de})=(1\pm\de)\ell q$ and $\Emin$ otherwise. Since $T$ is increasing on $I_1$ we have $E_\ell^{1,-\de}\le E_{\ell}^{1}\le E_{\ell}^{1,\de}$. It then follows that for any $\ell\in\Z_\ast$ 
\begin{align}
&0<\frac{T(E)}{\ell q}-1<\delta \ \ \text{ if and only if } \ E \in (E_{\ell}^{1},E_{\ell}^{1,\de}) \label{E:INT1}\\
&0<1-\frac{T(E)}{\ell q}<\delta \ \ \text{ if and only if } \ E \in (E_{\ell}^{1,-\de},E_{\ell}^{1}) \label{E:INT2} \\
&\frac{T(E)}{\ell q}-1<0 \text { if } \  E \in (\Emin,E_{\ell}^{1}).
\end{align}

We now note that for any $\ell>0$, condition~\eqref{eq:COND_DELTA} implies that the interval $[E_{\ell}^1,E_{\ell}^{1,\de}]$ is either empty (if $\Eml=\Emld=\Emin$) or contained entirely inside the elliptic region $\Iell$ or inside the hyperbolic region $\Ihypall\cap I_1$ (and possibly in both). Analogous comment applies to $[E_{\ell}^{1,-\de},E_{\ell}^1]$. We let 
\begin{align}
\ell^\ast=\lfloor \frac{T(E^\ast)}{q}\rfloor+1
\end{align} 
be the smallest value of $\ell>0$, condition such that $[E_{\ell}^1,E_{\ell}^{1,\delta}]\subset(E^\ast,0)=\Ihyp\cap I_1$. Then for $\ell<\ell^\ast$ we have $[E_{\ell}^1,E_{\ell}^{1,\delta}]\subset[E_{\min},E^\ast/2]=\Iell$. 

{\em Outside region.}
Analogously we let $E_{\ell}^2\in I_2, E_{\ell}^{2,\delta}\in I_2$ be two unique positive energy levels such that 
\begin{align}
T(E_{\ell}^2)= \ell q, \ \ T(E_{\ell}^{2,\delta})=(1+\delta)\ell q.
\end{align}
Since $T$ is decreasing on $I_2$ we have $E_{\ell}^{2,\delta}<E_{\ell}^{2}$. It then follows that for any $\ell\in\Z_\ast$ 
\begin{align}
&0<\frac{T(E)}{\ell q}-1<\delta \ \ \text{ if and only if } \ E \in (E_{\ell}^{2,\de},E_{\ell}^{2}), \\
& \frac{T(E)}{\ell q}-1<0 \text { if } \  E\in(E_{\ell}^{2},\infty). \label{E:INT2.2}
\end{align}
We now note that for any $\ell>0$,  condition~\eqref{eq:COND_DELTA} implies that the interval $[E_{\ell}^{2,\delta},E_{\ell}^{2}]$ is contained entirely inside the hyperbolic region $\Ihypall\cap I_2$ or inside the exterior region $\Iext$ (and possibly in both). We let 
\begin{align}
L^\ast=\lfloor \frac{T(|E^\ast|)}{q}\rfloor+1 \label{E:BIGLMSTARDEF}
\end{align} 
be the smallest value of $\ell>0$ such that $[E_{\ell}^{2,\delta},E_{\ell}^2]\subset(0,|E^\ast|)=\Ihypall\cap I_2$. Then for $\ell<L^\ast$ we have $[E_{\ell}^{2,\delta},E_{\ell}^2] \subset[\frac{|E^\ast|}{2},\infty)=\Iext$. 


\subsection{Energy estimates}\label{S:EE}


\begin{lemma}\label{L:ENERGYBASIC}
Under the assumptions of Theorem~\ref{T:MAIN} assume that $\lambda=\frac{2\pi i }{q}$ is an eigenvalue of the operator $\mathcal L$ and assume without loss that $q>0$. Denote the associated eigenfunction by $g\in\mathcal H$.
\begin{enumerate}
\item (Monotone case) If $\mu: I\to \mathbb R_+$ is monotonically decreasing then
\begin{equation}\label{E:BOUNDONE}
    \|\partial_x\varphi_g\|_{L^2}^2\leq \eps^3\left(\mathcal{J}_{\text{nr}}+\mathcal{J}_{\text{trap}}+\mathcal{J}_{\text{ext}}\right), 
\end{equation}
where
\begin{align}
    \mathcal{J}_{\text{nr}}& =\sum_{\ell\in\mathbb{Z}_\ast^{>0}}\int_{E_{\ell}^{1,\de}}^{0} \frac{T(E)|\mu'(\k^2E)|}{\frac{T(E)}{\ell q}-1}|\widehat{\varphi}_g(\ell, E)|^2\diff E \notag\\
    & \ \ \ \ + 2 \sum_{\ell>0}\int_{0}^{E_{\ell}^{2,\de}} \frac{T(E) |\mu'(\k^2E)|}{\frac{T(E)}{\ell q}-1}\big( |\widehat{\varphi}_g(\ell, E)|^2+|\widehat{\varphi}_g(-\ell, E)|^2 \big) \diff E,\label{E:JNRDEF}\\
    \mathcal{J}_{\text{trap}} & =\sum_{\ell>0}\int_{E_{\ell}^{1}}^{E_{\ell}^{1,\de}}\frac{T(E)|\mu'(\k^2E)|}{\frac{T(E)}{\ell q}-1}|\widehat{\varphi}_g(\ell, E)|^2\diff E, \label{E:JTRAPDEF}\\
    \mathcal{J}_{\text{ext}} &= 2 \sum_{\ell>0}\int_{E_{\ell}^{2,\de}}^{E_{\ell}^2} \frac{T(E) |\mu'(\k^2E)|}{\frac{T(E)}{\ell q}-1}\big( |\widehat{\varphi}_g(\ell, E)|^2+|\widehat{\varphi}_g(-\ell, E)|^2 \big) \diff E.
    \label{E:JEXTDEF}
\end{align}
\item (Non-monotone case)  If $\mu:I_1\to\mathbb R_+$ is monotone increasing and $\mu:I_2\to\mathbb R_+$ monotone decreasing, then
\begin{equation}\label{E:BOUNDTWO}
    \|\partial_x\varphi_g\|_{L^2}^2\leq \eps^3\left(\mathcal{J}_{\text{nr},\uparrow}+\mathcal{J}_{\text{trap},\uparrow}+\mathcal{J}_{\text{ext}}\right), 
\end{equation}
where
\begin{align}
    \mathcal{J}_{\text{nr},\uparrow}& =\sum_{\ell>0}\int_{\Emin}^{E_\ell^{1,-\de}} \frac{T(E)|\mu'(\k^2E)|}{1-\frac{T(E)}{\ell q}}|\widehat{\varphi}_g(\ell, E)|^2\diff E \notag\\
    & \ \ \ \  +\sum_{\ell>0}\int_{I_1}\frac{T(E)|\mu'(\k^2E)|}{\frac{T(E)}{\ell q}+1}|\widehat{\varphi}_g(-\ell, E)|^2\diff E \notag\\
    & \ \ \ \ + 2 \sum_{\ell>0}\int_{0}^{E_{\ell}^{2,\de}} \frac{T(E) |\mu'(\k^2E)|}{\frac{T(E)}{\ell q}-1}\big( |\widehat{\varphi}_g(\ell, E)|^2+|\widehat{\varphi}_g(-\ell, E)|^2 \big) \diff E,\label{E:JNRDEFUP}\\
    \mathcal{J}_{\text{trap},\uparrow} & =\sum_{\ell>0}\int_{E_{\ell}^{1,-\de}}^{E_{\ell}^{1}}\frac{T(E)|\mu'(\k^2E)|}{1-\frac{T(E)}{\ell q}}|\widehat{\varphi}_g(\ell, E)|^2\diff E,\label{E:JTRAPDEFUP}
    \end{align}
\end{enumerate}
Here quantities $E_{\ell}^i$, $E_{\ell}^{i,\pm\de}$, $i=1,2$ are introduced in Section~\ref{S:FS}.
\end{lemma}


\begin{proof}
Integrating-by-parts and using Lemma~\ref{L:FOURIERFORMULA} we obtain
\begin{align}
    \|\partial_x\varphi_g\|_{L^2}^2
    &=-\int_{0}^{1}\varphi_g\partial_x^2\varphi_g\diff x
    =-\k\iint_{\Omega}\varphi_g(x)g(x,v)\diff v\diff x\\
    & = - \k\iint_{\Omega\cap\{E<0\}}\varphi_g(x)g(x,v)\diff v\diff x \notag\\
    & \ \ \  -\k\iint_{\Omega\cap\{E(x,v)>0, v>0\}}\varphi_g(x)g_+(x,v)\diff v\diff x
    -\k\iint_{\Omega\cap\{E(x,v)>0,v<0\}}\varphi_g(x)g_-(x,v)\diff v\diff x \\
    &=-\k\int_{\Emin}^0\int_{\mathbb S^1}\varphi_g(\theta,E)g(\theta,E)T(E)\diff\theta \diff E -\k\int_{0}^\infty\int_{\mathbb S^1}\varphi_g(\theta,E)g_+(\theta,E)T(E)\diff\theta \diff E  \notag \\
    & \ \ \ \ -\k\int_{0}^\infty\int_{\mathbb S^1}\varphi_g(\theta,E)g_-(\theta,E)T(E)\diff\theta \diff E, \notag\\
    \end{align}
    where we have changed variables from phase space to action-angle coordinates $\diff v\diff x= T(E)\diff\theta\diff E$. We next
    use Plancharel's identity in $\th$-variable, and formulas~\eqref{E:GFORMULA}--\eqref{E:GMINUSFORMULA} to get
    \begin{align}
  \|\partial_x\varphi_g\|_{L^2}^2  &=-\varepsilon^3\sum_{\ell\in\Z_\ast}\int_{I_1}\frac{T(E)\mu'(\k^2E)}{\frac{T(E)}{\ell q}-1}|\widehat{\varphi}_g(\ell, E)|^2\diff E \notag\\
& \ \ \ \    -\varepsilon^3\sum_{\ell\in\Z_\ast}\int_{I_2} T(E) \mu'(\k^2E)  \big(\frac{1}{\frac{T(E)}{\ell q}-1}-\frac{1}{\frac{T(E)}{\ell q}+1}\big)
    |\widehat{\varphi}_g(\ell, E)|^2 \diff E \notag\\
    & = -\varepsilon^3\sum_{\ell\in\Z_\ast}\int_{I_1}\frac{T(E)\mu'(\k^2E)}{\frac{T(E)}{\ell q}-1}|\widehat{\varphi}_g(\ell, E)|^2\diff E \notag\\
& \ \ \ \    -2\varepsilon^3\sum_{\ell\in\Z_\ast}\int_{I_2} \frac{T(E) \mu'(\k^2E)}{(\frac{T(E)}{\ell q}-1)(\frac{T(E)}{\ell q}+1)} |\widehat{\varphi}_g(\ell, E)|^2 \diff E \notag \\
& =- \varepsilon^3\sum_{\ell\in\Z_\ast}\int_{I_1}\frac{T(E)\mu'(\k^2E)}{\frac{T(E)}{\ell q}-1}|\widehat{\varphi}_g(\ell, E)|^2\diff E \notag\\
& \ \ \ \    -2\varepsilon^3\sum_{\ell>0}\int_{I_2} \frac{T(E) \mu'(\k^2E)}{(\frac{T(E)}{\ell q}-1)(\frac{T(E)}{\ell q}+1)}\big( |\widehat{\varphi}_g(\ell, E)|^2+|\widehat{\varphi}_g(-\ell, E)|^2 \big) \diff E. \label{E:ENERGYONE}
\end{align}

If $E\in I_2$, then by definition $\frac{T(E)}{\ell q}-1<0$ if and only if $E\in[E_{\ell}^2,\infty)$. Moreover, since $\frac{T(E)}{\ell q}+1\ge 1$ for all $\ell>0$ and $-\mu'(\k^2 E)=|\mu'(\k^2E)|$ on $I_2$, we may bound the last line of~\eqref{E:ENERGYONE} by 
\begin{align}
&2\varepsilon^3\sum_{\ell>0}\int_{0}^{E_{\ell}^2} \frac{T(E) |\mu'(\k^2E)|}{\frac{T(E)}{\ell q}-1}\big( |\widehat{\varphi}_g(\ell, E)|^2+|\widehat{\varphi}_g(-\ell, E)|^2 \big) \diff E \notag\\
&=2\varepsilon^3\sum_{\ell>0}\big(\int_{0}^{E_{\ell}^{2,\de}}+\int_{E_{\ell}^{2,\de}}^{E_{\ell}^2}\big) \frac{T(E) |\mu'(\k^2E)|}{\frac{T(E)}{\ell q}-1}\big( |\widehat{\varphi}_g(\ell, E)|^2+|\widehat{\varphi}_g(-\ell, E)|^2 \big) \diff E. \label{E:ENERGYAUX2}
\end{align}

{\em Case 1: Assume that $\mu$ is decreasing on $I_1$.}
Note that $\frac{T(E)}{\ell q}-1<0$ for all $\ell<0$. Clearly if $E\in I_1$ then by definition $\frac{T(E)}{\ell q}-1<0$ if $E\in[\Emin,E_{\ell}^1]$ and since $\mu'<0$ on $I_1$ the first line on the right-most side of~\eqref{E:ENERGYONE} is bounded by
\begin{align}
 \varepsilon^3\sum_{\ell\in\Z_\ast}\int_{I_1}\frac{T(E)|\mu'(\k^2E)|}{\frac{T(E)}{\ell q}-1}|\widehat{\varphi}_g(\ell, E)|^2\diff E 
 &\le  \varepsilon^3\sum_{\ell>0}\int_{E_{\ell}^1}^0\frac{T(E)|\mu'(\k^2E)|}{\frac{T(E)}{\ell q}-1}|\widehat{\varphi}_g(\ell, E)|^2\diff E \notag\\
 & = \varepsilon^3\sum_{\ell>0}\big(\int_{E_{\ell}^1}^{E_{\ell}^{1,\de}}+\int_{E_{\ell}^{1,\de}}^0\big)\frac{T(E)|\mu'(\k^2E)|}{\frac{T(E)}{\ell q}-1}|\widehat{\varphi}_g(\ell, E)|^2\diff E. \label{E:ENERGYAUX1}
\end{align} 
Using~\eqref{E:ENERGYONE},~\eqref{E:ENERGYAUX1}, and~\eqref{E:ENERGYAUX2} the bound~\eqref{E:BOUNDONE} follows.

{\em Case 2: Assume now that $\mu$ is increasing on $I_1$.} For negative $\ell<0$ we have
\[
- \varepsilon^3\sum_{\ell<0}\int_{I_1}\frac{T(E)\mu'(\k^2E)}{\frac{T(E)}{\ell q}-1}|\widehat{\varphi}_g(\ell, E)|^2\diff E
= \varepsilon^3\sum_{\ell>0}\int_{I_1}\frac{T(E)\mu'(\k^2E)}{\frac{T(E)}{\ell q}+1}|\widehat{\varphi}_g(-\ell, E)|^2\diff E>0,
\]
and this accounts for the second line of~\eqref{E:JNRDEFUP}. For $\ell>0$ only the terms with $\frac{T(E)}{\ell q}-1<0$ give a positive contribution. We split such indices in the range 
$-\de<\frac{T(E)}{\ell q}-1<0$ and $\frac{T(E)}{\ell q}-1\le-\de$, which leads to $E\in(\Emin,E^{1,-\de}_{\ell}]$ and $E\in(E^{1,-\de}_{\ell},E^1_{\ell})$ respectively. This accounts for the first line of~\eqref{E:JNRDEFUP} and~\eqref{E:JTRAPDEFUP}.
\end{proof}

\begin{remark}
This lemma splits the analysis into the non-resonant contribution $\mathcal{J}_{\text{nr}}$, and the two resonant contributions $\mathcal{J}_{\text{trap}}$ and $\mathcal{J}_{\text{ext}}$. These two contributions are respectively concentrated in the trapped energy region $I_1$ and the outside energy region $I_2$. The key challenge in our analysis is to estimate these terms, as they will ``feel" different dynamics of the characteristic system near the two critical points $\Emin$ and $E=0$,  and at high energies $E=\infty$.
\end{remark}


\subsubsection{Estimates in the trapped region: $\mathcal J_{\text{trap}}$-term}


The hardest part of our analysis is the treatment of the term $\mathcal J_{\text{trap}}$~\eqref{E:JTRAPDEF} that involves the estimates near both the elliptic point $\Emin$ and the hyperbolic point $E=0$.
To facilitate the ``small denominators" analysis, 
for any $\ell\in\Z_\ast$ and any $t>\ell q$  we define
\begin{equation}\label{E:PMLDEF}
    p_{m,\ell}(t):=\ell q    \log\left(\frac{t}{\ell q}-1\right), 
\end{equation}
so that the following relation holds
\begin{equation*}
    \frac{1}{T'(E)}\frac{d}{dE} \big(p_{m,\ell}(T(E))\big)=\frac{1}{\frac{T(E)}{\ell q}-1}.
\end{equation*}

The central result of this section is the following lemma.

\begin{lemma}\label{L:TRAPLEMMA}
Under the assumptions of Lemma~\ref{L:ENERGYBASIC}, there exists a constant $C>0$ such that
\begin{align}\label{E:TRAPBOUND}
\mathcal J_{\text{trap}} & \le C \k^2 \|\partial_x\varphi_g\|_{L^2}^2 \ \ \text{ if $\mu$ is decreasing on $I_1$}, \\
\mathcal J_{\text{trap}}^{\uparrow} & \le C \k^2 \|\partial_x\varphi_g\|_{L^2}^2 \ \ \text{ if $\mu$ increasing on $I_1$}, \label{E:TRAPBOUNDUP}
\end{align}
where we recall definitions~\eqref{E:JTRAPDEF} and~\eqref{E:JTRAPDEFUP} of $\mathcal J_{\text{trap}}$ and $\mathcal J_{\text{trap}}^{\uparrow}$ respectively.
\end{lemma}


\begin{proof}
We observe that by assumption~\eqref{E:MUMINUSPLUSREG} we have the uniform bound
\begin{align}\label{E:MUUNIF}
|\mu(\eps^2 E)|+ |\mu'(\eps^2 E)| + |\mu''(\eps^2 E)| \lesssim 1, \ \ E\in I_1,
\end{align}
which will be frequently used throughout this section. We only prove the estimate~\eqref{E:TRAPBOUND}, as~\eqref{E:TRAPBOUNDUP}
follows analogously.

We integrate by parts in~\eqref{E:JTRAPDEF} to obtain
\begin{align}
\mathcal{J}_{\text{trap}}  &=  \sum_{\ell>0}\int_{E_{\ell}^{1}}^{E_{\ell}^{1,\de}}\frac{T(E)|\mu'(\k^2E)|}{\frac{T(E)}{\ell q}-1}|\widehat{\varphi}_g(\ell, E)|^2\diff E \notag\\
  &=\sum_{\ell>0}\frac{T(E)p_{m,\ell}(T(E))}{T'(E)}|\mu'(\k^2E)||\widehat{\varphi}_g(\ell, E)|^2\Bigg|_{E_{\ell}^{1}}^{E_{\ell}^{1,\de}} \notag \\
    & \ \ \ \ -\sum_{\ell>0}\int_{E_{\ell}^{1}}^{E_{\ell}^{1,\de}}\frac{\diff}{\diff E}\left(\frac{T(E)|\mu'(\k^2E)||\widehat{\varphi}_g(\ell, E)|^2}{T'(E)}\right)p_{m,\ell}(T(E))\diff E \label{E:TRAP1}
\end{align}

Due to the negativity of $\pml(t)$ for $t\in(\ell q,2\ell q)$, $\ell,q>0$, we see that the upper boundary of the boundary-term contribution is simply negative for $\de<1$, since $T(E_{\ell}^{1,\de})=(1+\de)\ell q$ and therefore bounded from above by $0$. On the other hand
\begin{align}\label{E:BDRY1}
\lim_{E\to E_{\ell}^1}\frac{T(E)p_{m,\ell}(T(E))}{T'(E)}|\mu'(\k^2E)||\widehat{\varphi}_g(\ell, E)|^2 = 0.
\end{align}
The above follows from Lemma~\ref{L:FOURIERFORMULA} which implies $\widehat{\varphi}_g(\ell, E_\ell^1)=0$ since $T(E_\ell^1)=\ell q$ by definition. 
Moreover, the function $E\mapsto \widehat{\varphi}_g(\ell, E)$ is differentiable at any $E\in I_1$ and therefore vanishes algebraically at $E_\ell^1$. Therefore $\lim_{E\to E_\ell^1}\Big[\log(\frac{T(E)}{\ell q}-1)\widehat{\varphi}_g(\ell, E)\Big]=0$, thus implying~\eqref{E:BDRY1}.

Therefore, it follows from~\eqref{E:TRAP1} that 
\begin{align}
   \mathcal{J}_{\text{trap}}
  & \leq
    -\sum_{\ell>0}\int_{E_{\ell}^{1}}^{E_{\ell}^{1,\de}}\frac{\diff}{\diff E}\left(\frac{T(E)|\mu'(\k^2E)||\widehat{\varphi}_g(\ell, E)|^2}{T'(E)}\right)p_{m,\ell}(T(E))\diff E \notag \\
    & =\sum_{\ell>0}(A_\ell+B_\ell), \label{E:JTRAPBOUND}
\end{align}
where for any $\ell>0$
\begin{align}
    &A_\ell=-2\int_{E_{\ell}^{1}}^{E_{\ell}^{1,\de}}\frac{T(E)|\mu'(\k^2E)|}{T'(E)}\text{Re}\left(\overline{\widehat{\varphi}_g(\ell, E)}\frac{\diff}{\diff E}\widehat{\varphi}_g(\ell, E)\right)p_{m,\ell}(T(E))\diff E, \label{E:AELLDEF}\\
    &B_\ell=-\int_{E_{\ell}^{1}}^{E_{\ell}^{1,\de}}\frac{\diff}{\diff E}\left(\frac{T(E)|\mu'(\k^2E)|}{T'(E)}\right)|\widehat{\varphi}_g(\ell, E)|^2p_{m,\ell}(T(E))\diff E. \label{E:BELLDEF}
\end{align}

Observe that for any $\ell >0$ we have
\begin{align*}
    |A_\ell|&\leq 2\int_{E_{\ell}^{1}}^{E_{\ell}^{1,\de}}\frac{T(E)|\mu'(\k^2E)|}{T'(E)}\left|\text{Re}\left(\overline{\widehat{\varphi}_g(\ell, E)}\frac{\diff}{\diff E}\widehat{\varphi}_g(\ell, E)\right)\right| |p_{m,\ell}(T(E))|\diff E\\
    &\leq 2\int_{E_{\ell}^{1}}^{E_{\ell}^{1,\de}}\frac{T(E)|\mu'(\k^2E)|}{T'(E)}\left|\widehat{\varphi}_g(\ell, E)\right|\left|\frac{\diff}{\diff E}\widehat{\varphi}_g(\ell, E)\right||p_{m,\ell}(T(E))|\diff E.
\end{align*}

We use Cauchy-Schwarz  to get the following bound
\begin{align*}
    |A_\ell|^2&\leq 4\left(\int_{E_{\ell}^{1}}^{E_{\ell}^{1,\de}}\frac{T(E)|\mu'(\k^2E)|}{T'(E)}\left|\widehat{\varphi}_g(\ell, E)\right|\left|\frac{\diff}{\diff E}\widehat{\varphi}_g(\ell, E)\right||p_{m,\ell}(T(E))|\diff E\right)^2\\
    &\leq 4\int_{E_{\ell}^{1}}^{E_{\ell}^{1,\de}}\frac{p_{m,\ell}(T(E))^2|\widehat{\varphi}_g(\ell, E)|^2T(E)|\mu'(\k^2E)| \tE^\beta}{T'(E)^{1-\alpha}}\diff E\int_{E_{\ell}^{1}}^{E_{\ell}^{1,\de}}\frac{\left|\frac{\diff}{\diff E}\widehat{\varphi}_g(\ell, E)\right|^2 T(E)|\mu'(\k^2E)|}{T'(E)^{1+\alpha}\tE^{\beta}} \diff E,
\end{align*}
where $\alpha\in(\frac12,2)$ and $\beta\in(0,1)$ are given parameters, and we remind the reader that $\tE = E-\Emin$.
By Cauchy-Schwarz again we obtain 
\begin{equation}\label{E:AELLBOUND}
    |A_\ell|\leq 2\sqrt{A_{\ell,1}A_{\ell,2}}\leq A_{\ell,1}+ A_{\ell,2}
\end{equation}
where 
\begin{equation}
    A_{\ell,1}:=\int_{E_{\ell}^{1}}^{E_{\ell}^{1,\de}}\frac{p_{m,\ell}(T(E))^2|\widehat{\varphi}_g(\ell, E)|^2T(E)|\mu'(\k^2E)|\tE^\beta}{T'(E)^{1-\alpha}}\diff E \label{E:AELL1}
\end{equation}
and
\begin{equation}
    A_{\ell,2}:=\int_{E_{\ell}^{1}}^{E_{\ell}^{1,\de}}\frac{\left|\frac{\diff}{\diff E}\widehat{\varphi}_g(\ell, E)\right|^2 T(E)|\mu'(\k^2E)|}{T'(E)^{1+\alpha}\tE^\beta}\diff E\label{E:AELL2TRAP}
\end{equation}

{\bf Estimates on $A_{\ell,2}$.} We recall that for $\ell\ge\ell^\ast$ $[\Eml,\Emld]\subset \Ihypall\cap I_1=(E^\ast,0)$ which implies $\E\gtrsim1$ on $(E^\ast,0)$ and therefore
\begin{align}
\sum_{\ell\ge\ell^\ast}|A_{\ell,2}| &\le 
C \sum_{\ell\in\Z_\ast} \int_{E^\ast}^0 \lv\frac{\diff}{\diff E}\widehat{\varphi_g}(\ell,E)\rv^2 \frac{T(E)}{T'(E)^{1+\alpha}}\, |\mu'(\k^2E)| \,\diff E \notag\\
&= C\int_{(E^\ast,0)\times\mathbb S^1} \lv\frac{\diff}{\diff E}\varphi_g(\theta,E)\rv^2 \frac{T(E)}{T'(E)^{1+\alpha}}\, |\mu'(\k^2E)| \,\diff (\theta,E), \label{E:AELL21}
\end{align}
where we have used Plancharel in the last line.
We next use Lemma~\ref{T:DXDE} to obtain the bound
\begin{align}\label{E:PAEPHIBOUND}
|\pa_E\varphi_g|= |\pa_x\varphi_g|\left|\frac{\pa x}{\pa E}\right| \le \frac{C}{|E|} |\pa_x\varphi_g|,
\ \ E\in (E^\ast,0),
\end{align}
which together with the bound $\frac1{|E|}\lesssim T'(E)\lesssim \frac1{|E|}$ from Theorem~\ref{T:PERIODASYMP} in turn implies
\begin{align}
\sum_{\ell\ge\ell^\ast} |A_{\ell,2}|
& \le C \int^1_0 |\pa_x\varphi_g(x)|^2 \int_{\max\{E^\ast,-\varphi_0(x)\}}^0 \frac{|\mu'(\k^2E)||E|^{\alpha-1}}{ \sqrt{2(E+\varphi_0(x))}} \diff E \diff x,\label{E:ABOUND0}
\end{align}
where we have changed variables $(\theta,E)\mapsto (x,E)$ and used $\frac{d\th}{dx}= \frac1{T(E)\sqrt{E+\varphi_0(x)}}$. We now change vaiables
\begin{align}
|E|=\varphi_0(x)s, \ \ \diff E = - \varphi_0(x)\diff s,
\end{align}
to obtain
\begin{align}
 \int_{\max\{E^\ast,-\varphi_0(x)\}}^0 \frac{|\mu'(\k^2E)||E|^{\alpha-1}}{ \sqrt{2(E+\varphi_0(x))}} \diff E 
& \le C  \int_{\max\{E^\ast,-\varphi_0(x)\}}^0 \frac{|E|^{\alpha-1}}{ \sqrt{2(E+\varphi_0(x))}} \diff E  \notag\\
& =\varphi_0(x)^{\al-\frac12} \int_0^{\min\{1,\frac{|E^\ast|}{\varphi_0(x)}\}} \frac{s^{\al-1}}{\sqrt{1-s}}\diff s \notag\\
& \le C \varphi_0(x)^{\al-\frac12},
\end{align}
where we have used the bound $|\mu'(\eps^2E)|\lesssim1$.
Therefore for any $\al>\frac12$, due to boundedness of $|\varphi_0|$ we conclude from~\eqref{E:ABOUND0} that 
\begin{align}
\sum_{\ell\ge\ell^\ast} |A_{\ell,2}| \le C\|\pa_x\varphi_g\|_{L^2}^2. \label{E:ABOUND}
\end{align}
On the other hand, an analogous argument using Lemma~\ref{T:DXDE} and Theorem~\ref{T:PERIODASYMP} in the elliptic region gives
\begin{align}
\sum_{\ell=1}^{\ell^\ast} |A_{\ell,2}|
& \le C \int^1_0 |\pa_x\varphi_g(x)|^2 \int_{\Emin}^{\min\{\frac{E^\ast}{2},-\varphi_0(x)\}} \frac{|\mu'(\k^2E)|\tE^{\beta-\frac12}}{ \sqrt{2(E+\varphi_0(x))}} \diff E \diff x \notag \\
&\le C\|\pa_x\varphi_g\|_{L^2}^2, \label{E:ABOUNDELL}
\end{align}
where we have used $\beta>0$ to guarantee the integrability of the interior $E$-integral.

{\bf Estimates on $A_{\ell,1}$.} If we change variables $\theta\mapsto x$ and noticing that $\frac{\partial\varphi_g}{\partial\theta}=\frac{\partial\varphi_g}{\partial x}\frac{\partial x}{\partial\theta}$ 
we have
\begin{align*}
    \sum_{\ell>0}A_{\ell,1} &=\sum_{\ell>0}\int_{E_{\ell}^{1}}^{E_{\ell}^{1,\de}}\frac{p_{m,\ell}(T(E))^2|\widehat{\varphi}_g(\ell, E)|^2T(E)|\mu'(\k^2E)|}{T'(E)^{1-\alpha}\tE^\beta} \diff E\\
    &\leq\sum_{\ell>0}\frac{1}{4\pi^2\ell^2}\int_{E_{\ell}^{1}}^{E_{\ell}^{1,\de}}\int_{\mathbb{S}^1}\frac{p_{m,\ell}(T(E))^2|\partial_\theta\varphi_g(\theta,E)|^2T(E)|\mu'(\k^2E)|}{T'(E)^{1-\alpha}\tE^\beta}\diff E\diff\theta  \\
    &=\sum_{\ell>0}\frac{1}{4\pi^2\ell^2}\int_{x_{-}(E_{\ell}^{1,\de})}^{x_{+}(E_{\ell}^{1,\de})}|\partial_x\varphi_g(x)|^2\int_{\max\{E_{\ell}^{1},-\varphi_0(x)\}}^{E_{\ell}^{1,\de}}\frac{p_{m,\ell}(T(E))^2|\mu'(\k^2E)|T(E)^2\sqrt{2E+2\varphi_0(x)}}{T'(E)^{1-\alpha}\tE^\beta}\diff E \diff x\\
    &\leq \sum_{\ell>0}\frac{C}{4\pi^2\ell^2}\int_{0}^{1}|\partial_x\varphi_g(x)|^2\int_{\max\{E_{\ell}^{1},-\varphi_0(x)\}}^{E_{\ell}^{1,\de}}\frac{p_{m,\ell}(T(E))^2|\mu'(\k^2E)|T(E)^2}{T'(E)^{1-\alpha}\tE^\beta}\diff E\diff x\\
    &\leq C \|\partial_x\varphi_g\|_{L^2}^2\sum_{\ell>0}\frac{1}{\ell^2}\int_{E_{\ell}^{1}}^{E_{\ell}^{1,\de}}\frac{p_{m,\ell}(T(E))^2|\mu'(\k^2E)|T(E)^2}{T'(E)^{1-\alpha}\tE^\beta}\diff E,
\end{align*}
where $C>0$ is an $\ell$-independent constant. We have in the process used the identity $\frac{\pa x}{\pa\th}=\pm T(E)\sqrt{E+\varphi_0(x)}$ and the trivial bound $\sqrt{E+\varphi_0(x)}\lesssim 1$.
Thus we have
\begin{equation}
    \sum_{\ell>0}A_{\ell,1}\leq C\|\partial_x\varphi_g\|_{L^2}^2\sum_{\ell>0}\frac{1}{\ell^2}I_{m,\ell}, \label{E:AELLONEX}
\end{equation}
where we define
\begin{equation*}
    I_{m,\ell}:=\int_{E_{\ell}^{1}}^{E_{\ell}^{1,\de}}\frac{p_{m,\ell}(T(E))^2|\mu'(\k^2E)|T(E)^2}{T'(E)^{1-\alpha}\tE^\beta} \diff E
    =\ell^2q^2\int_{E_{\ell}^{1}}^{E_{\ell}^{1,\de}}\frac{\left|\log\left(\frac{T(E)}{\ell q}-1\right)\right|^2|\mu'(\k^2E)|T(E)^2}{T'(E)^{1-\alpha}\tE^\beta}    \diff E,
\end{equation*}
where we recall $q$ from Lemma~\ref{L:ENERGYBASIC}.

Just like above, for $\ell\ge\ell^\ast$ $[\Eml,\Emld]\subset \Ihypall\cap I_1=(E^\ast,0)$ which implies $\E\gtrsim1$ on $(E^\ast,0)$.
Moreover, by Theorem~\ref{T:PERIODASYMP} we also have
\begin{align}\label{E:KEYHYP}
\exp(T(E)) \lesssim T'(E) \lesssim  \exp(T(E)), \ \ E\in(E^\ast,0).
\end{align}
Therefore
\begin{align}
\sum_{\ell\ge\ell^\ast}\frac{1}{\ell^2}I_{m,\ell}&
\le C q^2\sum_{\ell\ge\ell^\ast}\int_{E_{\ell}^{1}}^{E_{\ell}^{1,\de}}\frac{\left|\log\left(\frac{T(E)}{\ell q}-1\right)\right|^2|\mu'(\k^2E)|T(E)^2}{T'(E)^{1-\alpha}}    \diff E \notag \\
& =   C q^5\sum_{\ell\ge\ell^\ast}\ell^3  \int_{1}^{1+\delta}\left|\log\left(t-1\right)\right|^2 t^2 \, e^{-(2-\alpha) \ell q t}\diff t \notag \\
&\leq C q^5\sum_{\ell\ge\ell^\ast}\ell^3  e^{-(2-\alpha) \ell q }\int_{1}^{1+\de}\left|\log(t-1)\right|^2t^2\diff t \notag \\
& \le C  q^5\sum_{\ell\ge\ell^\ast}\ell^3  e^{-\kappa_\al \ell q }, \label{E:AL1}
\end{align}
where $\kappa_\al:=2-\al>0$, we have used $|\mu'(\eps^2E)|\lesssim1$, and we have changed variables via
\be\label{E:CH}
\frac{T(E)}{\ell q} = t.
\ee
We assume that $\al<2$ so that $\kappa_\al>0$. 
By the integral criterion
\begin{align}
\sum_{\ell>0} \ell^3  e^{-\kappa_\al \ell q }& \le e^{-\kappa_\al q}+\int_1^\infty x^3 e^{-\kappa_\al q x}\diff x=  e^{-\kappa_\al  q }+ q^{-4} \int_{q}^\infty x^3 e^{-\kappa_\al x}\diff x\notag\\
&\le e^{-\kappa_\al q} +C_\al q^{-4} \int_{q}^\infty y^{-2} \diff y  \le e^{-\kappa_\al q}+ C_\alpha q^{-5} . \label{E:IC1}
\end{align}
We feed this back into~\eqref{E:AL1} and conclude
$
\sum_{\ell\ge\ell^\ast}\frac{1}{\ell^2}I_{m,\ell} \le C 
$
and therefore from~\eqref{E:AELLONEX} we get
\begin{align} \label{E:AL1BOUND}
\sum_{\ell\geq\ell^\ast}A_{\ell,1} \le C \|\partial_x\varphi_g\|_{L^2}^2.
\end{align}  

In the range $1\le\ell<\ell^\ast$, due to \eqref{eq:COND_DELTA} we have $E_{\ell}^{1,\delta}<E_2^\ast$ and we are thus in the elliptic region.
Therefore
\begin{align}
    \sum_{0<\ell<\ell_{m}^{\ast}}\frac{1}{\ell^2}I_{m,\ell}
    &= \sum_{0<\ell<\ell_{m}^{\ast}}q^2\int_{E_{\ell}^{1}}^{E_{\ell}^{1,\de}}\frac{\log\left(\frac{T(E)}{\ell q}-1\right)^2|\mu'(\k^2E)|T(E)^2}{T'(E)^{1-\alpha}\tE^\beta} \diff E\notag \\
    & \le C \sum_{0<\ell<\ell_{m}^{\ast}}q^2\int_{E_{\ell}^{1}}^{E_{\ell}^{1,\de}}\frac{\log\left(\frac{T(E)}{\ell q}-1\right)^2}{(E-\Emin)^\beta} \diff E \notag\\
    & \le C \sum_{0<\ell<\ell_{m}^{\ast}}q^2\int_{E_{\ell}^{1}}^{E_{\ell}^{1,\de}}\frac{\log\left(\frac{T(E)}{\ell q}-1\right)^2}{(T(E)-T(\Emin))^\beta} \diff E,
    \end{align}
    where in the last line we used the bound $\tE\gtrsim T(E)-T_{\min}$ for $E\in \Iell$.
     We now use the change of variables~\eqref{E:CH} and the above bound to obtain
  \begin{align}
    \sum_{0<\ell<\ell_{m}^{\ast}}\frac{1}{\ell^2}I_{m,\ell}
    &\leq C\sum_{0<\ell<\ell_{m}^{\ast}} q^2(\ell q)^{1-\beta}\int_{\max\{1,\frac{\Tmin}{\ell q}\}}^{1+\de}\frac{\log\left(t-1\right)^2}{(t-\frac{\Tmin}{\ell q})^\beta} \diff t \notag\\
    &\le C\sum_{0<\ell<\ell_{m}^{\ast}} \frac1{\ell^2} \int_{\max\{1,\frac{\Tmin}{\ell q}\}}^{1+\de}\frac{\log\left(t-1\right)^2}{(t-\frac{\Tmin}{\ell q})^\beta} \diff t \le C, \label{E:AL1BOUND2}
    \end{align}
    where we have used the bound $\ell q \le T(E^\ast)+1$ for all $\ell\in\{1,\dots,\ell^\ast\}$ and the assumption $\beta<1$.
Thus \eqref{E:AL1BOUND} with \eqref{E:AL1BOUND2} together imply that
\begin{equation*}
    \sum_{\ell>0}A_{\ell,1}\leq \|\partial_x\varphi_g\|_{L^2}^2
\end{equation*}
Together with bounds~\eqref{E:ABOUND},~\eqref{E:ABOUNDELL}, and~\eqref{E:AELLBOUND} we finally conclude that
\begin{equation}
    \sum_{\ell>0}|A_\ell|\leq C\|\partial_x\varphi_g\|_{L^2}^2 \label{E:ALFINAL}
\end{equation}


{\bf Estimates on $B_{\ell}$.} Recalling~\eqref{E:BELLDEF} we have 
\begin{align}
  |B_\ell| & 
    \leq\int_{E_{\ell}^{1}}^{E_{\ell}^{1,\de}}\left|\frac{\diff}{\diff E}\left(\frac{T(E)|\mu'(\k^2E)|}{T'(E)}\right)\right||\widehat{\varphi}_g(\ell, E)|^2|p_{m,\ell}(T(E))|\diff E \notag \\
    &\leq C\int_{E_{\ell}^{1}}^{E_{\ell}^{1,\de}}\left(\left|\frac{T'(E)^2-T(E)T''(E)}{T'(E)^2}\right||\mu'(\eps^2 E)|+\k^2\left|\frac{T(E)}{T'(E)}\right||\mu''(\eps^2E)|\right)|\widehat{\varphi}_g(\ell, E)|^2|p_{m,\ell}(T(E))|\diff E \notag \\
    &\le C(B_{\ell,1}+\k^2 B_{\ell,2}), \label{E:BLDECOMP}
\end{align}
where
\begin{align}  
B_{\ell,1}&:=\int_{E_{\ell}^{1}}^{E_{\ell}^{1,\de}}\left|\frac{T'(E)^2-T(E)T''(E)}{T'(E)^2}\right||\widehat{\varphi}_g(\ell, E)|^2|p_{m,\ell}(T(E))|\diff E, \label{E:BELL1DEF}\\   
B_{\ell,2}&:=\int_{E_{\ell}^{1}}^{E_{\ell}^{1,\de}}\left|\frac{T(E)}{T'(E)}\right||\widehat{\varphi}_g(\ell, E)|^2|p_{m,\ell}(T(E))|\diff E,
\end{align}
and we have used~\eqref{E:MUUNIF}.
Arguing analogously as in the bounds for $A_{\ell}$, the Cauchy-Schwarz inequality gives
\begin{align}
    |B_{\ell,2}|
    &\leq C\left(\int_{E_{\ell}^{1}}^{E_{\ell}^{1,\de}}\frac{p_{m,\ell}(T(E))^2|\widehat{\varphi}_g(\ell, E)|^2T(E)}{T'(E)^{1-\ga}}\diff E+\int_{E_{\ell}^{1}}^{E_{\ell}^{1,\de}}|\widehat{\varphi}_g(\ell, E)|^2\frac{T(E)}{T'(E)^{1+\ga}} \diff e\right), \label{E:BELL2}
\end{align}
for a parameter $\ga\in(0,1)$ to be determined.
We note that
\begin{align*}
    \int_{E_{\ell}^{1}}^{E_{\ell}^{1,\de}}|\widehat{\varphi}_g(\ell, E)|^2\frac{T(E)}{T'(E)^{1+\ga}} \diff E&\leq\frac{1}{4\pi^2\ell^2}\int_{E_{\ell}^{1}}^{E_{\ell}^{1,\de}}\int_{\mathbb{S}^1}|\partial_\theta\varphi_g(\theta,E)|^2\frac{T(E)}{T'(E)^{1+\ga}}\diff\theta\diff E\\
    &\leq\frac{C}{\ell^2}\int_{x_{-}(E_{\ell}^{1,\de})}^{x_{+}(E_{\ell}^{1,\de})}|\partial_x\varphi_g|^2\int_{E_{\ell}^{1,\delta}}^{E_{\ell}^1}\frac{T(E)^2}{T'(E)^{1+\ga}}\diff E\diff x\\
    &\leq \frac{C}{\ell^2}\|\partial_x\varphi_g\|_{L^2}^2\int_{E_{\ell}^{1}}^{E_{\ell}^{1,\de}}\frac{T(E)^2}{T'(E)^{1+\ga}}\diff E,
\end{align*}
where we changed variables $\theta\mapsto x$ and noticed that $\frac{\partial\varphi_g}{\partial\theta}=\frac{\partial\varphi_g}{\partial x}\frac{\partial x}{\partial\theta}$, together
with  $\frac{\pa x}{\pa\th}=\pm T(E)\sqrt{E+\varphi_0(x)}$.
In the hyperbolic range $\ell\geq\ell^\ast$ we use the change of variables~\eqref{E:CH} and~\eqref{E:KEYHYP} to get
\begin{align*}
    \int_{E_{\ell}^{1}}^{E_{\ell}^{1,\de}}\frac{T(E)^2}{T'(E)^{1+\ga}}\diff E&\le C\ell^3 q^3\int_{1}^{1+\de}t^2e^{-(2+\ga)\ell qt}\diff t
    \leq C\ell^3 q^3e^{-(2+\ga)\ell q} .
\end{align*}
The two previous bounds then give
\begin{align}   \sum_{\ell\geq\ell^\ast}\int_{E_{\ell}^{1}}^{E_{\ell}^{1,\de}}|\widehat{\varphi}_g(\ell, E)|^2\frac{T(E)}{T'(E)^{1+\ga}} \diff E
&\leq Cq^3\|\partial_x\varphi_g\|^2_{L^2}\sum_{\ell>0}\ell e^{-\kappa_\ga\ell q}  \le C\|\partial_x\varphi_g\|^2_{L^2}, \label{E:EXPARG}
\end{align}
where $\kappa_\ga:=2+\ga$ and we have used the integral criterion in the very last line by analogy to~\eqref{E:IC1}.
Note further that the first term on the right-hand side of~\eqref{E:BELL2} is in fact $A_{\ell,1}$~\eqref{E:AELL1}, so by~\eqref{E:AL1BOUND}
we conclude
\begin{equation}\label{E:BL2HYP}
    \sum_{\ell\geq\ell^\ast}|B_{\ell,2}|\leq C\|\partial_x\varphi_g\|^2_{L^2}.
\end{equation}

To estimate $B_{\ell,1}$~\eqref{E:BELL1DEF}  
we first note that by Theorem~\ref{T:PERIODASYMP} 
\begin{align}
\left|\frac{T'(E)^2-T(E)T''(E)}{T'(E)^2}\right|&\leq C |\log |E|| \chh(E) + C \che(E) \notag\\
& \le C T(E) \chh(E) + C\che(E).\label{E:USE}
\end{align}
Therefore, in the hyperbolic range $\ell\geq\ell^\ast$
\begin{align}  
    |B_{\ell,1}|&\leq\frac{C}{\ell^2}\|\partial_x\varphi_g\|^2_{L^2}\int_{E_{\ell}^{1}}^{E_{\ell}^{1,\de}}T(E)|p_{m,\ell}(T(E))|\diff E \notag \\
    &=\frac{C}{\ell^2}\|\partial_x\varphi_g\|^2_{L^2}\int_{E_{\ell}^{1}}^{E_{\ell}^{1,\de}}T(E)\ell q\left|\log\left(\frac{T(E)}{\ell q}-1\right)\right|\diff E \notag\\
   &\le  \frac{C}{\ell^2}\|\partial_x\varphi_g\|^2_{L^2}\ell^3 q^3\int_{1}^{1+\de}t\left|\log\left(t-1\right)\right|e^{-\ell q t}\diff t\notag\\
    &\leq C\ell q^3e^{-\ell q}\|\partial_x\varphi_g\|^2_{L^2},
\end{align}
where we have used the change of variables~\eqref{E:CH} and~\eqref{E:KEYHYP}.
Applying the integral criterion as before, we conclude
\begin{align}
    \sum_{\ell\geq\ell^\ast}|B_{\ell,1}|
    &\leq Cq^3\|\partial_x\varphi_g\|^2_{L^2}\sum_{\ell>0}\ell e^{-\ell q}  
    \leq Cq^3\|\partial_x\varphi_g\|^2_{L^2}\int_{1}^\infty xe^{-q x}dx \notag \\
    &=Cq\|\partial_x\varphi_g\|^2_{L^2}\int_{ q}^\infty xe^{-x}\diff x \le C\|\partial_x\varphi_g\|^2_{L^2}. \label{E:BL1HYP}
\end{align}

In the elliptic range $1\le \ell < \ell^\ast$, using~\eqref{E:USE}, we have by the same argument
\begin{align}
\sum_{\ell=1}^{\ell^\ast} |B_{\ell,1}| &\le \frac{C}{\ell^2}\|\partial_x\varphi_g\|^2_{L^2}\int_{E_{\ell}^{1}}^{E_{\ell}^{1,\de}}T(E)|p_{m,\ell}(T(E))|\diff E \notag\\
& \le \sum_{\ell=1}^{\ell^\ast}\frac{C}{\ell^2}\|\partial_x\varphi_g\|^2_{L^2} \ell^3 q^3 \int_{\max\{1, \frac{\Tmin}{\ell q}\}}^{1+\de} t |\log(t-1)|^2 \diff t \notag\\
& \le C \|\partial_x\varphi_g\|^2_{L^2}, \label{E:BL1ELL}
\end{align}
where we have used $\ell q\lesssim 1$ in the elliptic range. An analogous bound gives
\begin{align}
\sum_{\ell=1}^{\ell^\ast} |B_{\ell,2}|\le C \|\partial_x\varphi_g\|^2_2. \label{E:BL2ELL}
\end{align}
Combining~\eqref{E:BLDECOMP},~\eqref{E:BL1HYP},~\eqref{E:BL2HYP},~\eqref{E:BL1ELL},~\eqref{E:BL2ELL}, we conclude
\begin{equation}
    \sum_{\ell>0}|B_\ell|\leq C\|\partial_x\varphi_g\|_{L^2}^2. \label{E:BLFINAL}
\end{equation}
From~\eqref{E:JTRAPBOUND},~\eqref{E:ALFINAL}, and~\eqref{E:BLFINAL} we finally obtain~\eqref{E:TRAPBOUND}.
\end{proof}

\subsubsection{Estimates outside the separatrix: $\mathcal J_{\text{ext}}$-term.}


The main result of this section is the following lemma, which treats the particles located outside the separatrix.


\begin{lemma}\label{L:EXTLEMMA}
Under the assumptions of Lemma~\ref{L:ENERGYBASIC}, there exists a constant $C>0$ such that
\begin{align}\label{E:EXTBOUND}
\mathcal J_{\text{ext}} \le C \eps^{-2} \|\partial_x\varphi_g\|_{L^2}^2,
\end{align}
where we recall the definition~\eqref{E:JEXTDEF} of $\mathcal J_{\text{ext}}$.
\end{lemma}


\begin{proof}
We recall from Section~\ref{S:FS} and~\eqref{E:BIGLMSTARDEF} that for frequencies $\ell\ge L^m_\ast$ we are in the hyperbolic zone, i.e.
$[E^{2,\de}_{m,\ell},E^{2}_{m,\ell}]\subset \Ihypall\cap I_2=(0,|E^\ast|)$. 
If $\mu$ is regular (i.e. if in addition to~\eqref{E:MUMINUSPLUSREG} $\mu\in C^2$ in a neighbourhood of $E=0$ as well)
in this range we may apply the exact same reasoning as in the hyperbolic zone inside the trapped region (i.e. the range $\ell\ge\ell^\ast$ which implies $[E^{1}_{m,\ell},E^{1,\de}_{\ell,m}]\subset \Ihypall\cap I_1=(E^\ast,0)$). These bounds give by analogy to~\eqref{E:ABOUND} and~\eqref{E:AL1BOUND}
\begin{align}
\sum_{\ell\ge L^\ast} \int_{E_{\ell}^{2,\de}}^{E_{\ell}^{2}}\frac{T(E)|\mu'(\k^2E)|}{\frac{T(E)}{\ell q}-1}\big(|\widehat{\varphi}_g(\ell, E)|^2+|\widehat{\varphi}_g(-\ell, E)|^2\big)\diff E
\le C \|\pa_x\varphi_g\|_{L^2}^2.
\end{align}

However, to include the Schamel-like steady states satisfying merely~\eqref{E:MUBLOWUP}, we must explain the arguments in the hyperbolic zone more carefully.
Without loss of generality, we provide the bounds only for the integrals involving $|\widehat{\varphi}_g(\ell, E)|^2$, as the integrals containing $|\widehat{\varphi}_g(-\ell, E)|^2$ are estimated identically. 

We proceed with the decomposition as in~\eqref{E:TRAP1}, where we integrate by parts to remove the singularity in the denominator. We obtain
\begin{align}
 &\sum_{\ell>0}\int_{E_{\ell}^{2,\de}}^{E_{\ell}^{2}}\frac{T(E)|\mu'(\k^2E)|}{\frac{T(E)}{\ell q}-1}|\widehat{\varphi}_g(\ell, E)|^2\diff E  = \sum_{\ell>0} \frac{T(E)p_{m,\ell}(T(E))}{T'(E)}|\mu'(\k^2E)||\widehat{\varphi}_g(\ell, E)|^2\Bigg|_{E_{\ell}^{2,\de}}^{E_{\ell}^{2}} \notag \\
    & \qquad - \sum_{\ell>0}\int_{E_{\ell}^{2,\de}}^{E_{\ell}^{2}}\frac{\diff}{\diff E}\left(\frac{T(E)|\mu'(\k^2E)||\widehat{\varphi}_g(\ell, E)|^2}{T'(E)}\right)p_{m,\ell}(T(E))\diff E, \label{E:TRAP1NEW}
\end{align}
where we recall~\eqref{E:PMLDEF}.
Due to the positivity of $\pml(t)$ for $t\in(\ell q,\infty)$ we see that the upper boundary of the boundary-term contribution is simply negative for $\delta<1$, since $T(E_{\ell}^{2,\de})=(1+\de)\ell q$ and $T'(E)<0$ for $E>0$. On the other hand
\begin{align*}
\lim_{E\to E_{\ell}^2}\frac{T(E)p_{m,\ell}(T(E))}{T'(E)}|\mu'(\k^2E)||\widehat{\varphi}_g(\ell, E)|^2 = 0,
\end{align*}
which follows from an argument analogous to~\eqref{E:BDRY1}.
Therefore, it follows from~\eqref{E:TRAP1NEW} that 
\begin{align}
 \sum_{\ell>0}\int_{E_{\ell}^{2,\de}}^{E_{\ell}^{2}}\frac{T(E)|\mu'(\k^2E)|}{1- \frac{T(E)}{\ell q}}|\widehat{\varphi}_g(\ell, E)|^2\diff E  
& \leq
    -\sum_{\ell>0}\int_{E_{\ell}^{2,\de}}^{E_{\ell}^{2}}\frac{\diff}{\diff E}\left(\frac{T(E)|\mu'(\k^2E)||\widehat{\varphi}_g(\ell, E)|^2}{T'(E)}\right)p_{m,\ell}(T(E))\diff E \notag \\
    & =\sum_{\ell>0}(A_\ell^{\text{ext}}+B_\ell^{\text{ext}}), \label{E:JTRAPBOUNDNEW}
\end{align}
where for any $\ell>0$
\begin{align}
    &A_\ell^{\text{ext}}=-2\int_{E_{\ell}^{2,\de}}^{E_{\ell}^{2}}\frac{T(E)|\mu'(\k^2E)|}{T'(E)}\text{Re}\left(\overline{\widehat{\varphi}_g(\ell, E)}\frac{\diff}{\diff E}\widehat{\varphi}_g(\ell, E)\right)p_{m,\ell}(T(E))\diff E, \label{E:AELLDEFEXT}\\
    &B_\ell^{\text{ext}}=-\int_{E_{\ell}^{2,\de}}^{E_{\ell}^{2}}\frac{\diff}{\diff E}\left(\frac{T(E)|\mu'(\k^2E)|}{T'(E)}\right)|\widehat{\varphi}_g(\ell, E)|^2p_{m,\ell}(T(E))\diff E. \label{E:BELLDEFEXT}
\end{align}

Observe next that
\begin{align*}
    \left|\frac{\diff}{\diff E}\left(\frac{T(E)|\mu'(\k^2E)|}{T'(E)}\right)\right|&\leq\left|\frac{T'(E)^2-T(E)T''(E)}{T'(E)^2}\right||\mu'(\k^2E)|+\k^2\left|\frac{T(E)}{T'(E)}\right|\left|\mu''(\k^2E)\right|.
\end{align*}
Therefore, using~\eqref{E:MUBLOWUP},~\eqref{E:MUTAIL}, and Theorem~\ref{T:PERIODASYMP} we obtain
\begin{align}
\left|\frac{\diff}{\diff E}\left(\frac{T(E)|\mu'(\k^2E)|}{T'(E)}\right)\right|
\le C\frac{|\log E|}{\eps\sqrt{E}} \chi_{0<E\le |E^\ast|} + \frac C{\eps^2 E}\chi_{E\ge |E^\ast|}. \label{E:EXTCRUCIAL}
\end{align}

{\em Exterior zone: $0<\ell<L^\ast$, i.e. $E\in \Iext$.}
So away from the hyperbolic zone, i.e. for indices $\ell\in\{1,\dots L^\ast-1\}$ we  
have
\begin{align}
    \sum_{0<\ell<L^\ast}B_{\ell}^{\text{ext}}
    &\leq  C\eps^{-2}\|\partial_x\varphi_g\|^2_{L^2}\sum_{0<\ell<L^\ast}\frac{1}{\ell^2}\int_{E_{\ell}^{2}}^{E_{\ell}^{2,\delta}}\frac1{E}|p_{m,\ell}(T(E))|\diff E. \notag
\end{align}
We now use the change of variables~\eqref{E:CH}
and keep in mind that for $E\in\Iext\subset\mathbb R^+$ we have $E^{\frac12}\lesssim T(E)\lesssim E^{\frac12}$ and thus $\frac1{E}\lesssim \frac1{T(E)^2}$, $E>|E_\ast|$.
Keeping in mind that by Theorem~\ref{T:PERIODASYMP} $T(E)^{-3} \lesssim |T'(E)| \lesssim T(E)^{-3}$, 
we have
\begin{align*}
   & C\|\partial_x\varphi_g\|^2_{L^2}\eps^{-2}\sum_{0<\ell<L^\ast}\frac{1}{\ell^2}\int_{E_{\ell}^{2}}^{E_{\ell}^{2,\delta}}\frac1{E}|p_{m,\ell}(T(E))|\diff E.\notag\\
    &\leq C\|\partial_x\varphi_g\|_{L^2}^2\k^{-2}\sum_{0<\ell<L^\ast}\frac{1}{\ell^2}\int_{E_{\ell}^{2}}^{E_{\ell}^{2,\delta}}T(E)^{-2}|p_{m,\ell}(T(E))|\diff E\\
    &\leq C\|\partial_x\varphi_g\|_{L^2}^2\k^{-2}\sum_{0<\ell<L^\ast}\frac{1}{\ell^2}\int_{1-\delta}^{1} (t\ell q)^{-2} \ell q |\log(1-t)| \ell q (t\ell q)^3 \diff t\\
    &=  C\k^{-2}\|\partial_x\varphi_g\|_{L^2}^2\sum_{0<\ell<L^\ast}\ell q^3\int_{1-\delta}^{1}t |\log(1-t)|\diff t \\
    &\leq C\k^{-2}\|\partial_x\varphi_g\|_{L^2}^2\sum_{0<\ell<L^\ast}\frac1{\ell^2}\\
    &\leq C\k^{-2}\|\partial_x\varphi_g\|_{L^2}^2,
\end{align*}
where we have used the bound $\ell q \le C$ for all $0<\ell<L^\ast$ (see~\eqref{E:BIGLMSTARDEF}). 
As in~\eqref{E:AELLBOUND}--\eqref{E:AELL2TRAP} we apply Cauchy-Schwarz to obtain
\begin{align}
\big|A_\ell^{\text{ext}}\big|
\le A_{\ell,1}^{\text{ext}} + A_{\ell,2}^{\text{ext}}, 
\end{align}
where for $\ell\in\{1,\dots,L^\ast-1\}$ we let
\begin{align}
A_{\ell,1}^{\text{ext}} &= \int_{E_{\ell}^{2}}^{E_{\ell}^{2,\delta}}\frac{p_{m,\ell}(T(E))^2|\widehat{\varphi}_g(\ell, E)|^2T(E)|\mu'(\k^2E)|}{|T'(E)|}\diff E \label{E:AELLONEEXT}\\
A_{\ell,2}^{\text{ext}} &=\int_{E_{\ell}^{2}}^{E_{\ell}^{2,\delta}}\frac{\left|\frac{\diff}{\diff E}\widehat{\varphi}_g(\ell, E)\right|^2 T(E)|\mu'(\k^2E)|}{|T'(E)|} \diff E. \label{E:AELL2}
\end{align}

By Theorem~\ref{T:DXDE} we have the bound $|\pa_E\varphi_g|= |\pa_x\varphi_g|\left|\frac{\pa x}{\pa E}\right| \le \frac C{|E|}|\pa_x\varphi_g|$, $E\in\Iext$. Therefore by Theorem~\ref{T:PERIODASYMP}, Plancharel, the change of variables $(\th,E)\mapsto (x,E)$, and~\eqref{E:MUBLOWUP}
\begin{align}
    \sum_{0<\ell < L^\ast}|A^{\text{ext}}_{\ell,2}|
&\leq C \int_{0}^{1}\int_{\frac{|E^\ast|}{2}}^{\infty} |\pa_x\varphi_g|^2|E|^{-2}  |\mu'(\k^2E)| |E|^{\frac32} (E+\varphi_0(x))^{-\frac12}\diff ( E,x) \notag \\
&\leq C \|\pa_x\varphi_g\|_{L^2}^2 \eps^{-2}\int_{\frac{|E^\ast|}{2}}^\infty|E|^{-2}\diff E \notag\\
& \le C \eps^{-2}\|\pa_x\varphi_g\|_{L^2}^2,
\end{align}
where we have used the assumption~\eqref{E:MUTAIL} to control $|\mu'(\eps^2 E)|$.
To control~\eqref{E:AELLONEEXT}, we proceed analogously to~\eqref{E:AELL21}. Using the bound~\eqref{E:SIMPLE},~\eqref{E:MUTAIL}, and Theorem~\ref{T:PERIODASYMP}, we conclude
\begin{align}
   \sum_{0<\ell < L^\ast}|A^{\text{ext}}_{\ell,1}|
   &\le C \|\pa_x\varphi_g\|_{L^2}^2 \sum_{0<\ell < L^\ast} \frac{1}{\ell^2} (\ell q)^2\int_{E_{\ell}^{2}}^{E_{\ell}^{2,\delta}} \frac{|\log(\frac{T(E)}{\ell q}-1)|^2|E|^{-\frac12}\eps^{-2}|E|^{-1}}{|E|^{-\frac32}} \diff E \notag\\
   &\le C\eps^{-2} \|\pa_x\varphi_g\|_{L^2}^2 \sum_{0<\ell < L^\ast} \frac{1}{\ell^2} (\ell q)^2\int_{E_{\ell}^{2}}^{E_{\ell}^{2,\delta}} |\log(\frac{T(E)}{\ell q}-1)|^2 \diff E \notag \\
   & \le  C\eps^{-2} \|\pa_x\varphi_g\|_{L^2}^2 \sum_{0<\ell < L^\ast} \frac{1}{\ell^2} (\ell q)^2 \int_1^{1+\de} |\log (t-1)|^2 (\ell q)^{-2} t^{-3} \diff t\notag \\
   & \le  C\eps^{-2} \|\pa_x\varphi_g\|_{L^2}^2 \sum_{0<\ell < L^\ast} \frac{1}{\ell^2} \notag\\
   & \le C\eps^{-2} \|\pa_x\varphi_g\|_{L^2}^2,
\end{align}
where we have used the change of variables~\eqref{E:CH} in the third line together with the observation $\diff E \approx (\ell q)^{-2} t^{-3}\diff t$, which follows from~\eqref{E:CH} and
the sharp asymptotics of Theorem~\ref{T:PERIODASYMP} in the exterior region $\Iext$.

{\em Hyperbolic zone: $\ell\ge L^\ast$, i.e. $E\in \Ihypall\cap(0,\infty)$.}
Observe that for $E\in (0,|E^\ast|)$, by analogy to~\eqref{E:KEYHYP} we have the bound
\begin{align}\label{E:KEYHYPEXT}
\exp(T(E)) \lesssim |T'(E)| \lesssim  \exp(T(E)), \ \ E\in(0,|E^\ast|).
\end{align}

Using~\eqref{E:BELLDEFEXT}--\eqref{E:EXTCRUCIAL}, the change of variables~\eqref{E:CH}, and Theorem~\ref{T:PERIODASYMP} we obtain
\begin{align}
    \sum_{\ell\ge L^\ast}B_{\ell}^{\text{ext}}
    &\leq  C\eps^{-1}\|\partial_x\varphi_g\|^2_{L^2}\sum_{\ell\ge L^\ast}\frac{1}{\ell^2}\int_{E_{\ell}^{2}}^{E_{\ell}^{2,\delta}}\frac{|\log E|}{\sqrt E}|p_{m,\ell}(T(E))|\diff E \notag \\
& \le     C\eps^{-1}\|\partial_x\varphi_g\|^2_{L^2}\sum_{\ell\ge L^\ast}\frac{1}{\ell^2}\int_{E_{\ell}^{2}}^{E_{\ell}^{2,\delta}}\frac{T(E)}{\exp(-\frac12 T(E))}|p_{m,\ell}(T(E))|\diff E \notag \\
& \le C\eps^{-1}\|\partial_x\varphi_g\|^2_{L^2}\sum_{\ell\ge L^\ast}\frac{1}{\ell^2}(\ell q)^3\int_{1}^{1+\de}\frac{t}{\exp(-\frac12 \ell q t)}|\log(t-1)| \exp(-\ell q t)\diff t \notag\\
& = C\eps^{-1}\|\partial_x\varphi_g\|^2_{L^2}\sum_{\ell\ge L^\ast}\frac{1}{\ell^2}(\ell q)^3\int_{1}^{1+\de}t|\log(t-1)| \exp(-\frac12\ell q t)\diff t \notag\\
& \le C\eps^{-1}\|\partial_x\varphi_g\|^2_{L^2}\sum_{\ell\ge L^\ast}\frac{1}{\ell^2}(\ell q)^3 e^{-\frac{1+\de}{2}\ell q} \le C\eps^{-1}\|\partial_x\varphi_g\|^2_{L^2}, \label{E:EXPEXT}
\end{align} 
where in the last line we have used an argument analogous to~\eqref{E:EXPARG}.

In the hyperbolic we need a refinement of the decomposition~\eqref{E:AELLONEEXT}--\eqref{E:AELL2}. We introduce a positive parameter $1<\nu<\frac32$. By 
 Cauchy-Schwarz we have
\begin{align}
\big|A_\ell^{\text{ext}}\big|
\le A_{\ell,1}^{\text{ext}} + A_{\ell,2}^{\text{ext}}, 
\end{align}
where for $\ell\ge L_\ast$ we let
\begin{align}
A_{\ell,1}^{\text{ext}} &= \int_{E_{\ell}^{2}}^{E_{\ell}^{2,\delta}}\frac{p_{m,\ell}(T(E))^2|\widehat{\varphi}_g(\ell, E)|^2T(E)|\mu'(\k^2E)|}{|T'(E)|^{1-\nu}}\diff E \label{E:AELLONEEXT}\\
A_{\ell,2}^{\text{ext}} &=\int_{E_{\ell}^{2}}^{E_{\ell}^{2,\delta}}\frac{\left|\frac{\diff}{\diff E}\widehat{\varphi}_g(\ell, E)\right|^2 T(E)|\mu'(\k^2E)|}{|T'(E)|^{1+\nu}} \diff E. \label{E:AELL2}
\end{align}

By Theorem~\ref{T:DXDE} we have the bound $|\pa_E\varphi_g|= |\pa_x\varphi_g|\left|\frac{\pa x}{\pa E}\right| \le \frac C{|E|}|\pa_x\varphi_g|$, $E>0$. Therefore by Theorem~\ref{T:PERIODASYMP}, Plancharel, the change of variables $(\th,E)\mapsto (x,E)$, and~\eqref{E:MUBLOWUP}, we have
\begin{align}
    \sum_{\ell \ge L^\ast}|A^{\text{ext}}_{\ell,2}|
&\leq C \int_{0}^{1}\int_{0}^{|E^\ast|} |\pa_x\varphi_g|^2|E|^{-2}  |\mu'(\k^2E)| |E|^{1+\nu}(E+\varphi_0(x))^{-\frac12}\diff ( E,x) \notag \\
&\leq C \|\pa_x\varphi_g\|_{L^2}^2 \eps^{-1}\int_{0}^{|E^\ast|} E^{\nu-2}\diff E \notag\\
& \le C \eps^{-1}\|\pa_x\varphi_g\|_{L^2}^2,
\end{align}
since $\nu>1$. 
On the other hand, note that $\int_{\mathbb T^1}|\pa_\th\varphi_g|^2 \diff\th \le \int_{\mathbb T^1} |\pa_x\varphi_g|^2|\pa_\th x| \diff x$. Therefore, 
\begin{align}
   \sum_{\ell \ge L^\ast}|A^{\text{ext}}_{\ell,1}|
   &\le C \|\pa_x\varphi_g\|_{L^2}^2 \sum_{\ell 
  \ge L^\ast} \frac{1}{\ell^2} (\ell q)^2\int_{E_{\ell}^{2}}^{E_{\ell}^{2,\delta}} \frac{|\log(\frac{T(E)}{\ell q}-1)|^2T(E)^2|\mu'(\eps^2E)|}{T'(E)^{1-\nu}} \diff E. \notag
  \end{align}
  We now use~\eqref{E:KEYHYPEXT}, the relation $\diff E =\frac{\ell q \diff t}{T'(E)}$, and $\frac1{\sqrt E}\lesssim |T'(E)|^{\frac12}$ to obtain
   \begin{align}
  \sum_{\ell \ge L^\ast}|A^{\text{ext}}_{\ell,1}|
 &\le C\eps^{-1} \|\pa_x\varphi_g\|_{L^2}^2 \sum_{\ell \ge L^\ast} \frac{1}{\ell^2} (\ell q)^2\int_{E_{\ell}^{2}}^{E_{\ell}^{2,\delta}} |\log(\frac{T(E)}{\ell q}-1)|^2 \frac{T(E)^2 }{T'(E)^{1-\nu}E^{\frac12}} \diff E \notag \\
 & \le  C\eps^{-1} \|\pa_x\varphi_g\|_{L^2}^2 \sum_{\ell \ge L^\ast}  \frac{1}{\ell^2} (\ell q)^5\int_{1}^{1+\de} |\log(t-1)|^2 t^2 \exp(-(\frac32-\nu) q\ell t) \diff t \notag\\
   & \le C\eps^{-1} \|\pa_x\varphi_g\|_{L^2}^2,
\end{align}
where we have again used the exponential decay as in~\eqref{E:EXPEXT} and the assumption $\nu<\frac32$.
\end{proof}

\begin{remark}
The above proof shows that the Schamel-type singularity~\eqref{E:MUBLOWUP} allows our techniques to carry through with some additional room reflected in the freedom to choose any $\nu\in(1,\frac32)$ in the above argument. In other words one can relax the assumption~\eqref{E:MUBLOWUP} to a more singular behaviour. We do not pursue this question here.
\end{remark}


\subsubsection{Non-resonant estimates: $\mathcal J_{\text{nr}}$-term}


The goal of this section is to deal with the non-resonant contribution $\mathcal J_{\text{nr}}$~\eqref{E:JNRDEF}, which is particularly easy to bound as
the denominators never get small.


\begin{lemma}\label{L:NRLEMMA}
Under the assumptions of Lemma~\ref{L:ENERGYBASIC}, there exists a constant $C>0$ such that
\begin{enumerate}
\item if $\mu$ is monotone we have the bound
\begin{equation}\label{E:NRBOUND}
\mathcal J_{\text{nr}} \le \frac{C}{\delta} \|\partial_x\varphi_g\|_{L^2}^2,
\end{equation}
where we recall the definition~\eqref{E:JNRDEF} of $\mathcal J_{\text{nr}}$;
\item if $\mu$ has a local maximum at $E=0$ we have the bound
\begin{equation}\label{E:NRBOUNDUP}
\mathcal J_{\text{nr}}^{\uparrow} \le \frac{C}{\delta} \|\partial_x\varphi_g\|_{L^2}^2,
\end{equation}
where we recall the definition~\eqref{E:JNRDEFUP} of $\mathcal J_{\text{nr}}^{\uparrow}$.
\end{enumerate}
\end{lemma}


\begin{proof}
We only prove~\eqref{E:NRBOUND}, as~\eqref{E:NRBOUNDUP} follows analogously. Note first that 
\begin{equation}\label{E:SIMPLE}
    4\pi^2\ell^2|\widehat{\varphi}_g(\ell, E)|^2\leq\int_{\mathbb{S}^1}|\partial_\theta\varphi_g(\theta,E)|^2d\theta, \ \ (\ell,E)\in\Z_\ast\times I.
\end{equation}
Therefore, upon changing variables $\theta\mapsto x$ and using $\frac{\partial\varphi_g}{\partial\theta}=\frac{\partial\varphi_g}{\partial x}\frac{\partial x}{\partial\theta}$ we have
\begin{align*}
&\sum_{\ell\in\Z_\ast}\int_{E_{\ell}^{1,\delta}}^{0}\frac{T(E)|\mu'(\k^2E)|}{\frac{T(E)}{\ell q}-1}|\widehat{\varphi}_g(\ell, E)|^2\diff e \\
&\leq\frac{1}{\delta}\sum_{\ell\in\Z_\ast}\int_{E_{\ell}^{1,\delta}}^{0}|\mu'(\eps^2 E)|T(E)|\widehat{\varphi}_g(\ell, E)|^2dE\\
&\leq\frac{1}{4\pi^2\delta}\sum_{\ell\in\Z_\ast}\frac{1}{\ell^2}\int_{E_{\ell}^{1,\delta}}^{0}\int_{\mathbb{S}^1}|\partial_\theta\varphi_g(\theta,E)|^2|\mu'(\eps^2 E)|T(E)\diff\theta \diff E\\
&= \frac{1}{4\pi^2\delta}\sum_{\ell\in\Z_\ast}\frac{1}{\ell^2}\int_{0}^{1}|\partial_x\varphi_g|^2\int_{E_{\min}}^{0}T(E)^2|\mu'(\eps^2 E)|\sqrt{2E+2\varphi_0(x)}\diff E\diff x\\
&\leq \frac{C}{\delta}\sum_{\ell>0}\frac{1}{\ell^2}\int_{0}^{1}|\partial_x\varphi_g|^2\int_{E_{\min}}^{0}T(E)^2\diff E\diff x\\
&=\frac{C}{\delta}\|\partial_x\varphi_g\|_{L^2}^2\int_{E_{\min}}^{0}T(E)^2\diff E \\
& \le \frac{C}{\delta}\|\partial_x\varphi_g\|_{L^2}^2,
\end{align*}
where we have used~\eqref{E:TASYMPBOUND} in the last line and the bound $|\mu'(\eps^2E)|\lesssim 1$ for $E\in I_1$. We proceed analogously with the second term on the right-hand side of~\eqref{E:JNRDEF} and obtain
\begin{align}
&2 \sum_{\ell>0}\int_{0}^{E_{\ell}^{2,\de}} \frac{T(E) |\mu'(\k^2E)|}{\frac{T(E)}{\ell q}-1}\big( |\widehat{\varphi}_g(\ell, E)|^2+|\widehat{\varphi}_g(-\ell, E)|^2 \big) \diff E \notag \\
& \le \frac4{4\pi^2\de} \sum_{\ell>0}\frac1{\ell^2}\int_{0}^{E_{\ell}^{2,\de}}\int_{\mathbb S^1}T(E) |\mu'(\k^2E)| |\pa_\th\varphi_g(\th,E)|^2\diff\th \diff E \notag\\
& = \frac1{\pi^2 \de}\sum_{\ell>0} \frac1{\ell^2} \int_{\mathbb T^1} |\pa_x\varphi_g|^2 \int_{0}^{E_{\ell}^{2,\de}}T(E)^2\sqrt{2E+2\varphi_0(x)} |\mu'(\k^2E)| \diff E \diff x.
\end{align}
We observe now, using Theorem~\ref{T:PERIODASYMP},
\begin{align}
&\int_{0}^{E_{\ell}^{2,\de}}T(E)^2\sqrt{2E+2\varphi_0(x)} |\mu'(\k^2E)| \diff E \notag\\
& \le \Big( \int_{0}^{1} + \int_1^\infty\Big) T(E)^2\sqrt{2E+2\varphi_0(x)} |\mu'(\k^2E)| \diff E \notag\\
& \le \frac{\sqrt{2+2\varphi_0(x)}}{\eps}\int_0^1 \frac{|\log E|^2}{\sqrt E} \diff E + C\int_{1}^\infty \frac1{E}\sqrt{2E+2|\Emin|} \frac1{\eps^2 E} \diff E \notag\\
& \le \frac C{\eps^2},
\end{align}
where in the next-to-last line above we used~\eqref{E:MUBLOWUP} to bound the argument of the $\int_0^1 \diff E$-integral, and~\eqref{E:MUTAIL} to bound the argument of the
$\int_1^\infty \diff E$-integral.
\end{proof}

\subsection{Proof of the main theorem and Corollary~\ref{C:LANDAU}} \label{S:PROOFS}


{\em Proof of Theorem~\ref{T:MAIN}.}
By Lemma~\ref{L:ESSENTIAL} it suffices to show that there are no non-zero pure imaginary eigenvalues of $\mathscr L\Big|_{\overline{\text{ran}\,\D}}$.  Assume by contradiction that there exists such an eigenvalue of the form $\l=\frac{2\pi i }{q}$, where without loss we assume that $q>0$. Denote the associated eigenfunction by $g\in\overline{\text{ran}\D}\cap D(\D)$, $g\neq0$. Then by Lemmas~\ref{L:ENERGYBASIC},~\ref{L:NRLEMMA},~\ref{L:TRAPLEMMA}, and~\ref{L:EXTLEMMA}, we conclude that there exists a $q$-independent constant $C>0$ such that
\begin{align}
  \|\partial_x\varphi_g\|_{L^2}^2 \le C_\ast \eps   \|\partial_x\varphi_g\|_{L^2}^2.
\end{align}
If $\eps<\frac1{C_\ast}$ it follows that necessarily $\pa_x\varphi_g=0$. Since $\l g=\mathscr L g = \D g -\k^2 |\mu'(\k^2 E)| v\pa_x \varphi_g =\D g$,
it follows that $\l\neq0$ is in the pure point spectrum of $\D$, a contradiction to the fact that $0$ is the only eigenvalue of $\D$.

\bigskip

{\em Proof of Corollary~\ref{C:LANDAU}.}
We recall from~\eqref{E:VLASOVLIN}--\eqref{E:POISSONLIN},~\eqref{E:DL}
that the linearised VP-system can be written in the form
\begin{align}\label{E:LINEQN}
\pa_tf = -\mathscr{L}f
\end{align}
where $\mathscr{L}:=\mathcal{D}\circ L$, $\mathcal{D}$, with the operators $\D$ and $L$ defined in~\eqref{E:DDEF} and~\eqref{E:LDEF} respectively.
The operator $\mathscr{L}$ is densely defined
and since it is skew-adjoint with respect to the inner product $(\cdot,\cdot)_L$, Stone's Theorem implies that $e^{-t\mathscr L}$ generates a $C^0$-semigroup in $\H$ so that for any $f_{in}\in D(\D)$ the unique solution
of~\eqref{E:LINEQN} is given by $e^{-t\mathscr L}f_0$. 
We define the operator $K:\H\to \H$
\be\label{E:KDEF}
Kf(x,v):=|\mu'(\k^2E(x,v))|^{\frac12}\pa_x\varphi_f(x), \ \ (x,v)\in[0,1]\times\mathbb R,
\ee
and claim that $K$ is indeed compact. To see this, first note that from the definition of the inner product~\eqref{E:INNERPRODL} and from Lemma~\ref{L:FORCEBOUND} the norm $\|\cdot\|_{L}$ induced by~\eqref{E:INNERPRODL} is equivalent to $\|\cdot\|_{\H}$. So we need to prove compactness of $K$ with respect to $\|\cdot\|_\H$. Let $(f_n)_{n\in\mathbb N}\subset \H$ be a bounded sequence so that $\|f_n\|_\H\le C$, $n\in\mathbb N$. 
Denote $\varphi_{f_n}$ by $\varphi_n$. Then for any $n$,
\begin{align}
|\pa_x\varphi_n(y)-\pa_x\varphi_n(z)| & = |\int_z^y\pa_{xx}\varphi_n(\tilde x)\diff\tilde x| \le \eps \Big( \int_z^y\int |\mu'(\k^2E(\tilde x,v))| \diff v \diff \tilde x\Big)^{\frac12} \|f_n\|_{\H}.
\end{align}
By assumption~\eqref{E:MUTAIL} we estimate
\begin{align}
\int |\mu'(\k^2E(\tilde x,v))| \diff v \le C \int \frac1{1 + \k^2|\frac12 v^2 - \varphi_0(\tilde x)|} \diff v \le C,
\end{align}
where we recall that $0\le \varphi_0 \le \varphi_0(x_0)$. Therefore by the two previous estimates we conclude
\begin{align}
|\pa_x\varphi_n(y)-\pa_x\varphi_n(x)| & \le C \eps|y-z|^{\frac12}.
\end{align}
This implies the equicontinuity of the sequence $\{\pa_x\varphi_n\}$. Testing $\pa_{xx}\varphi_n=\eps\int f_n\diff v$ against $x$ and integrating-by-parts over $x\in[0,1]$, we infer the identity $\pa_x\varphi_f(0)=\eps\int_{[0,1]\times\mathbb R}f(x,v)x\diff(x,v)$. Therefore, upon expressing
\[
\pa_x\varphi_n(x) = \eps\int_{[0,1]\times\mathbb R}f_n(x,v)x\diff(x,v)+\int_0^x\int_{\mathbb R} f_n(x,v)\diff x \diff v, \ \ x\in[0,1],
\]
it is easy to see that $\|\pa_x\varphi_n\|_{C^0([0,1])}$ is uniformly bounded. Using Arzela-Ascoli we pass to a uniform limit $\Phi=\lim_{k\to\infty}\pa_x\varphi_{n_k}$ along a subsequence $n_k$. It is then easily checked that $\||\mu'(\k^2E)|^{\frac12}\pa_x\varphi_{n_k}-|\mu'(\k^2E)|^{\frac12}\Phi\|_\H$ converges to $0$. Therefore the operator $K$ is compact.
 
Since by Theorem~\ref{T:MAIN} the point spectrum of $\mathscr L\Big|_{\rD}$ is empty and $\rD$ is an invariant subspace for the dynamics generated by~\eqref{E:LINEQN} , by the RAGE theorem~\cite{CyFrKiSi} we conclude
\begin{align}\label{E:RAGETH}
	0=\lim_{T\rightarrow 0}\frac{1}{T}\int_0^T\|Ke^{t\mathscr L} P_c(\mathscr L)f_{in}\|_{L}^2\diff t&\geq C\lim_{T\rightarrow 0}\frac{1}{T}\int_0^T\|Ke^{t\mathscr L} P_c(\mathscr L)f_{in}\|_{H}^2\diff t\\
&=C\lim_{T\rightarrow 0}\frac{1}{T}\int_0^T\|\pa_x\varphi_{f(t,\cdot)}\|_{L^2}^2\diff t,
\end{align}
where the inequality comes from the equivalence between $\|\cdot\|_L$ and $\|\cdot\|_\H$.


\appendix


\section{Proof of Theorem~\ref{T:PERIODASYMP}}\label{A:PERIOD}


\begin{lemma}
\label{lema:turn_points}
Assume that the assumptions $(\varphi1)$--$(\varphi4)$ are valid and let $x_{\pm}(\cdot)$ be the turning points of the characteristic flow defined through~\eqref{E:TURPOINTS}. Then the functions $x_\pm:I_1\rightarrow (0,1)$ are continuous and can be extended continuously to $\overline{I_1}=[\Emin,0]$. Moreover
the left turning point $x_{-}(E)$ satisfies
\begin{equation}\label{E:LEFTASY}
    x_{-}(E)=\sqrt{\frac{2}{\varphi_0''(0)}}|E|^{\frac{1}{2}}\big(1+o_{E\to0}(1)\big).
\end{equation}
Similarly, the right turning point $x_{+}(E)$ satisfies
\begin{equation}\label{E:RIGHTASY}
    x_{+}(E)=1-\sqrt{\frac{2}{\varphi_0''(1)}}|E|^{\frac{1}{2}}\big(1+o_{E\to0}(1)\big)
\end{equation}
and
\begin{equation}\label{E:EMINLIM}
	\lim_{E\rightarrow \Emin^+} x_{\pm}(E) = x_0
\end{equation}
\end{lemma}


\begin{proof}
The continuity of $x_{\pm}$ and \eqref{E:EMINLIM} follow from the regularity of $\varphi_0$ and the fact that $\varphi_0|_{[0,x_0]}$ and $\varphi_0|_{[x_0,1]}$ are monotone and continuous. The proof of the asymptotic formulas  \eqref{E:LEFTASY} and \eqref{E:RIGHTASY} follows from the second order Taylor expansion of $\varphi_0$ around $0$ and $1$ and $\varphi_0(x_{-}(E))=\varphi_0(x_{+}(E))=|E|$.
\end{proof}


\begin{lemma}[Sharp bounds on the period function]\label{L:ASYM_PER}
Let $\varphi_0$ satisfy the assumptions $(\varphi1)$--$(\varphi4)$. Then
there exists $\Emin<E^\ast<0$ satisfying~\eqref{E:PROP_E} and a constant $C>0$  such that for any $E\in I$ the following bounds hold:
\begin{align}
        \frac{1}{C} \Big(\big|\log|E|\big| \chh +  \che + |E|^{-\frac12} \chext \Big)\leq T(E)\leq C\Big( \big|\log|E|\big|\chh+  \che + |E|^{-\frac12}\chext\Big),\label{E:TASYMPBOUNDAPP}
\end{align}
where we recall~\eqref{E:CHISDEF} and the period function~\eqref{E:PERIODDEF}.
\end{lemma}

\begin{proof}
{\em Step 1. Elliptic region: $E\in \Iell\subset I_1$.} This result is well known and follows verbatim from Appendix B of~\cite{HRS2021}.

{\em Step 2. Hyperbolic region: $E\in \Ihypall$.}
Recall that $\Emin<E^\ast<0$ and  definition~\eqref{E:IHYPDEF}.
For any $E\in(E^\ast,0)\subset I_1$ we introduce the decomposition
\begin{equation*}
	T(E)=T_{-}(E)+T_{+}(E)+T_{\ast}(E),
\end{equation*}
where
\begin{align*}
	T_{-}(E)&:=2\int_{x_{-}(E)}^{x_{-}(E^\ast)}\frac{1}{\sqrt{2E+2\varphi_0(x)}}\diff x,\\
	T_{+}(E)&:=2\int_{x_{+}(E^\ast)}^{x_{+}(E)}\frac{1}{\sqrt{2E+2\varphi_0(x)}}\diff x,\\
	T_{\ast}(E)&:=2\int_{x_{-}(E^\ast)}^{x_{+}(E^\ast)}\frac{1}{\sqrt{2E+2\varphi_0(x)}}\diff x.
\end{align*}
We Taylor expand around $x_{-}(E)$ to obtain
\begin{equation*}
\varphi_0(x)=\varphi_0(x_{-}(E))+\varphi_0'(x_{-}(E))(x-x_{-}(E))+\frac{\varphi_0''(\xi(x,E))}{2}(x-x_{-}(E))^2
\end{equation*}
for some $\xi(x,E)\in[x_{-}(E),x]\subseteq [0,x_{-}(E^\ast)]$. By our choice of $E^\ast$, there exist $C_1,C_2>0$ constants (which depend on $E^\ast$ in general) such that
\begin{align}
	\varphi_0'(x_{-}(E))(x-x_{-}(E))+C_1(x-x_{-}(E))^2\notag &\geq\varphi_0(x)-\varphi_0(x_{-}(E)) \notag \\
	&\geq\varphi_0'(x_{-}(E))(x-x_{-}(E))+C_2(x-x_{-}(E))^2. \label{E:TAYLORIMP}
\end{align}
Recalling the relation $\varphi_0(x_-(E))=-E$ we conclude
\begin{equation*}
	T_{-}(E)\geq C\int_{x_{-}(E)}^{x_{-}(E^\ast)}\frac{1}{\sqrt{\varphi_0'(x_{-}(E))(x-x_{-}(E))+C_1(x-x_{-}(E))^2}}\diff x
\end{equation*}
and
\begin{equation*}
T_{-}(E)\leq C\int_{x_{-}(E)}^{x_{-}(E^\ast)}\frac{1}{\sqrt{\varphi_0'(x_{-}(E))(x-x_{-}(E))+C_2(x-x_{-}(E))^2}}\diff x
\end{equation*}
We now change variables 
\begin{align}\label{E:CHVAR}
y=C_i\frac{x-x_-(E)}{\varphi_0'(x_{-}(E))}, \ \ i\in\{1,2\},
\end{align}
respectively, and obtain
\begin{equation*}
	\frac{C}{\sqrt{C_1}}\int_{0}^{\frac{C_1(x_{-}(E^\ast)-x_{-}(E))}{\varphi_0'(x_{-}(E))}}\frac{\diff y}{\sqrt{y+y^2}}\leq T_{-}(E)\leq\frac{C}{\sqrt{C_2}}\int_{0}^{\frac{C_2(x_{-}(E^\ast)-x_{-}(E))}{\varphi_0'(x_{-}(E))}}\frac{\diff y}{\sqrt{y+y^2}}.
\end{equation*}
We now observe the identity
\begin{equation}\label{E:INTID}
	\int_0^b\frac{1}{\sqrt{x+x^2}}\diff x = \frac{1}{2}\log\left(b\right)+\log\left(1+\sqrt{\frac{b+1}{b}}\right).
\end{equation}
We note the second term above converges to $\log(2)$ as $b\rightarrow+\infty$. Therefore, by Lemma~\ref{lema:turn_points}, there exists a constant $C>0$ (which depends on $E^\ast$ in general) such that
\begin{equation*}
	\frac{1}{C}|\log(|E|)|\leq T_{-}(E)\leq C|\log(|E|)|
\end{equation*} 
for every $E\in(E^\ast,0)$. A completely analogous argument yields
\begin{equation*}
	\frac{1}{C}|\log(|E|)|\leq T_{+}(E)\leq C|\log(|E|)|, \ \ E\in(E^\ast,0).
\end{equation*}
Since $T$ is increasing on $I_1$ it suffices to show that $T_\ast(E)$ is uniformly bounded from above on $(E^\ast,0)$.
To that end we change variables $-E'=\varphi_0(x)$ so that 
\begin{align*}
T_{\ast}(E)&=\frac{2}{\sqrt{2}} \int_{x_{-}(E^\ast)}^{x_0}\frac{1}{\sqrt{\varphi_0(x)+E}}\diff x+\frac{2}{\sqrt{2}} \int_{x_0}^{x_{+}(E^\ast)}\frac{1}{\sqrt{\varphi_0(x)+E}}\diff E'\\
&=-\frac{2}{\sqrt{2}} \int_{E^\ast}^{E_{\min}}\frac{1}{\varphi_0'(\varphi_0^{-1}(-E'))\sqrt{E-E'}}\diff E'-\frac{2}{\sqrt{2}} \int_{E_{\min}}^{E^\ast}\frac{1}{\varphi_0'(\varphi_0^{-1}(-E'))\sqrt{E-E'}}\diff E'\\
&=\frac{2}{\sqrt{2}} \int_{E_{\min}}^{E^\ast}\frac{1}{\varphi_0'(x_{-}(E'))\sqrt{E-E'}}\diff E'-\frac{2}{\sqrt{2}} \int_{E_{\min}}^{E^\ast}\frac{1}{\varphi_0'(x_{+}(E'))\sqrt{E-E'}}\diff E'.
\end{align*}
From Lemma~\ref{lema:turn_points}, the Taylor expansion for $\varphi_0$ and $\varphi_0'$ around $x_0$, and $\varphi_0''(x_0)<0$ we conclude that that there exists a constant $C>0$ (which depends on $E^\ast$) such that
\begin{equation*}
\frac{2}{\sqrt{2}} \int_{E_{\min}}^{E^\ast}\frac{1}{\varphi_0'(x_{-}(E'))\sqrt{E-E'}}\diff E'\leq C \int_{E_{\min}}^{E^\ast}\frac{1}{\sqrt{E'-E_{\min}}\sqrt{E-E'}}\diff E'.
\end{equation*}
Changing variables $E'-E_{\min}=(E^\ast-E_{\min})\theta$ with $\theta\in[0,1]$ we get
\begin{align*}
	\int_{E_{\min}}^{E^\ast}\frac{1}{\sqrt{E'-E_{\min}}\sqrt{E-E'}}\diff E'&\leq\int_{0}^{1}\frac{\diff\theta}{\sqrt{\theta}\sqrt{\frac{E-E_{\min}}{E^\ast-E_{\min}}-\theta}} \le \int_{0}^{1}\frac{\diff\theta}{\sqrt{\theta}\sqrt{1-\theta}} <\infty.
\end{align*}
Therefore $T_\ast$ is bounded from above and the claim follows for all $E\in(E^\ast,0)$.  A similar argument yields for the second integral.

Now let $E\in (0,|E^{\ast}|)\subset I_2$, and recall that the period function in this region is given by
\begin{equation}\label{E:PERIODEXT}
	T(E)=\int_{0}^{1}\frac{1}{\sqrt{2E+2\varphi_0(x)}}\diff x.
\end{equation}
By the Taylor expansion of $\varphi_0$ around $0$ and $1$ we obtain that there exist $\xi=\xi(x)\in[0,x]$ and $\zeta=\zeta(x)\in[x,1]$ such that
\begin{equation*}
	\varphi_0(x)=\frac{\varphi_0''(0)}{2}x^2\left(1+\frac{2\varphi_0'''(\xi(x))}{6\varphi_0''(0)}x\right)
\end{equation*}
and
\begin{equation*}
	\varphi_0(x)=\frac{\varphi_0''(1)}{2}(1-x)^2\left(1-\frac{2\varphi_0'''(\zeta(x))}{6\varphi_0''(1)}(1-x)\right).
\end{equation*}
Since $\varphi_0''(0)>0$ and the terms in parentheses converge to $1$ as $x\rightarrow 0$ and $x\rightarrow 1$ respectively, then there exists $\delta>0$ and a constant $C>0$ (which depends on $\delta$ and $E^\ast$) such that
\begin{equation*}
	\frac{1}{C}\int_{0}^{\delta}\frac{1}{\sqrt{E+x^2}}\diff x\leq T(E)\leq C\left(\int_{0}^{\delta}\frac{1}{\sqrt{E+x^2}}\diff x+1\right).
\end{equation*} 
Now note that
\begin{equation*}
	\int_{0}^{\delta}\frac{1}{\sqrt{E+x^2}}\diff x=\int_{0}^{\frac{\delta}{\sqrt{E}}}\frac{1}{\sqrt{1+x^2}}\diff x=\log\left(\sqrt{\frac{\delta^2}{E}+1}+\frac{\delta}{\sqrt{E}}\right),
\end{equation*}
which implies that for small enough values of $|E^\ast|>0$ the following inequalities hold
\begin{equation}
	\label{eq:hypestim_out}
	\frac{1}{C}|\log(E)|\leq T(E)\leq C|\log(E)|, \ \ E\in(0,|E^\ast|).
\end{equation}
This concludes the proof of~\eqref{E:TASYMPBOUNDAPP} for all $E\in \Ihypall$.

{\em Step 3: Exterior region: $E\in \Iext$.} Since we assume $E\in(\frac{|E^\ast|}{2},\infty)$ here, this bound is very simple. The upper bound follows by letting $\varphi_0(x)\ge0$ in~\eqref{E:PERIODEXT} and the lower bound  follows from $\varphi_0(x)\le-\Emin$, which gives
\begin{equation}\label{E:TEXT}
	\frac{1}{(2E-2E_{\min})^{\frac{1}{2}}}\leq T(E)\leq\frac{1}{(2E)^{\frac{1}{2}}}
\end{equation}
\end{proof}


\begin{lemma}[Sharp bounds on the first derivative of the period function]\label{L:ASYM_DERPER}
Let $\varphi_0$ satisfy the assumptions $(\varphi1)$--$(\varphi4)$. Then
there exists $\Emin<E^\ast<0$ satisfying~\eqref{E:PROP_E} and a constant $C>0$ such that for any $E\in I$ the following bounds hold:
\begin{align}
       \frac{1}{C}\Big(\frac{1}{|E|}\chh+  \che + |E|^{-\frac32}\chext\Big)\leq |T'(E)| \leq C \Big(\frac{1}{|E|}\chh+  \che+ |E|^{-\frac32}\chext\Big),\label{E:TPRIMEASYMPBOUNDAPP}
\end{align}
where we recall~\eqref{E:CHISDEF} and the period function~\eqref{E:PERIODDEF}.
\end{lemma}


\begin{proof}
{\em Step 1. Elliptic region: $E\in \Iell$.} This result is well known and follows verbatim the proof from Appendix B of~\cite{HRS2021}.

{\em Step 2. Hyperbolic region: $E\in \Ihypall$.} Assume first $E\in\Ihyp\cap I_1=(E_\ast,0)$. From \cite[Lemma A.10]{HRSS2023}, we have the following expression for the derivative of the period function
\begin{equation}
T'(E)=\frac{1}{E-E_{\min}}\int_{x_{-}(E)}^{x_{+}(E)}\frac{G(x)}{\sqrt{2E+2\varphi_0(x)}}\diff x \label{E:TPRIMEFORM}
\end{equation}
where
\begin{equation}
        G(x)=\begin{cases}
			\frac{\varphi_0'(x)^2-2\varphi_0''(x)(\varphi_0(x)-\varphi_0(x_0))}{\varphi_0'(x)^2} & \text{if $x\neq x_0$}\\
            0 & \text{if $x=x_0$} 
		 \end{cases}
		 \label{E:BIGGDEF}
\end{equation}
For each $E\in\Ihypall$ we decompose $T'$ as follows
\begin{equation*}
	T'(E)=T_{-}'(E)+T_{+}'(E)+T_{\ast}'(E)
\end{equation*}
where
\begin{align*}
	T_{-}'(E)&:=\frac{1}{E-E_{\min}}\int_{x_{-}(E)}^{x_{-}(E^\ast)}\frac{G(x)}{\sqrt{2E+2\varphi_0(x)}}\diff x,\\
	T_{+}'(E)&:=\frac{1}{E-E_{\min}}\int_{x_{+}(E^\ast)}^{x_{+}(E)}\frac{G(x)}{\sqrt{2E+2\varphi_0(x)}}\diff x,\\
	T_{\ast}'(E)&:=\frac{1}{E-E_{\min}}\int_{x_{-}(E^\ast)}^{x_{+}(E^\ast)}\frac{G(x)}{\sqrt{2E+2\varphi_0(x)}}\diff x.
\end{align*}
From the Taylor expansion of $\varphi_0,\varphi_0'$, and $\varphi_0''$ around $0$ it follows that 
\begin{align*}
    G(x)=\frac{2\varphi_0(x_0)}{\varphi_0''(0)x^2}(1+o_{x\rightarrow 0}(1)).
\end{align*}

Note that since $G>0$ for small values of $x$ we choose $\Emin<E^\ast<0$ close enough to $0$ such that $G>0$ in $[x_{-}(E),x_{-}(E^\ast)]$ for all $E\in(E^\ast,0)$. This implies that there exist constants $C,C'>0$ such that for every $E\in(E^\ast,0)$ the following inequality holds
\begin{equation*}
\frac{C}{E-E_{\min}}\int_{x_{-}(E)}^{x_{-}(E^\ast)}\frac{1}{x^2\sqrt{E+\varphi_0(x)}}\diff x\leq T_{-}'(E)\leq\frac{C'}{E-E_{\min}}\int_{x_{-}(E)}^{x_{-}(E^\ast)}\frac{1}{x^2\sqrt{E+\varphi_0(x)}}\diff x
\end{equation*}
and in a similar way as in the proof of Lemma~\ref{L:ASYM_PER} we can conclude that
\begin{equation*}
	T_{-}'(E)\geq \frac{C}{E-E_{\min}}\int_{x_{-}(E)}^{x_{-}(E^\ast)}\frac{1}{x^2\sqrt{\varphi_0'(x_{-}(E))(x-x_{-}(E))+C_1(x-x_{-}(E))^2}}\diff x
\end{equation*}
and
\begin{equation}\label{E:TMPRIME}
T_{-}'(E)\leq \frac{C'}{E-E_{\min}}\int_{x_{-}(E)}^{x_{-}(E^\ast)}\frac{1}{x^2\sqrt{\varphi_0'(x_{-}(E))(x-x_{-}(E))}}\diff x
\end{equation}
Using the  change of variables~\eqref{E:CHVAR}  we get
\begin{equation*}
T_{-}'(E)\geq\frac{C}{x_{-}(E)^{2}(E-E_{\min})}\int_{0}^{\frac{C_1(x_{-}(E^\ast)-x_-(E))}{\varphi_0'(x_{-}(E))}}\frac{1}{(P_1(E)y+1)^2\sqrt{y^2+y}}\diff y
\end{equation*}
where we introduce the notation 
\begin{align}\label{E:PONEDEF}
P_1(E):=\varphi_0'(x_{-}(E))/(C_1x_{-}(E)). 
\end{align}
Then note that the integral on the right can be rewritten as 
\begin{equation*}
    \int_{0}^{\infty}\frac{\diff y}{(P_1(E)y+1)^2\sqrt{y+y^2}}-\int_{\frac{C_1(x_{-}(E^\ast)-x_-(E))}{\varphi_0'(x_{-}(E))}}^{\infty}\frac{\diff y}{(P_1(E)y+1)^2\sqrt{y+y^2}}.
\end{equation*}
From Lemma~\ref{lema:turn_points} and $P_1(E)=O_{E\to0}(1)$ the error converges to $0$ as follows
\begin{align*}
    \left|\int_{\frac{C_1(x_{-}(E^\ast)-x_-(E))}{\varphi_0'(x_{-}(E))}}^{\infty}\frac{\diff y}{(P_1(E)y+1)^2\sqrt{y+y^2}}\right|&\leq \frac{C}{P_1(E)^2}\int_{\frac{C_1(x_{-}(E^\ast)-x_-(E))}{\varphi_0'(x_{-}(E))}}^\infty\frac{1}{y^{\frac{5}{2}}}\diff y \\
    &\leq C|E|^{\frac{3}{4}} (1+o_{E\to0}(1)),
\end{align*}
where we have used $\varphi_0'(x_{-}(E))\approx |E|^{\frac{1}{2}}$ for small $|E|$, which follows from~Lemma~\ref{lema:turn_points}. Thus
\begin{equation}
	\label{eq:bound_der_per1}
	T_{-}'(E)\geq \frac{C}{|E|}(1+o_{E\to0}(1)),
\end{equation}
where we also used $|E-\Emin|\le\Emin$ for $E\in \Ihyp\cap I_1$.
On the other hand, again from Lemma~\ref{lema:turn_points},~\eqref{E:TMPRIME}, change of variables $x\to |x'(E)|x$, and $|E-\Emin|\ge |E_\ast-\Emin|$ we have
\begin{align}
	T_{-}'(E)
&\leq\frac{C'}{x_{-}(E)^{2}}\int_{1}^{\infty}\frac{1}{x^2\sqrt{x-1}}\diff x 
\leq\frac{C'}{|E|}\int_{1}^{\infty}\frac{1}{x^2\sqrt{x-1}}\diff x.  \label{eq:bound_der_per2}
\end{align}
Then from \eqref{eq:bound_der_per1} and \eqref{eq:bound_der_per2} we have that there exist $C>0$ such that for every $E\in(E^\ast,0)$
\begin{equation}
	\frac{1}{C}\frac{1}{|E|}\leq T_{-}'(E)\leq \frac{C}{|E|}.
\end{equation}
 A completely analogous argument yields
\begin{equation}
	\frac{1}{C}\frac{1}{|E|}\leq T_{+}'(E)\leq \frac{C}{|E|}, \ \ E\in(E^\ast,0).
\end{equation}
 Also note that
\begin{equation*}
	|T_{\ast}'(E)|\leq\frac{1}{E-E_{\min}}\int_{x_{-}(E^\ast)}^{x_{+}(E^\ast)}\frac{|G(x)|}{\sqrt{2E+2\varphi(x)}}\diff x.
\end{equation*}
Since $G$ is continuous and $|G(x)|\lesssim 1$ in $[x_{-}(E^\ast),x_{+}(E^\ast)]$, by repeating the argument from the proof of Lemma~\ref{L:ASYM_PER} in this region, there exist a constant $C>0$ (which might depend on $E^\ast$) such that 
\begin{equation*}
	|T_{\ast}'(E)|\leq C, \ \ E\in(E^\ast,0).
\end{equation*}
Now let $E\in (0,|E^{\ast}|)\in I_2$, and recall that the period function in this region is given by
\begin{equation}\label{E:DERPERIODEXT}
	T'(E)=-\int_{0}^{1}\frac{1}{(2E+2\varphi_0(x))^{\frac{3}{2}}}\diff x.
\end{equation}
The same argument as in Lemma~\ref{L:ASYM_PER} implies than there exist $\delta>0$ and a constant $C>0$ (which depends on $\delta$ and $E^\ast$) such that
\begin{equation*}
	-\frac{1}{C}\int_{0}^{\delta}\frac{1}{(E+x^2)^{\frac{3}{2}}}\diff x\leq T'(E)\leq -C\left(\int_{0}^{\delta}\frac{1}{(E+x^2)^{\frac{3}{2}}}\diff x+1\right).
\end{equation*} 
Note that 
\begin{equation*}
\int_{0}^{\delta}\frac{1}{(E+x^2)^{\frac{3}{2}}}\diff x=\frac{1}{E}\int_{0}^{\frac{\delta}{\sqrt{E}}}\frac{1}{(1+y^2)^{\frac{3}{2}}}\diff y\leq\frac{1}{E}
\end{equation*}
which implies that there exists a constant $C>0$ such that
\begin{equation*}
-\frac{1}{CE}\leq T'(E)\leq -\frac{C}{E}.
\end{equation*}
This concludes the proof of~\eqref{E:TPRIMEASYMPBOUNDAPP} for all $E\in \Ihypall$.

{\em Step 3: Exterior region: $E\in \Iext$.} Again this argument is simpler. We differentiate~\eqref{E:PERIODEXT} with respect to $E$ and then similarly to the proof of~\eqref{E:TEXT} we conclude
\begin{equation}
-\frac{1}{(2E)^{\frac{3}{2}}}\leq T'(E)\leq-\frac{1}{(2E-2E_{\min})^{\frac{3}{2}}}.
\end{equation}
\end{proof}


\begin{lemma}[Sharp bounds on the second derivative of the period function]\label{L:ASYM_DER2PER}
Let $\varphi_0$ satisfy the assumptions $(\varphi1)$--$(\varphi4)$. Then
there exists $\Emin<E^\ast<0$ satisfying~\eqref{E:PROP_E} and a constant $C>0$ such that for any $E\in I$ the following bounds hold:
\begin{align}
       \frac{1}{C}\Big(\frac{1}{|E|^2}\chh+  \che + |E|^{-\frac52}\chext\Big)\leq |T''(E)| \leq C \Big(\frac{1}{|E|^2}\chh+  \che+ |E|^{-\frac52}\chext\Big),\label{E:TPRIMEPRIMEASYMPBOUNDAPP}
\end{align}
where we recall~\eqref{E:CHISDEF} and the period function~\eqref{E:PERIODDEF}.
\end{lemma}


\begin{proof}
{\em Step 1. Elliptic region: $E\in \Iell$.} This result is well known and follows from Appendix B in~\cite{HRS2021}.

{\em Step 2. Hyperbolic region: $E\in \Ihypall$.} We first focus on the case $E\in\Ihyp\cap I_1=(E_\ast,0)$. Using \cite[eqn.(A.28)]{HRSS2023} we have the following expression for the second derivative of the period function
\begin{equation*}
    T''(E)=\frac{1}{(E-E_{\min})^2}\int_{x_{-}(E)}^{x_{+}(E)}G_1'(x)\frac{E+E_{\min}+2\varphi_0(x)}{\sqrt{2E+2\varphi_0(x)}}\diff x,
\end{equation*} 
where we defined
\begin{equation*}
        G_1(x)=\begin{cases}
			\frac{G(x)}{\varphi_0'(x)} & \text{if $x\neq x_0$},\\
            -\frac{\varphi_0'''(x_0)}{3\varphi_0''(x_0)^2} & \text{if $x=x_0$} ,
		 \end{cases}
\end{equation*}
and $G$ is given by~\eqref{E:BIGGDEF}.
It is straightforward to see that if $x\neq x_0$, then
\begin{equation*}
    G_1'(x)=\frac{6(\varphi_0(x)-\varphi_0(x_0))\varphi_0''(x)^2-2\varphi_0'''(x)(\varphi_0(x)-\varphi_0(x_0))\varphi_0'(x)-3\varphi_0'(x)^2\varphi_0''(x)}{\varphi_0'(x)^4}.
\end{equation*}
We next split $T''$
\begin{equation*}
	T''(E)=T_{-}''(E)+T_{+}''(E)+T_{\ast}''(E),
\end{equation*}
where
\begin{align*}
	T''_{-}(E)&=\frac{1}{(E-E_{\min})^2}\int_{x_{-}(E)}^{x_{-}(E^\ast)}G_1'(x)\frac{E+E_{\min}+2\varphi_0(x)}{\sqrt{2E+2\varphi_0(x)}}\diff x,\\
	T''_{+}(E)&=\frac{1}{(E-E_{\min})^2}\int_{x_{+}(E^\ast)}^{x_{+}(E)}G_1'(x)\frac{E+E_{\min}+2\varphi_0(x)}{\sqrt{2E+2\varphi_0(x)}}\diff x,\\
	T''_{\ast}(E)&=\frac{1}{(E-E_{\min})^2}\int_{x_{-}(E^\ast)}^{x_{+}(E^\ast)}G_1'(x)\frac{E+E_{\min}+2\varphi_0(x)}{\sqrt{2E+2\varphi_0(x)}}\diff x.
\end{align*}
Analogously to the proofs of Lemmas~\ref{L:ASYM_PER}--\ref{L:ASYM_DERPER}, from the Taylor expansion we can prove that
\begin{equation*}
    \label{eq:exp_taylor_g1}
    G_1'(x)=-\frac{6\varphi_0(x_0)}{\varphi_0''(0)^2x^4}(1+o_{x\to0}(1)).
\end{equation*}
If we choose $E^\ast>E_{\min}/2$, then $E_{\min}<E+E_{\min}+2\varphi_0(x)\leq E_{\min}-2E^\ast<0$ for every $E\in(E^\ast,0)$ and every $x\in [x_{-}(E),x_{-}(E^\ast)]$. This implies that there exist positive constants $C$ and $C'$ such that
\begin{equation*}
	\frac{C'}{(E-E_{\min})^2}\int_{x_{-}(E)}^{x_{-}(E^\ast)}\frac{1}{x^4\sqrt{E+\varphi_0(x)}}\diff x\leq |T_{-}''(E)|\leq \frac{C}{(E-E_{\min})^2}\int_{x_{-}(E)}^{x_{-}(E^\ast)}\frac{1}{x^4\sqrt{E+\varphi_0(x)}}\diff x.
\end{equation*}
From our choice of $E^\ast$, there exists a constant $C>0$ (which again might depend on $E^\ast$) such that for every $E\in(E^\ast,0)$
\begin{equation*}
   |T_{-}''(E)|\leq C\varphi_0'(x_{-}(E))^{-\frac{1}{2}}\int_{x_{-}(E)}^{x_{-}(E^\ast)}\frac{1}{x^4\sqrt{x-x_{-}(E)}}\diff x,
\end{equation*}
where we have used~\eqref{E:TAYLORIMP} and $|E-\Emin|\ge |E_\ast-\Emin|$, similarly to the previous two lemmas.
Note that
\begin{equation*}
    \int_{x_{-}(E)}^{x_{-}(E^\ast)}\frac{1}{x^4\sqrt{x-x_{-}(E)}}\diff x=\frac{1}{x_{-}(E)^{\frac{7}{2}}}\int_{1}^{\frac{x_{-}(E^{\ast})}{x_{-}(E)}}\frac{1}{y^4\sqrt{y-1}}\diff y\leq\frac{1}{x_{-}(E)^{\frac{7}{2}}}\int_{1}^{\infty}\frac{1}{y^4\sqrt{y-1}}\diff y,
\end{equation*}
and thus 
\begin{align*}
    |T_{-}''(E)|&\leq\frac{C}{(x_{-}(E)^4}=\frac{C}{|E|^2}(1+o_{E\to0}(1)).
\end{align*}
Similarly to the proof of Lemma~\ref{L:ASYM_PER} we can prove that there exist a constant $C_1>0$ such that
\begin{equation*}
	|T_{-}''(E)|\geq \frac{C'}{(E-E_{\min})^2}\int_{x_{-}(E)}^{x_{-}(E^\ast)}\frac{1}{x^4\sqrt{\varphi_0'(x_{-}(E))(x-x_{-}(E))+C_1(x-x_{-}(E))^2}}\diff x.
\end{equation*}
Using the notation~\eqref{E:PONEDEF} as in the proof of Lemma~\ref{L:ASYM_DERPER} and $|E-\Emin|\le |\Emin|$, we conclude
\begin{equation*}
|T_{-}''(E)|\geq\frac{C}{x_{-}(E)^{4}}\int_{0}^{\frac{C_1(x_{-}(E^\ast)-x_-(E))}{\varphi_0'(x_{-}(E))}}\frac{1}{(P_1(E)y+1)^4\sqrt{y^2+y}}\diff y.
\end{equation*}
The integral above is now rewritten as 
\begin{equation*}
    \int_{0}^{\infty}\frac{\diff y}{(P_1(E)y+1)^4\sqrt{y+y^2}}-\int_{\frac{C_1(x_{-}(E^\ast)-x_-(E))}{\varphi_0'(x_{-}(E))}}^{\infty}\frac{\diff y}{(P_1(E)y+1)^4\sqrt{y+y^2}};
\end{equation*}
the error converges to $0$ since
\begin{align*}
    \left|\int_{\frac{C_1(x_{-}(E^\ast)-x_-(E))}{\varphi_0'(x_{-}(E))}}^{\infty}\frac{\diff y}{(P_1(E)y+1)^4\sqrt{y+y^2}}\right|&\leq \frac{C}{P_1(E)^4}\int_{\frac{C_1(x_{-}(E^\ast)-x_-(E))}{\varphi_0'(x_{-}(E))}}^\infty\frac{1}{y^{\frac{9}{2}}}\diff y = O_{E\to0}(|E|^{\frac{7}{4}}).
\end{align*}
It then follows that there exists $C>0$ such that for every $E\in(E^\ast,0)$
\begin{equation}
	\frac{1}{C}\frac{1}{|E|^2}\leq T_{\pm}''(E)\leq \frac{C}{|E|^2},
\end{equation}
where the bound on $T_+''$ follows analogously.
Finally,
\begin{equation*}
|T_{\ast}''(E)|\leq\frac{C}{(E-E_{\min})^2}\int_{x_{-}(E^\ast)}^{x_{+}(E^\ast)}\frac{|G_1'(x)|}{\sqrt{2E+2\varphi_0(x)}}\diff x.
\end{equation*}
Since $G_1'$ is continuous we conclude with the same argument as in Lemma~\ref{L:ASYM_PER} that there exists a constant $C>0$ such that
$|T_{\ast}''(E)|\leq C$.
Now let $E\in \Ihyp\cap I_2= (0,|E^{\ast}|)$, and recall that the period function in this region is given by
\begin{equation}\label{E:SECDERPERIODEXT}
	T''(E)=3\int_{0}^{1}\frac{1}{(2E+2\varphi_0(x))^{\frac{5}{2}}}\diff x.
\end{equation}
The same argument as in Lemma~\ref{L:ASYM_PER} implies than there exist $\delta>0$ and a constant $C>0$ (which depends on $\delta$ and $E^\ast$) such that
\begin{equation*}
	\frac{1}{C}\int_{0}^{\delta}\frac{1}{(E+x^2)^{\frac{5}{2}}}\diff x\leq T''(E)\leq C\left(\int_{0}^{\delta}\frac{1}{(E+x^2)^{\frac{5}{2}}}\diff x+1\right).
\end{equation*} 
and note that 
\begin{equation*}
\int_{0}^{\delta}\frac{1}{(E+x^2)^{\frac{5}{2}}}\diff x=\frac{1}{E^2}\int_{0}^{\frac{\delta}{\sqrt{E}}}\frac{1}{(1+y^2)^{\frac{5}{2}}}\diff y\approx \frac{C}{E^2}
\end{equation*}
which implies that there exists a constant $C>0$ such that
$\frac{1}{CE^2}\leq T''(E)\leq\frac{C}{E^2}$.
This concludes the proof of~\eqref{E:TPRIMEPRIMEASYMPBOUNDAPP} for all $E\in \Ihypall$.

{\em Step 3: Exterior region: $E\in \Iext$.} Here the estimates come from the fact that
\begin{equation*}
(2E-2E_{\min})^{-\frac{5}{2}}\leq T'(E)\leq (2E)^{-\frac{5}{2}}
\end{equation*}
\end{proof}


\begin{lemma}\label{L:RSLEMMA}
Let $\varphi_0$ satisfy the assumptions $(\varphi1)$--$(\varphi4)$.
Let 
\[
I_{\pm}(E):=(x_{-}(E),x_{+}(E)), \ \ E\in I.
\] 
Then there exists $\Emin<E^\ast<0$ satisfying~\eqref{E:PROP_E} and a constant $C>0$ such that for any $E\in I$ and any $x\in I_{\pm}(E)$ the following bound holds
\begin{align}
\big|\pa_E\big(T(E)\th(x,E)\big)\big|\lesssim \mathcal{I}_{\text{ell}}(x,E)+\mathcal{I}_{\text{hyp}}(x,E)+\mathcal{I}_{\text{ext}}(x,E),\notag
\end{align}
where
\begin{align*}
\mathcal{I}_{\text{ell}}(x,E)&:=\frac{\chi_{I_\pm(E)}(x)\che(E)}{(E-\Emin)^{\frac{1}{2}}\sqrt{E+\varphi_0(x)}},\\
\mathcal{I}_{\text{hyp}}(x,E)&:=T'(E)\left(
\frac{\chi_{I_{\pm}(E)\setminus\overline{I_{\pm}(E^\ast)}}(x)\chi_{\Ihypall\cap I_1}(E)}{\sqrt{E+\varphi_0(x)}}+\chi_{I_{\pm}(E)}(x)\chi_{\Ihypall\cap I_2}(E)\right),\\
\mathcal{I}_{\text{ext}}(x,E)&:=T'(E)\chi_{I_{\pm}(E)}(x)\chext(E).
\end{align*}
\end{lemma}


\begin{proof}
Note that
\begin{align}\label{E:IDEF}
\big|\pa_E\big(T(E)\th(x,E)\big)\big| =\left|\frac{\pa}{\pa E}\int_{x_{-}(E)}^{x}\frac{\diff y}{\sqrt{2E+2\varphi_0(y)}}\right|,
\end{align}
where we recall from~\eqref{E:TURN_POINTS_EXT} that $x_{-}(E)=0$ for $E>0$. We may now adapt the derivation of the formula for the derivative of the period function~\cite[Lemma A.10]{HRSS2023}
to obtain the formula
\begin{align}\label{E:RSDEF}
    \big|\pa_E\big(T(E)\th(x,E)\big)\big| =\begin{cases} \displaystyle\mathcal R(x,E) + \mathcal S(x,E) & \text{if }E\in I_1, x\in I_{\pm}(E),\\
\displaystyle\mathcal{T}(x,E) & \text{if }E\in I_2,x\in[0,1],
\end{cases}
\end{align}
    where
    \begin{align}
    \mathcal R(x,E)&:=\frac{\varphi_0(x_0)-\varphi_0(x)}{(E-E_{\min})\varphi_0'(x)\sqrt{2E+2\varphi_0(x)}}, \ \ x\in I_\pm (E), \ E\in I_1, \label{E:RDEF} \\
\mathcal S(x,E)&:=\frac{1}{2(E-E_{\min})}\int_{x_{-}(E)}^{x}\frac{G(y)} {\sqrt{2E+2\varphi_0(y)}}\diff y, \ \ x\in I_\pm (E), \ E\in I_1,\label{E:SDEF}\\
\mathcal T(x,E)&:=-\int_{0}^x\frac{1}{(2E+2\varphi_0(y))^{\frac{3}{2}}}\diff y, \ \ E\in I_2, \label{E:SDEF}
\end{align}
where we recall the definition~\eqref{E:BIGGDEF} of $G(y)$.
From the Taylor expansion of $\varphi_0$ around $x_0$ we have
\begin{equation*}
\left|\frac{\varphi_0(x_0)-\varphi_0(x)}{\varphi_0'(x)}\right|=\frac12\left|x-x_0\right|(1+o_{x\to x_0}(1))
\end{equation*}
and
\begin{equation*}
E-E_{\min}\geq\varphi_0(x_0)-\varphi_0(x)\geq-\frac{\varphi_0''(x_0)}{2}(x_0-x)^2(1+o_{x\to x_0}(1)).
\end{equation*}
This gives 
\begin{equation*}
\left|\frac{\varphi_0(x_0)-\varphi_0(x)}{\varphi_0'(x)}\right|\leq C\sqrt{E-E_{\min}},
\end{equation*}
and therefore
\begin{equation}\label{E:BOUNDR1}
|\mathcal R(x,E)| \leq \frac{C}{\tE^{\frac{1}{2}}\sqrt{E+\varphi_0(x)}}, \ \ \tE = E-\Emin, \ x\in I_\pm (E), \ \ E\in \Iell,
\end{equation}
where we recall the definition~\eqref{E:IELLDEF} of $\Iell$.
On the other hand, since $E^\ast$ satisfies~\eqref{E:PROP_E}, then $\varphi_0$ is convex on $I_{\pm}(E)\setminus \overline{I_{\pm}(E^\ast)}$, for every $E\in (E^\ast,0)$. This implies that
\begin{equation}\label{E:BOUNDR2}
|\mathcal R(x,E)| \leq\frac{C}{|\varphi_0'(x_{-}(E))|\sqrt{E+\varphi_0(x)}}\leq\frac{C}{|E|^{\frac{1}{2}}\sqrt{E+\varphi_0(x)}}\leq \frac{CT'(E)}{\sqrt{E+\varphi_0(x)}}
\end{equation}
for every $x\in I_{\pm}(E)\setminus \overline{I_{\pm}(E^\ast)}$ and $E\in (E^\ast,0)$, where we used that $\tE^{-\frac{1}{2}}\lesssim 1$ in this region, and in the last inequality we used Theorem~\ref{T:PERIODASYMP} and the fact that $|E|^{-\frac{1}{2}}\leq|E|^{-1}\lesssim T'(E)$ if $|E^\ast|\leq 1$. Thus from \eqref{E:BOUNDR1} and \eqref{E:BOUNDR2} we get
\begin{align}\label{E:BOUNDR}
|\mathcal R(x,E)| &\lesssim \frac{\chi_{I_{\pm}(E)}(x)\che(E)}{\tE^{\frac{1}{2}}\sqrt{E+\varphi_0(x)}}+\frac{\chi_{I_{\pm}(E)\setminus\overline{I_{\pm}(E^\ast)}}(x)\chi_{\Ihypall\cap I_1}(E)T'(E)}{\sqrt{E+\varphi_0(x)}}
\end{align}
To bound $\mathcal{S}$ we note that for any $E\in I_1$ and $x\in I_\pm(E)$ we have
\begin{align}\label{E:BOUNDS}
|\mathcal S(x,E)|\le \frac12 \frac1{E-\Emin}\int_{x_-(E)}^{x_+(E)}\frac{G(y)}{\sqrt{E+\varphi_0(y)}}\diff y = \frac12 T'(E) \chi_{I_{\pm}(E)}(x),
\end{align}
where we have used~\eqref{E:TPRIMEFORM}.
For $E\in I_2$ and $x\in[0,1]$ we use~\eqref{E:DERPERIODEXT} to see that
\begin{align}\label{E:BOUNDT}
|\mathcal T(x,E)|\leq |T'(E)|.
\end{align}
 The result now follows from~\eqref{E:RSDEF} and the bounds \eqref{E:BOUNDR}, \eqref{E:BOUNDS}, and \eqref{E:BOUNDT}.
\end{proof} 


\section{Spectral properties}\label{A:SPECTRAL}


\begin{lemma}\label{L:FORCEBOUND}
Let $(f_0,\varphi_0)$ be a steady state of the VP-system~\eqref{E:VLASOV2}--\eqref{E:POISSON2}, so that 
$(f_0,\varphi_0)$ is of the form~\eqref{E:EDEF}, 
the microscopic equation of state $\mu(\cdot)$ satisfies properties~\eqref{E:MUREG}--\eqref{E:MU4}
and the associated 
electrostatic potential $\varphi_0$ satisfies properties $(\varphi1)$--$(\varphi4)$. Then the following statements hold.
\begin{enumerate}
\item[(a)]
For any $\beta\in[0,\frac12)$ there exists a constant $C_\beta>0$ such that for any $f\in\H$ we have the bound
\begin{align}
\|\pa_x\varphi_f\|_{L^\infty(\mathbb T^1)} \le C_\beta \eps^{\frac{\beta}{2}} \|f\|_\H,\notag
\end{align}
where we recall~\eqref{E:HILBERTSPACE}.
\item[(b)]
The operator $\H \ni f\mapsto \eps^2 \mu'(\eps^2E) v\pa_x \varphi_f \in \H$ is well-defined, compact  and satisfies the bound
\begin{align}\label{E:WP}
\|\eps^2 \mu'(\eps^2E) v\pa_x \varphi_f \|_\H \le C \eps^{\frac12} \|f\|_\H.
\end{align}
\item[(c)]
Assume additionally that $\mu\in C^2(h,\infty)$. Then the operator~\eqref{E:LDEF} is bounded on $\H$.
\end{enumerate}
\end{lemma}


\begin{proof}
{\em Proof of part (a).}
From the Poisson equation we obtain the formula
\begin{align}\label{E:FORCE}
\pa_x\varphi_f(x) = - \eps \int_0^x\int f (t,y,v)\diff v\diff y - \eps \iint y f(y,v)\diff v \diff y.
\end{align}
Note that 
\begin{align}
& \Big|\int_0^x\int f (t,y,v)\diff v\diff y\Big| + \Big|  \iint y f(y,v)\diff v \diff y\Big| \notag\\
& \le \big(\int_0^1(\chi_{[0,x]}(y)+|y|)\int  |\mu'|\diff v\diff y\big)^{\frac12} \big(\int_0^x \int\frac{f^2}{|\mu'|}\diff y\diff v\big)^{\frac12} \label{E:BOUNDPOT}
\end{align}
Let $\beta\in[0,\frac12)$ be given. Then by~\eqref{E:MUMINUSPLUSREG}--\eqref{E:MUTAIL} we obtain
\begin{align}
&\int|\mu'(\eps^2E(x,v))| \diff v \notag\\
 &\le C \int_{E(x,v)\le0}1\diff v + C \int_{1\ge E(x,v)>0} \frac{1}{\eps\sqrt{E(x,v)}}\diff v + C \int_{E(x,v)>1} \frac1{1+ \eps^2 E(x,v)}\diff v \notag\\
& \le C + \frac C{\eps} \int_{\sqrt{2\varphi_0(x)}}^{\sqrt{2(1+\varphi_0(x))}} \frac1{\sqrt{\frac12 v^2-\varphi_0(x)}} \diff v 
+ \frac C{\eps^{2-\beta}} \int_{\sqrt{2(1+\varphi_0(x))}}^\infty \frac1{(\frac12 v^2-\varphi_0(x))^{1-\beta}} \diff v \notag\\
& = \frac C{\eps} \int_{1}^{\frac{\sqrt{2(1+\varphi_0(x))}}{\sqrt{2\varphi_0(x)}}} \frac1{\sqrt{w^2-1}} \diff w 
+ \frac C{\eps^{2-\beta}} \varphi_0(x)^{\beta-\frac12}\int_{\frac{\sqrt{2(1+\varphi_0(x))}}{\sqrt{2\varphi_0(x)}}}^\infty \frac1{(w^2-1)^{1-\beta}} \diff v \notag\\
& \le \frac C{\eps} |\log \varphi_0(x)| + \frac {C_\beta}{\eps^{2-\beta}}. \label{E:MUPRIMEBOUND}
\end{align}
Therefore
\begin{align}
\Big|\int_0^1(\chi_{[0,x]}(y)+|y|)\int  |\mu'|\diff v\diff y\Big| \le \frac {C_\beta}{\eps^{2-\beta}}.\notag
\end{align}
This bound and~\eqref{E:BOUNDPOT} imply the statement of the lemma.

{\em Proof of part (b).}
Taking the square of the left-hand side of~\eqref{E:WP} we obtain
\begin{align}
\eps^4 \int |\mu'(\eps^2 E(x,v))| v^2 |\pa_x\varphi_f|^2 \diff (x,v) &\le \eps^4 \|\pa_x\varphi_f\|_{L^\infty}^2 \sup_{x\in\mathbb T^1}\int |\mu'(\eps^2 E(x,v))| v^2 \diff v \notag\\
& \le \eps^4 \|\pa_x\varphi_f\|_{L^\infty}^2 \Big(C + \sup_{x\in\mathbb T^1}\int_{E>0} |\mu'(\eps^2 E(x,v))| v^2 \diff v\Big). \notag
\end{align}
Letting $\eps^2 E(x,v)= e$, we have $v^2 =\frac{2}{\eps^2}(e+\eps^2\varphi_0(x)) $ and $v \diff v = \frac {\diff e}{\eps^2}$. We may therefore bound the above integral by 
\begin{align}
\sup_{x\in\mathbb T^1}\int_{E>0} |\mu'(\eps^2 E(x,v))| v^2 \diff v \le C\eps^{-3} \int_0^\infty|\mu'(e)|\sqrt {e+\eps^2 \varphi_0(x_0)} \diff e \le C\eps^{-3},
\end{align}
where we have used  ($\mu3$), specifically the assumption~\eqref{E:COMP}. Using part (a) with $\beta=0$, we conclude~\eqref{E:WP}. To prove compactness, in the same way as in the proof of Corollary~\ref{C:LANDAU}, given a bounded sequence $(f_n)_{n\in\N}$ in $\H$, from Arzela-Ascoli theorem, after passing through subsequence we get a uniform limit $\Phi=\lim_{k\rightarrow\infty}\partial_x\varphi_{n_k}$. Similarly to the proof of the bound \eqref{E:WP} we estimate
\begin{align*}
\|\k^2\mu'(\k^2E)v\partial_x\varphi_{n_k}-\k^2\mu'(\k^2E)v\Phi\|_{\H}^2&\leq\k^4\sup_{x\in\mathbb T^1}\|\partial_x\varphi_{n_k}-\Phi\|^2\sup_{x\in\mathbb T^1}\int|\mu'(\k^2E(x,v)|v^2\diff v\\
&\leq C\k\sup_{x\in\mathbb T^1}\|\partial_x\varphi_{n_k}-\Phi\|^2,
\end{align*}
where the compactness follows from the fact that the last inequality converges to $0$ as $k\rightarrow\infty$.

{\em Proof of part (c).}
We note that 
\begin{align}
\| \eps^2\mu'(\eps^2 E) (\varphi_f-\int\varphi_f\diff x)\|_\H^2 & = \eps^4\iint |\mu'(\eps^2 E)| |\varphi_f-\int\varphi_f\diff \tilde x)|^2\diff (x,v) \notag\\
& \le C\eps^4  \int (\varphi_f-\int\varphi_f\diff \tilde x)|^2 \diff x \le C \|\pa_x\varphi_f\|_{L^2}^2 \notag\\
& \le C \eps^{4+\beta} \|f\|_\H^2,\notag
\end{align}
where we note that the $|\log\varphi_0(x)|$-term does not appear in~\eqref{E:MUPRIMEBOUND} for regular $\mu$.
\end{proof}

\bigskip 

{\em Proof of Lemma~\ref{L:ESSENTIAL}}
Using the representation~\eqref{E:DPOS}--\eqref{E:DNEG}, it is standard to see that $ \sigma_{\text{ess}}(\D) = i \, \mathbb R$, see for example~\cite{GuoLin,HRSS2023},
where we use that $\text{ran}\,T=(0,\infty)$ on $I=I_1\cup I_2$.  Since the operator $\H\ni f\mapsto \eps^2 |\mu'(\eps^2 E)|v\pa_x \varphi_f\in\H$ is compact and $\D$ is closed, the essential spectra of $\D$ and $\mathscr L$ are the same by Weyl's theorem. Description~\eqref{E:KERD} is equivaelnt to the statement that the zero-Fourier coefficient in $\th$ of its elements vanishes, which then follows easily from  the description~\eqref{E:KERD0}.


\section{Multi-branched microscopic equations of state}\label{A:SCHAMEL}

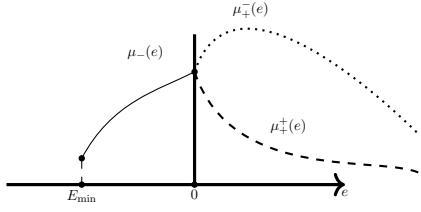
\begin{figure}
\begin{center}
\begin{tikzpicture}
\begin{scope}[scale=0.5, transform shape]


\coordinate [label=below:$0$] (E) at (10,7){};
\coordinate [] (G) at (8,11){};
\coordinate [label=below:$e$] (H) at (14,7){};
\coordinate [label=above:$\mu_-(e)$] (I) at (8.7,10.2){};
\coordinate [label=below:$\Emin$] (M) at (7,7){};
\coordinate [label=above:$\mu^{-}_+(e)$]  (K) at (11.5,11.3){};
\coordinate [label=above:$\mu^{+}_+(e)$]  (L) at (12.5,8.2){};

\draw[very thick] (E)--(5,7);
\draw[-{>},very thick] (E)--(H);
\draw[very thick] (E)--(10,11);
\draw (7,7.7) .. controls +(0.8,1.5) and +(-1,-0.5) .. (10,10);
\draw[dotted, thick] (10,10) .. controls +(1,3) and +(-1,1) .. (16,8.3);
\draw[dashed, thick] (10,10) .. controls +(1,-3) and +(-1,0.5) .. (16,7.3);
\draw[dashed] (7,7) -- (7,7.7);

\draw[fill=black] (10,7) circle (2pt);
\draw[fill=black] (10,10) circle (2pt);
\draw[fill=black] (7,7.7) circle (2pt);
\draw[fill=black] (7,7) circle (2pt);

\end{scope}
\end{tikzpicture}
\caption{Schematic depiction of a Schamel-like profile. Note that neither of the two branches $\mu_+^\pm$ matches the left branch in a $C^1$-fashion.}
\label{F:SCHAMEL}
\end{center}
\end{figure}

In order to allow for profiles which are not even in $v$ (at positive energies) we consider  a more general class of steady states
of the form
\begin{align}\label{E:FSCH}
f_0(x,v) = \begin{cases} \mu_-(\eps^2 E(x,v)) & E\le0\\ \mu_+^{\text{sign}(v)}(\eps^2 E(x,v)),  & E>0, \end{cases}
\end{align}
where $\text{sign}(v)=\pm$ if $\pm v>0$. Note that~\eqref{E:FSCH} formally solves the Vlasov equation since $\text{sign}(v)$ is a conserved quantity along the flow in the region $E>0$. In this case we let
\begin{align}\label{E:RHOPLUSNEW}
\rho_+(x) : = -\pa_{xx}\varphi_0(x) + \eps \int f_0(x,v)\diff v >0, \ \ x\in\mathbb T^1,
\end{align}
where the strict positivity needs to be shown with a suitable choice of $\varphi_0$.
We then require that each of the two functions $\mu^\pm(e):=\mu_-(e)\chi_{[\Emin,0]}(e)+\mu_+^\pm(e)\chi_{(0,\infty)}(e)$  satisfies $(\mu1)$--$(\mu3)$ while $(\mu4)$ is replaced by the assumption
\begin{enumerate}
\item[$(\mu4')$] $\mu_-$ is either strictly decreasing or strictly increasing, $\mu_+$ has at most one critical point, which is a local maximum.
\end{enumerate}  
A particular family of examples which motivates this class of steady states 
are 
the Schamel-type distributions~\cite{HSc1986,Hutch} given by
\begin{align}\label{E:SCHDEF}
S_{\al,\beta}(x,v) =  \begin{cases} \exp(\beta \eps^2 E(x,v)-\al^2) & E(x,v)\in(-\infty,0], \\ \exp(-(\text{sgn}(v) \eps \sqrt{E(x,v)}+\al)^2) & E(x,v)>0, \end{cases}
\end{align}
where $\al\ge0$ and $\beta>0$ are given parameters.
We may therefore write $\mu_-(e)= e^{\beta e-\al^2}$ and $\mu^{\pm}_+(e)=e^{-(\pm \sqrt e+\al)^2}$ (here it is understood that the choice of sign in the $\pm$-dependence corresponds to the sign of $v$) and express the  steady state $S_{\al,\beta}$ as a multi-branched function of $e$: $S_{\al,\beta}(x,v)=\mu_-(\eps^2 E)\chi_{[\Emin,0]}(e)+\mu_+^\pm(\eps^2 E)\chi_{(0,\infty)}(e)$. The two branches are illustrated in Figure~\ref{F:SCHAMEL}. We note that the lower brunch is decreasing on $I_2$ but the higher branch has a local maximum at some $e_\ast>0$.
A simple corollary of our method is the analogue of Theorem~\ref{T:MAIN} for profiles~\eqref{E:SCHDEF}.


\begin{corollary}\label{C:SCHAMEL}
Let $(f_0,\varphi_0)$ be a steady state of the VP-system~\eqref{E:VLASOV2}--\eqref{E:POISSON2}
of the form~\eqref{E:FSCH} such that
the functions $\mu^\pm(\cdot)=\mu_-\chi_{e\le0}+\mu_+^\pm\chi_{e>0}$ satisfy $(\mu1)$--$(\mu3)$ and $(\mu 4')$,
and the associated 
electrostatic potential $\varphi_0$ satisfies properties $(\varphi1)$--$(\varphi4)$ (or $(\varphi1)$--$(\varphi3)$ and $(\varphi4^*)$).
Then there exists an $\eps_0>0$ such that for any $0<\eps<\eps_0$ the  operator~\eqref{E:DL} $\mathscr L:D(\mathscr L)\subset \mathcal H\to\mathcal H$ has no non-zero embedded eigenvalues.
\end{corollary}


{\em Sketch of the proof of Corollary~\ref{C:SCHAMEL}.}
The proof is nearly identical, except we must allow for two $\mu$-profiles for $E>0$, see~\eqref{E:FSCH} and Figure~\ref{F:SCHAMEL}. According to properties $(\mu1)$--$(\mu3)$ and $(\mu 4')$, $\mu_+^\pm$ is either a) decreasing and in that case $\mu^\pm(e):=\mu_-(e)\chi_{[\Emin,0]}(e)+\mu_+^\pm(e)\chi_{(0,\infty)}(e)$ satisfies the assumptions $(\mu1)$--$(\mu4)$, or b) $\mu_+^\pm$ has a local maximum at some $e>0$ (this is precisely the case with $\mu_+^+$ in Figure~\ref{F:SCHAMEL}).
In the former case Lemma~\ref{L:ENERGYBASIC} applies directly and so do Lemmas~\ref{L:NRLEMMA},~\ref{L:TRAPLEMMA}, and~\ref{L:EXTLEMMA} to control the contribution of that branch of the solution. In the latter case, there exists an $e_\ast>0$ such that $\mu_+'>0$ on $(0,e_\ast)$ and $\mu_+'<0$ on $(e_\ast,\infty)$. In this case we modify Lemma~\ref{L:ENERGYBASIC} acordingly, to account for this sign change. Assuming that $\mu'>0$ on $I_1$, a  simple modification of the part 1. of Lemma~\ref{L:ENERGYBASIC} leads to
\begin{equation}\label{E:BOUNDTHREE}
    \|\partial_x\varphi_g\|_{L^2}^2\leq \eps^3\left(\mathcal{J}_{\text{nr}}^\ast+\mathcal{J}_{\text{trap}}^\ast+\mathcal{J}_{\text{ext}}^\ast\right), 
\end{equation}
where
\begin{align}
    \mathcal{J}_{\text{nr}}^\ast& =\sum_{\ell\in\mathbb{Z}_\ast^{>0}}\int_{E_{\ell}^{1,\de}}^{0} \frac{T(E)|\mu'(\k^2E)|}{\frac{T(E)}{\ell q}-1}|\widehat{\varphi}_g(\ell, E)|^2\diff E \notag\\
    & \ \ \ \ + 2 \sum_{\ell>0}\int_{E_{\ell}^{2,-\de}}^{\min\{E_\ell^{2,-\de},\frac{e_\ast}{\eps^2}\}} \frac{T(E) |\mu'(\k^2E)|}{|\frac{T(E)}{\ell q}-1|}\big( |\widehat{\varphi}_g(\ell, E)|^2+|\widehat{\varphi}_g(-\ell, E)|^2 \big) \diff E \notag\\
    & \ \ \ \ + 2 \sum_{\ell>0}\int_{(0,E_\ell^{2,\de}]\cap [\frac{e_\ast}{\eps^2},\infty)} \frac{T(E) |\mu'(\k^2E)|}{\frac{T(E)}{\ell q}-1}\big( |\widehat{\varphi}_g(\ell, E)|^2+|\widehat{\varphi}_g(-\ell, E)|^2 \big) \diff E,\label{E:JNRDEFAST}\\
    \mathcal{J}_{\text{trap}}^\ast & =\sum_{\ell>0}\int_{E_{\ell}^{1}}^{E_{\ell}^{1,\de}}\frac{T(E)|\mu'(\k^2E)|}{\frac{T(E)}{\ell q}-1}|\widehat{\varphi}_g(\ell, E)|^2\diff E, \label{E:JTRAPDEFAST}\\
    \mathcal{J}_{\text{ext}}^\ast &= 2 \sum_{\ell>0}\int_{E_\ell^2}^{\min\{E_{\ell}^{2,-\de},\frac{e_\ast}{\eps^2}\}} \frac{T(E) |\mu'(\k^2E)|}{|\frac{T(E)}{\ell q}-1|}\big( |\widehat{\varphi}_g(\ell, E)|^2+|\widehat{\varphi}_g(-\ell, E)|^2 \big) \diff E \notag\\
    & \ \ \ \ + 2 \sum_{\ell>0}\int_{[E_\ell^{2,\de},E_\ell^2]\cap [\frac{e_\ast}{\eps^2},\infty)} \frac{T(E) |\mu'(\k^2E)|}{\frac{T(E)}{\ell q}-1}\big( |\widehat{\varphi}_g(\ell, E)|^2+|\widehat{\varphi}_g(-\ell, E)|^2 \big) \diff E.
    \label{E:JEXTDEFAST}
\end{align}
We may now reprove the analogues of Lemmas~\ref{L:NRLEMMA},~\ref{L:TRAPLEMMA}, and~\ref{L:EXTLEMMA} using the same ideas to conclude the absence of embedded eigenvalues.



\section{Scaling analysis}\label{A:SCALING}

In this section we provide the algebraic scaling analysis which describes the scaling freedom in the problem.
We let
\begin{align}
\pa_t F +v \pa_x F + \pa_x \varphi_F \pa_v F & = 0,  \ \ (x,v)\in P\mathbb T^1\times\mathbb R, \label{E:VLASOV3}\\
-\pa_{xx}\varphi_F & = \tilde \rho_+(x) - \int F\diff v, \ \ x\in P\mathbb T^1.\label{E:POISSON3}
\end{align}
with some general period length $P>0$.
Upon setting 
\[
F(t,x,v) = \l^a f(\l^bt,\l x, \l^c v), \ \  \varphi_F(t,x)=\l^d\tilde\varphi(\l^bt,\l x), \ \  \tilde\rho_+(x)=\l^e \tilde\rho_+(\l x),
\]
we plug this ansatz into~\eqref{E:VLASOV3}--\eqref{E:POISSON3} and demand that the Vlasov equation be left invariant. A straightforward analysis
shows that $b =1-c$, $d=-2c$, and with $e=2-2c$, equations reduce to
\begin{align}
 \pa_tf +v \pa_x f + \pa_ x\varphi_f \pa_vf & = 0 ,  \ \ (x,v)\in P\l \mathbb T^1\times\mathbb R,\label{E:VLASOVRESCALED2}\\
-\pa_{xx}\varphi_f & = \rho_+ - \l^{a+c-2} \int f(x,v)\diff v, \ \ x\in P\l\mathbb T^1\label{E:POISSONRESCALED2}
\end{align}
where we have abused the notation by continuing to refer to the rescaled coordinates as $(t,x,v)$.
This is a 2-parameter freedom. 
If we assume $a+c\neq 2$ we may introduce the ``small" parameter
\[
\eps = \l^{a+c-2}
\]
to finally obtain the rescaled system~\eqref{E:VLASOV2}--\eqref{E:POISSON2}. At the level of steady states
the above analysis implies the following simple lemma.


\begin{lemma}
Let $a+c\neq 2$. Then any steady state of the VP-system~\eqref{E:VLASOV2}--\eqref{E:POISSON2} of the form 
\begin{align}
f_0 = \mu (\eps^2 E), \ \ E= \frac12v^2 - \varphi_0( x), \notag
\end{align}
supported on the periodic interval $x\in [0,P\eps^{\frac1{a+c-2}}]$ with the background ion density
\begin{align}
\rho_+( x) : = -\pa_{xx}\varphi_0 +  \eps \int \mu(\eps^2 E(x,v)) \diff v, \notag
\end{align}
gives rise to a steady state of the original VP-system~\eqref{E:VLASOV3}--\eqref{E:POISSON3} of the form
\begin{align}
F_0(x,v) &=  \eps^{\frac a{a+c-2}}f_0(\eps^{\frac1{a+c-2}}x, \eps^{\frac{c}{a+c-2}} v) =\eps^{\frac a{a+c-2}} \mu (\eps^{2+\frac{2c}{a+c-2}}e(x,v)), \label{E:S1}\\
\varphi_{F_0}(x) &=\eps^{-\frac{2c}{a+c-2}}\tilde\varphi_0(\eps^{\frac1{a+c-2}}x), \label{E:S2}
\end{align}
on the periodic interval 
$ x\in [0,P] $
with the ion background density defined as
\begin{align}
\tilde\rho_+(x) &=\eps^{\frac{2-2c}{a+c-2}} \rho_+(\eps^{\frac1{a+c-2}} x) \label{E:S3}
\end{align}
\end{lemma}



\begin{remark}
The 2-parameter family of rescalings~\eqref{E:S1}--\eqref{E:S2} maps any stationary solution of the rescaled system~\eqref{E:VLASOVRESCALED2}--\eqref{E:POISSONRESCALED2} into a stationary solution of the original VP-system, however with a different choice of the reference density and the spatial domain.
\end{remark}


\begin{remark}[$L^\infty$-critical scaling]
The  choice of scaling that keeps the $L^\infty$-norm of $\rho_+$ invariant is $c=1$ and the one that keeps the $L^\infty$ norm of $F_0$ invariant is $a=0$.  This then leads to the transformation
\begin{align}
F_0(x,v) = f_0(\frac x\eps,  \frac v\eps) = \mu (e(x,v)), \ \ 
\varphi_{F_0}(x) =\eps^{2}\varphi_0(\frac x\eps), \ \
\tilde\rho_+(x) =\rho_+(\frac{x}{\eps}), \notag
\end{align}
where $e(x,v)=\frac12 v^2 - \varphi_{F_0}(x)$. This corresponds to the scaling~\eqref{E:SCALINGFF}. 
\end{remark}

\begin{remark}[Stretching the periodic domain]
Other choices of parameters $(a,c)$ can be made while keeping the constraint $\eps=\l^{a+c-2}\ll1$. If $a+c-2<0$ the solution of the original problem~\eqref{E:VLASOV}--\eqref{E:POISSON} has a small spatial period for $0<\eps\ll1$, and if $a+c-2>0$ the period is large for $0<\eps\ll1$. 
Different choices of parameters $(a,c)$ will also accordingly rescale the profile of $\tilde\rho_+$.
\end{remark}

\end{document}